\documentclass[12pt]{article}
\usepackage{amsmath, amssymb}
\begin{document}
\bigbreak\centerline{\bf A General Non-Vanishing Theorem and an
Analytic Proof of }\centerline{\bf the Finite Generation of the
Canonical Ring}

\bigbreak
\centerline{Yum-Tong Siu\ %
\footnote{Partially supported by a grant from the National Science
Foundation.} }

\bigbreak\noindent{\sc Abstract.}  On August 5, 2005 in the
American Mathematical Society Summer Institute on Algebraic
Geometry in Seattle and later in several conferences I gave
lectures on my analytic proof of the finite generation of the
canonical ring for the case of general type. After my lectures
many people asked me for a copy of the slides which I used for my
lectures.  Since my slides were quite sketchy because of the time
limitation for the lectures, I promised to post later on a
preprint server my detailed notes from which my slides were
extracted. Here are my detailed notes giving the techniques and
the proof.

\bigbreak\centerline{\bf Table of Contents}

\medbreak\noindent\S1. Technique of Skoda on Ideal Generation

\smallbreak\noindent\S2. Reduction of Algebraic Geometric Problems
to $L^{2}$ Estimates for Stein Domains Spread Over $\mathbb{C}^{n}$

\smallbreak\noindent\S3. Stable Vanishing Orders and Their
Achievement by Finite Sums.

\smallbreak\noindent\S4. Reduction of Achievement of Vanishing Order
to Non-Vanishing Theorem on Hypersurface by Fujita Conjecture Type
Techniques

\smallbreak\noindent\S5. Diophantine Approximation of Kronecker

\smallbreak\noindent\S6. A General Non-Vanishing Theorem

\smallbreak\noindent\S7. Holomorphic Family of Artinian Subschemes
and Achievement of Stable Vanishing Orders for the Case of Higher
Codimension

\smallbreak\noindent\S8. Remark on the Approach of Extension Using
Techniques of the Invariance of Plurigenera

\smallbreak\noindent\S9. Remark on Positive Lower Bound of Curvature
Current

\bigbreak On August 5, 2005 in the American Mathematical Society
Summer Institute on Algebraic Geometry in Seattle I gave a lecture
in which I first presented my analytic method of proving the finite
generation of the canonical ring for the case of general type. Later
in several conferences (the Birthday Conference for Skoda in Paris
on September 12, 2005; the Memorial Conference for Vitushkin in
Moscow on September 26, 2005; the Birthday Conference for Bogomolov
in Miami on December 18, 2005; the Birthday Conference for Toledo in
Salt Lake City on March 24, 2006; the Retirement Conference for
Kiselman in Uppsala on May 15, 2006; the Birthday Conference of Lu
Qikeng in Beijing on June 6, 2006; the Trento Conference on CR
Geometry and PDE on September 6, 2006; the Second Chinese-German
Conference on Complex Analysis in Shanghai on September 12, 2006) I
lectured on the same analytic proof of the finite generation of the
canonical ring for the case of general type. After my lectures many
people asked me for a copy of the slides which I used for my
lectures. Since my slides were quite sketchy because of the time
limitation for the lectures, I promised to post later on a preprint
server my detailed notes from which my slides were extracted. Here
are my detailed notes giving the techniques and the proof.

\medbreak These notes are selected and organized from the notes of
roughly one year old which I wrote for myself over a period of
several years as memoranda while I worked on the problem of the
finite generation of the canonical ring.  Unlike a formal preprint
which follows the traditional order of history, definitions,
lemmas, propositions, and theorems, these notes start out with key
ideas and techniques as the backbone and then flesh out with more
and more details and explanations on how to deal with the
difficulties which arise in the course of the implementation of
the key ideas and techniques until the complete solution is
reached.  I made some selections when there are several ways of
handling a difficulty and I unified the notations and terminology
and put in the numbering for sections, paragraphs, definitions,
lemmas, propositions, theorems, remarks, equations, {\it et
cetera}, but the presentation retains essentially the order and
the substance of the material in my original notes.  This style of
presentation actually makes the geometric ideas and the techniques
for the proof more transparent. I hope that the people interested
in the details of the techniques and the proof of the finite
generation of the canonical ring presented in my lectures in the
several conferences since August 2005 will find these notes easier
to read and understand than a formal preprint.

\medbreak The result on the finite generation of the canonical
ring for the case of general type which these notes give an
analytic proof for is the following.

\bigbreak\noindent{\it Main Theorem.}  Let $X$ be a compact
complex algebraic manifold of complex dimension $n$ which is of
general type in the sense that there exist a positive integer
$m_0$ and a positive number $c$ such that ${\rm dim}_{\mathbb
C}\Gamma\left(X, mK_X\right)\geq cm^n$ for $m\geq m_0$, where
$K_X$ is the canonical line bundle of $X$.  Then the canonical
ring $\bigoplus_{m=1}^\infty\Gamma\left(X, mK_X\right)$ is
finitely generated.

\bigbreak An important component in the analytic proof of the
finite generation of the canonical ring for the case of general
type is the general non-vanishing theorem in the title of these
notes which is stated in (6.2) with its proof given in \S6.

\medbreak Since the problem of the finite generation of the
canonical ring is a well-known problem in algebraic geometry, the
history of the problem is not repeated in these notes.

\medbreak It was brought to my attention that on October 5, 2006
Caucher Birkar, Paolo Cascini, Christopher D. Hacon, James
McKernan posted a preprint on the ``existence of minimal models
for varieties of log general type'' on the ``arXiv.org'' server.

\bigbreak I explain here the organization of these notes.  The key
ingredient in the analytic proof of the finite generation of the
canonical ring is the theorem of Skoda on ideal generation.  It is
presented in \S1.  With $L^2$ estimates problems in algebraic
geometry involving holomorphic sections of holomorphic line bundles
over compact complex algebraic manifolds can be reduced to problems
for Stein domains spread over ${\mathbb C}^n$.  This is presented in
\S2, with explanations on what such a reduction means in the problem
of the finite generation of the canonical ring.  In \S3 an infinite
sum $\Phi$ of the absolute-value squares of fractional powers of
pluricanonical sections is introduced and the problem of the finite
generation of the canonical ring is shown to be equivalent to the
precise achievement of stable vanishing orders in the sense that the
vanishing orders of the infinite sum $\Phi$ can be achieved by some
of its finite partial sums.  The proof is done by descending
induction on the dimension of the subvariety $V$ where the stable
vanishing order is not yet known to be precisely achieved. In \S4 by
techniques for Fujita conjecture type problems, the initial
induction step where $V$ is a hypersurface is reduced to a general
non-vanishing theorem.

\medbreak In \S5 we group together results derived from Kronecker's
theorem on diophantine approximation which will be needed later for
our general non-vanishing theorem.  In \S6 the general non-vanishing
theorem is presented which gives the existence of nonzero
holomorphic sections belonging locally to an appropriate multiplier
ideal sheaf, under some positive lower bound condition for the
curvature current of the line bundles involved. This general
non-vanishing theorem rules out the possibility of an infinite
number of components for the embedded stable base point set, which
is the major obstacle in getting the finite generation of the
canonical ring. Essential for the proof of the general non-vanishing
theorem is Shokurov's technique of comparing two applications of the
theorem of Riemann-Roch, one to a line bundle and another to its
twisting by a flat line bundle.

\medbreak In \S7 the method of continuous variation of an Artinian
subscheme without jump is introduced in order to prove the precise
achievement of the stable vanishing order at a generic point of a
subvariety of higher codimension.  The vanishing order of $\Phi$
across a hypersurface is a single number, but in the case of a
subvariety $V$ of higher codimension this r\^ole is played by
Artinian subschemes transversal to the subvariety $V$.  Each finite
partial sum of $\Phi$ provides one Artinian subscheme which varies
continuous along the subvariety $V$ without jump except at a
subvariety $E$ of codimension $\geq 1$ in $V$.  At points of $V$
outside the countable union of this kind of subvarieties $E$ the
stable vanishing order is precisely achieved, making it possible to
go to the next step in the induction process. The proof of the
finite generation of the canonical ring is completed in \S7. Also in
\S7 it is explained why the proof of precisely achieving the stable
vanishing order at a generic point of a subvariety of higher
codimension cannot simply be reduced to the hypersurface case by
blow-up and how the continuous variation of an Artinian subscheme
without jump handles the problem.

\medbreak The development of extension techniques for the problem
of the deformational invariance of plurigenera was originally
intended for application to the problem of the finite generation
of the canonical ring.  In \S8 the approach by such extension
techniques to the problem of the finite generation of the
canonical ring is compared to the analytic proof presented in
these notes. Difficulties with the approach by extension
techniques are analyzed.  Finally \S9 gives some remarks
concerning the condition on positive lower bounds for curvature
currents, including a remark about the problem of the finite
generation of the canonical ring without the general type
condition and the difficulty of artificially adding an ample
twisting first with the taking of root-limits at the end to get
rid of its contribution.

\bigbreak\noindent{\it Notations.} ${\mathbb N}$ is the set of all
positive integers.  ${\mathbb Z}$ is the set of all integers.
${\mathbb Q}$ is the set of all rational numbers.  ${\mathbb R}$
is the set of all real numbers.  ${\mathbb C}$ is the set of all
complex numbers.  For a subvariety $V$ we use ${\mathcal I}_V$ to
denote the coherent ideal sheaf of all germs of holomorphic
functions vanishing on $V$.  The structure sheaf of a complex
space $X$ is denoted by ${\mathcal O}_X$.  The maximum ideal of
$X$ at a point $P$ is denoted by ${\mathfrak m}_{X,P}$.  A
multi-valued holomorphic section $s$ of a ${\mathbb
Q}$-line-bundle $E$ means that $s^N$ is a holomorphic section of
the holomorphic line bundle $NE$ for some positive integer $N$.
For a divisor $Y$ we denote by $s_Y$ the canonical section of the
line bundle defined by $Y$.  When $Y$ is a ${\mathbb Q}$-divisor,
the canonical section $s_Y$ is a multi-valued holomorphic section
of the ${\mathbb Q}$-line-bundle defined by $Y$. The space of all
sections of a holomorphic bundle or a sheaf $E$ over $X$ is
denoted by $\Gamma\left(X, E\right)$.  The term ``generic'' is
also used in the sense of avoiding some countable union of
subvarieties of codimension $\geq 1$ (or even some countable union
of locally defined subvarieties of codimension $\geq 1$). The
round-down of a real number $u$ is denoted by $\left\lfloor
u\right\rfloor$ which is the largest integer $\leq u$.  The
round-up of a real number $u$ is denoted by $\left\lceil
u\right\rceil$ which is the smallest integer $\geq u$.

\bigbreak

\bigbreak\noindent{\bf \S1. Technique of Skoda on Ideal Generation}

\bigbreak The key ingredient in the analytic proof of the finite
generation of the canonical ring is the following result of Skoda on
ideal generation [Skoda 1972].  (Skoda's original statement is for a
Stein domain in ${\mathbb C}^n$, but for its application to
algebraic geometry we need the version of a Stein domain spread over
${\mathbb C}^n$.)

\bigbreak\noindent(1.1) {\it Theorem (Skoda on Ideal Generation).}
Let $\Omega $ be a domain spread over $\mathbb{C}^{n}$ which is
Stein. Let $\psi$ be a plurisubharmonic function on $\Omega $,
$g_{1},\ldots ,g_{p}$ be holomorphic functions on $\Omega $, $\alpha
>1$, $q=\min \left(n,p-1\right)$, and $f$ be a holomorphic function on $\Omega
$. Assume that
$$\int_{\Omega }\frac{\left\vert f\right|^{2}e^{-\psi }}{\left(
\sum_{j=1}^p\left| g_{j}\right| ^{2}\right)^{\alpha q+1}}<\infty.$$
Then there exist holomorphic functions $h_{1},\ldots ,h_{p}$ on
$\Omega$ with $f=\sum_{j=1}^{p} h_jg_j$ on $\Omega$ such that
$$\int_{\Omega }\frac{\left\vert h_{k}\right|^{2}e^{-\psi
}}{\left(\sum_{j=1}^p\left|g_{j}\right|^2\right)^{\alpha q}}\leq
\frac{\alpha }{\alpha -1}\int_{\Omega
}\frac{\left|f\right|^2e^{-\psi }}{\left( \sum_{j=1}^p\left|
g_{j}\right|^2\right)^{\alpha q+1}}$$ for $1\leq k\leq p$.

\bigbreak\noindent(1.2) {\it Remark on the Condition of Skoda's
Theorem on Ideal Generation.}  The condition of Skoda's theorem on
ideal generation comes from transplanting the obvious necessary
supremum condition to an $L^2$ condition.  It is clear that a
necessary condition to express $f$ in terms of $g_1,\cdots,g_p$ as
$f=\sum_{j=1}^p h_j g_j$ is that $\left|f\right|\leq
C\sum_{j=1}^p\left|g_j\right|$ on any compact subset of $\Omega$
with $C$ depending on the compact subset. However, we have to use
Hilbert space arguments instead of Banach space arguments with the
supremum norm. Thus we need formulation in $L^2$ bounds such as
$$
\int\frac{\left|f\right|^2}{\sum_{j=1}^p\left|g_j\right|^2}<\infty.
$$
The translation of the condition from supremum norm to $L^2$ norm
calls for modification in the formulation, because in the case of
$f=1$ and
$\left\{g_1,\cdots,g_p\right\}=\left\{z_1,\cdots,z_n\right\}$, the
integral is finite by polar coordinate argument for $n\geq 2$ and
yet $f$ cannot be so expressed.  Skoda's formulation modifies the
integral bound to
$$
\int\frac{\left|f\right|^2}{\left(\sum_{j=1}^p\left|g_j\right|^2\right)^{\alpha
q+1}}<\infty
$$
with $\alpha>1$ and $q=\min\left(n,p-1\right)$ for polar coordinate
reasons.

\medbreak We have to pay a price for translating the formulation of
the assumption to the Hilbert space context in that in the
denominator there is a gap between the exponent for the sufficiency
and the exponent for the necessity, illustrated in the two extreme
cases of (i) $p=n$ and $g_1=z_1,\cdots,g_n=z_n$ and (ii) $p=1$ and
$g_1=z_1$.

\bigbreak

\bigbreak\noindent {\bf \S2. Reduction of Algebraic Geometric
Problems to $L^{2}$ Estimates for Stein Domains Spread Over
$\mathbb{C}^{n}$}.

\bigbreak Since $L^2$ holomorphic functions can be extended across a
complex hypersurface, by representing an Zariski open subset of a
compact complex algebraic manifold as a Stein domain spread over
${\mathbb C}^n$ and using a meromorphic section of a holomorphic
line bundle, we can reduce a problem in algebraic geometry involving
line bundles to an analytic problem for $L^2$ estimates on a Stein
domain spread over ${\mathbb C}^n$ as follows.

\bigbreak\noindent(2.1) {\it Description of Reduction.} Let $X$ be
an $n$-dimensional complex manifold inside $\mathbb{P}_{N}$.  We
regard $\mathbb{P}_{N-n-1}\subset \mathbb{P}_{N}-X$ as a source of
light. Take any $\mathbb{P}_{n}\subset \mathbb{P}_{N}$. Fix $x\in
X$, define $\pi(x)$ as the only point in ${\rm span
}\left(x,\mathbb{P}_{N-n-1}\right)\cap\mathbb{P}_n$ to make
$\pi:X\to{\mathbb P}_n$ a branched cover, where ${\rm span
}\left(x,\mathbb{P}_{N-n-1}\right)$ means the projective linear
subspace of dimension $N-n$ in ${\mathbb P}_N$ which contains both
$x$ and $\mathbb{P}_{N-n-1}$.

\medbreak Let $L$ be a holomorphic line bundle over $X$ and $s$ be a
global meromorphic section with pole-set $A$ and zero-set $B$. Take
some hypersurface $Z$ inside $\mathbb{P}_{n}$ containing ${\mathbb
P}_n-{\mathbb C}^n$ and $\pi\left(A\cup B\right)$ such that $\pi
:X-\pi ^{-1}(Z)\rightarrow \mathbb{P}_{n}-Z$ is a local
biholomorphism.  Then $X-\pi^{-1}(Z)$ is a Stein domain spread over
${\mathbb C}^n$.

\medbreak Take a metric $e^{-\varphi}$ of $L$ with $\varphi$ locally
bounded from above.  For any open subset $\Omega$ of $X$ and any
holomorphic function $f$ on $\Omega-\pi^{-1}(Z)$ with
$$\int_{\Omega-\pi^{-1}(Z)}\left|f\right|^2e^{2\log|s|-\varphi}<\infty,$$
the section $fs$ of $L$ can be extended to a holomorphic section of
$L$ over $\Omega$.  So instead of dealing with $X$ and $L$ to find
elements of $\Gamma\left(X,L\right)$, we can deal with holomorphic
functions on a Stein domain spread over ${\mathbb C}^n$ which
satisfy certain $L^2$ estimates.  This method of reduction
translates Skoda's theorem on ideal generation for Stein domains
spread over ${\mathbb C}^n$ to results concerning holomorphic
sections of line bundles over compact complex algebraic manifolds.
We are going to use it to reduce the problem of finite generation of
the canonical ring to a problem for the precise achievement of
stable vanishing orders which will be explained in details later.

\bigbreak\noindent(2.2) {\it Reformation of the Ideal Generation of
Skoda for Compact Algebraic Manifolds and Holomorphic Line Bundles.}
We now apply, to the reformulation of the theorem of Skoda on ideal
generation, the above procedure of reducing algebraic geometric
problems to $L^{2}$ estimates for Stein domains spread over
$\mathbb{C}^{n}$.  First we introduce the following definition and
notation for multiplier ideal sheaves.

\bigbreak\noindent(2.3) {\it Definition of Multiplier Ideal
Sheaves.} For a local plurisubharmonic function $\varphi$ on an open
subset of ${\mathbb C}^n$, the multiplier ideal sheaf ${\mathcal
I}_\varphi$ is the sheaf of germs of holomorphic functions $f$ such
that $\left|f\right|^2e^{-\varphi}$ is locally integrable.

\bigbreak Using the notation of multiplier ideal sheaves, we can
formulate the theorem of Skoda on ideal generation for compact
algebraic manifolds and holomorphic line bundles as follows.

\bigbreak\noindent(2.4) {\it Theorem.} Let $X$ be a compact complex
algebraic manifold of complex dimension $n$, $L$ be a holomorphic
line bundle over $X$, and $E$ be a holomorphic line bundle on $X$
with metric $e^{-\psi }$ such that $\psi $ is plurisubharmonic. Let
$k\geq 1$ be an integer, $G_{1},\ldots ,G_{p}\in \Gamma (X,L)$, and
$\left|G\right|^{2}=\underset{j=1}{\overset{p}{\sum
}}\left|G_j\right|^{2}$. Let
$\mathcal{I=I}_{(n+k+1)\log\left|G\right|^{2}+\psi}$ and
$\mathcal{J=I}_{(n+k)\log\left|G\right|^2+\psi}$. Then
$$\displaylines{\qquad\qquad\Gamma\left(X,\mathcal{I}\otimes
\left((n+k+1)L+E+K_{X}\right)\right)\hfill\cr\hfill
=\underset{j=1}{\overset{p}{%
\sum }}G_{j}\,\Gamma \left(X,\mathcal{J}\otimes\left(
(n+k)L+E+K_{X}\right)\right).\qquad\qquad\cr}$$

\medskip\noindent {\it Proof}. Take $F\in\Gamma\left(X,\mathcal{I}\otimes
\left((n+k+1)L+E+K_{X}\right)\right)$.  Let $S$ be a meromorphic
section of $E$. Take a branched cover map $\pi:X\rightarrow{\bf
P}_n$.  Let $Z_0$ be a hypersurface in ${\bf P}_n$ which contains
the infinity hyperplane of ${\bf P}_n$ and the branching locus of
$\pi$ in ${\bf P}_n$ such that $Z:=\pi^{-1}(Z_0)$ contains the
divisor of $G_1$ and both the pole-set and zero-set of $S$ . Let
$\Omega=X-Z$. Let $g_j=\frac{G_j}{G_1}$ ($1\leq j\leq p$) and
$\left|g\right|^2=\sum_{j=1}^p\left|g_j\right|^2$.  Define $f$ by
$$
\frac{F}{G_1^{n+k+1}S}=fdz_1\wedge\cdots\wedge dz_n,
$$
where $z_1,\cdots,z_n$ are the affine coordinates of ${\mathbb
C}^n$.  Use $\alpha=\frac{n+k}{n}$.  Let $\psi=\varphi-\log|S|^2$.
It follows from $F\in{\cal I}_{(n+k+1)\log|G|^2+\varphi}$ locally
that
$$
\int_X\frac{|F|^2}{|G|^{2(n+k+1)}}e^{-\varphi}<\infty,
$$
which implies that
$$
\int_\Omega\frac{|f|^2}{|g|^{2(n+k+1)}}e^{-\psi}=
\int_\Omega\frac{\left|\frac{F}{G_1^{n+k+1}S}\right|^2}{\left|\frac{G}{
G_1}\right|^{2(n+k+1)}}e^{-\psi}=
\int_\Omega\frac{|F|^2}{|G|^{2(n+k+1)}}e^{-\varphi}<\infty.
$$
By Skoda's theorem on ideal generation (1.1) with $q=n$ (which we
assume by adding some $F_{p+1}\equiv\cdots\equiv F_{n+1}\equiv 0$ if
$p<n+1$) so that $2\alpha q+2=2\cdot\frac{n+k}{n}\cdot
n+2=2(n+k+1)$, we obtain holomorphic functions $h_1,\cdots,h_p$ on
$\Omega$ such that $f=\sum_{j=1}^p h_jg_j$ and
$$
\sum_{j=1}^p\int_\Omega\frac{|h_j|^2}{|g|^{2(n+k)}}e^{-\psi}<\infty.
$$
Define
$$
H_j=G_1^{n+k}h_j S dz_1\wedge\cdots\wedge dz_n.
$$
Then $F=\sum_{j=1}^p H_j G_j$ and
$$
\int_\Omega\frac{|H_j|}{|G|^{2(n+k)}}e^{-\varphi}=
\int_\Omega\frac{|h_j|}{|g|^{2(n+k)}}e^{-\psi}<\infty
$$
so that $H_j$ can be extended to an element of
$\Gamma(X,(n+k)L+E+K_X)$. Q.E.D.

\medskip
As an illustration, we give the following trivial immediate
consequence of Theorem (2.4).  It is on the effective finite
generation for the simple case of a numerically effective canonical
line bundle, which is interesting in its own right, but not useful
for our purpose.

\bigskip \noindent(2.5) {\it Corollary.}  Let $F$ be a holomorphic
line bundle over a compact complex algebraic manifold $X$ of complex
dimension $n$.  Let $a>1$ and $b\geq 0$ be integers such that $aF$
and $bF-K_X$ are globally free over $X$.  Then the ring
$\bigoplus_{m=1}^\infty\Gamma(X,mF)$ is generated by
$\bigoplus_{m=1}^{(n+2)a+b-1}\Gamma(X,mF)$.

\medskip
\noindent {\it Proof.} For $0\leq\ell<a$ let $E_\ell=(b+\ell)F-K_X$
and $L=aF$.  Let $G_1,\cdots,G_p$ be a basis of
$\Gamma(X,L)=\Gamma(X,aF)$.  Let $H_1,\cdots,H_q$ be a basis of
$\Gamma(X,bF-K_X)$.  We give $E_\ell$ the metric
$$
\frac{1}{(\sum_{j=1}^p|G_j|^{\frac{2\ell}{
a}})(\sum_{j=1}^q|H_j|^2)}.
$$
Since both ${\cal I}$ and ${\cal J}$ from (2.4) are unit ideal
sheaves, it follows that
$$
\Gamma(X,(n+k+1)L+E_\ell+K_X) =\sum_{j=1}^p
G_j\Gamma(X,(n+k)L+E+K_X)
$$
for $k\geq 1$ and $0\leq\ell<a$, which means that
$$
\Gamma(X,((n+k+1)a+\ell+b)F) =\sum_{j=1}^p
G_j\Gamma(X,((n+k)a+\ell+b)F)
$$
for $k\geq 1$ and $0\leq\ell<a$.  Thus
$\bigoplus_{m=1}^{(n+2)a+b-1}\Gamma(X,mF)$ generates the ring
$\bigoplus_{m=1}^\infty\Gamma(X,mF)$.  Q.E.D.

\bigbreak

\bigbreak\noindent{\bf \S3. Stable Vanishing Orders and Their
Achievement by Finite Sums}.

\bigbreak Let $X$ be a compact complex algebraic manifold of complex
dimension $n$ which is of general type.  Let
$$
\Phi=\sum_{m=1}^\infty\varepsilon_m\sum_{j=1}^{q_m}\left|s^{(m)}_j\right|^{\frac{2}{m}},
$$
where
$$
s^{(m)}_1,\cdots,s^{(m)}_{q_m}\in\Gamma\left(X,\,mK_X\right)
$$
form a basis over ${\mathbb C}$ and $\varepsilon_m>0$ approach $0$
so fast as $m\to\infty$ that locally the infinite series which
defines $\Phi$ converges uniformly.  Note that the expression
$\frac{1}{\Phi}$ defines a metric for $K_X$.  In these notes we will
reserve the symbol $\Phi$ only for this kind of infinite sums.

\medbreak Likewise, more generally, for a big holomorphic line
bundle over $X$, we can introduce
$$
\Psi=\sum_{m=1}^\infty\delta_m\sum_{j=1}^{r_m}\left|g^{(m)}_j\right|^{\frac{2}{m}},
$$
where
$$
g^{(m)}_1,\cdots,g^{(m)}_{r_m}\in\Gamma\left(X,\,m L\right)
$$
form a basis over ${\mathbb C}$ and $\delta_m>0$ approach $0$ so
fast as $m\to\infty$ that locally the infinite series which defines
$\Psi$ converges uniformly.  In this case $\frac{1}{\Psi}$ defines a
metric for $L$.

\medbreak What plays the most important r\^ole in the finite
generation of the canonical ring is the vanishing orders of the
local function $\Phi$.  Since $\Phi$ is an infinite sum of the
absolute-value squares of fractional powers of holomorphic
functions, the easiest way to define the vanishing orders of $\Phi$
rigorously is to use Lelong numbers.

\bigbreak\noindent(3.1) {\it Definition of Lelong Numbers.}  A
$(1,1)$-current $\Theta$ on an open subset $G$ of ${\mathbb C}^n$ is
said to be {\it positive} if for any smooth $(1,0)$-forms
$\sigma_1,\cdots,\sigma_{n-1}$ on $G$ with compact support the
$(n,n)$-current
$$
\Theta\wedge\left(\sqrt{-1}\,\sigma_1\wedge
\overline{\sigma_1}\right)\wedge\cdots\wedge\left(\sqrt{-1}\,\sigma_{n-1}\wedge
\overline{\sigma_{n-1}}\right)
$$
on $G$ is a nonnegative measure.  Note that though by convention the
adjective ``positive'' is used in this definition, it actually means
just ``nonnegative'' and not ``strictly positive.''

\medbreak For a complex hypersurface $Y$ in $G$,
$\theta\mapsto\int_{{\rm Reg\,}Y}\theta$ for any smooth
$(n-1,n-1)$-form $\theta$ on $G$ with compact support defines a
closed positive $(1,1)$-current which we denote by $\left[Y\right]$,
where ${\rm Reg\,}Y$ is the regular part of $Y$ consisting of all
nonsingular points of $Y$.

\medbreak Every closed positive $(1,1)$-current $\Theta$ can be
locally written as $\sqrt{-1}\partial\bar\partial\varphi$ for some
plurisubharmonic function $\varphi$ (sometimes known as a {\it
plurisubharmonic potential}) and, conversely,
$\sqrt{-1}\partial\bar\partial\varphi$ is a closed positive
$(1,1)$-current for any plurisubharmonic function $\varphi$.

\medbreak For a closed positive $(1,1)$-current $\Theta$ on some
open subset $G$ of ${\mathbb C}^n$, the {\it Lelong number} of
$\Theta$ at a point $P_0$ is the limit of
$$
\frac{\ \ \int_{B_n(P_0,r)}{\rm trace\,}\Theta\ \ }{{\rm
Vol\,}\left(B_{n-1}(0,r)\right)}
$$
as $r\to 0$, where $B_m\left(Q,r\right)$ is the open ball in
${\mathbb C}^m$ of radius $r$ centered at $Q$ and ${\rm
Vol\,}\left(B_m\left(Q,r\right)\right)$ is its volume and ${\rm
trace\,}\Theta$ is
$$\Theta\wedge\frac{1}{(n-1)!}\left(\sum_{j=1}^n\frac{\sqrt{-1}}{2}\,dz_j\wedge
d\overline{z_j}\right)^{n-1}$$ with $z_1,\cdots,z_n$ being the
coordinates of ${\mathbb C}^n$.

\medbreak For a complex hypersurface $Y$ in $G$ the Lelong number of
$\left[Y\right]$ at a point $P_0$ of $G$ is the multiplicity of $Y$
at $P_0$.

\medbreak For any $c>0$ and any closed positive $(1,1)$-current
$\Theta$, we denote by $E_c\left(\Theta\right)$ the set of points
where the Lelong number of $\Theta$ is no less than $c$.  The set
$E_c\left(\Theta\right)$ is always a subvariety.  An irreducible
branch $E$ of the subvariety $E_c\left(\Theta\right)$ is called an
{\it irreducible Lelong set} of $\Theta$. By the {\it generic Lelong
number} of an irreducible Lelong set $E$ for $\Theta$ is meant the
Lelong number of $\Theta$ at any generic point $P$ of $E$ and is
independent of the choice of $P$.  For any subvariety $V$ (not
necessarily an irreducible Lelong set), the generic Lelong number of
$\Theta$ at $V$ means the Lelong number of $\Theta$ at a generic
point of $V$.

\medbreak For a metric $e^{-\varphi}$ we use the notation
$$\Theta_\varphi=\frac{\sqrt{-1}}{2\pi}\partial\bar\partial\varphi$$
to denote the curvature current of $e^{-\varphi}$.  In the case of
the metric $\frac{1}{\Phi}$ for the canonical line bundle $K_X$ of
$X$ its curvature current is given by
$$\Theta_{\log\Phi}:=\frac{\sqrt{-1}}{2\pi}\partial\bar\partial\log\Phi.$$

\medbreak\noindent(3.1.1) {\it Additional Factor of Two in Vanishing
Order of $\Phi$ Compared to Lelong Number.} We use the Lelong number
of $2\Theta_{\log\Phi}$ at a point $P$ of $X$ to define the
vanishing order of $\Phi$ at $P$.  In this definition of the
vanishing order of $\Phi$, a factor of $2$ is introduced in
$2\Theta_{\log\Phi}$, because, in the definition of $\Phi$ as an
infinite series, absolute-value squares are taken of the fractional
powers of holomorphic functions in the individual terms instead of
just the absolute value.  Note that when the complex dimension $n$
of $X$ is $>1$, the Lelong number of
$2\,\frac{\sqrt{-1}}{2\pi}\partial\bar\partial\log\Phi$ at a point
$P$ only measures the vanishing order of the restriction of $\Phi$
to a generic local complex curve in $X$ passing through $P$.

\medbreak Because of the confusion which this additional factor of
$2$ may cause, when we assign a number to the vanishing order, we
always use the one computed as the Lelong number of the closed
positive $(1,1)$-current
$\frac{\sqrt{-1}}{2\pi}\partial\bar\partial\log\Phi$.  When we refer
to the vanishing order of $\Phi$, it will only be for comparison
with the vanishing order of its partial sum
$$
\Phi_{m_0}:=\sum_{m=1}^{m_0}\varepsilon_m\sum_{j=1}^{q_m}\left|s^{(m)}_j\right|^{\frac{2}{m}}
$$
so that the additional factor of $2$ would not make any difference.

\bigbreak\noindent(3.2) {\it Failure of the Noetherian Argument for
Irreducible Lelong Sets Defined by $\Phi$.}  Though $\Phi$ is
defined by an infinite sum of the absolute-value squares of
fractional powers of holomorphic functions, the use of fractional
powers make the Noetherian argument inapplicable and we cannot
conclude from the Noetherian argument that the number of irreducible
Lelong sets of $\Theta_{\log\Phi}$ is finite.

\medbreak To understand how this can happen, we illustrate the
situation with the following simple example.  On ${\mathbb C}^2$
with coordinates $(z,w)$, the vanishing order of $z$ is $1$ at $z=0$
(measured by the vanishing order of its restriction to a generic
complex line through $0$).  The vanishing order of
$\left|z\right|+\left|w-\frac{1}{\ell}\right|^{\frac{1}{\ell}}$ at
$(z,w)=\left(0,\frac{1}{\ell}\right)$ is $\frac{1}{\ell}$ (as
defined by using the Lelong number of the closed positive
$(1,1)$-current defined by it). The collection of all the
irreducible Lelong sets of the closed positive $(1,1)$-current
defined by
$$\sum_{m=1}^\infty\varepsilon_m\left|z\right|^{\frac{1}{m}}
\prod_{\ell=1}^m\left(\left|z\right|+\left|w-\frac{1}{\ell}\right|^{\frac{1}{\ell}}\right)\leqno{(3.2.1)}
$$
(for positive numbers $\varepsilon_m$ decreasing fast enough)
consists of $\left\{\left(0,\frac{1}{m}\right)\right\}$ for
$m\in{\mathbb N}$ with the Lelong number $\frac{1}{m}$ at the point
$\left\{\left(0,\frac{1}{m}\right)\right\}$.  When we use the finite
partial sum
$$\sum_{m=1}^N\varepsilon_m\left|z\right|^{\frac{1}{m}}
\prod_{\ell=1}^m\left(\left|z\right|+\left|w-\frac{1}{\ell}\right|^{\frac{1}{\ell}}\right),\leqno{(3.2.2)}
$$
instead of the infinite sum $(3.2.1)$, the collection of all the
irreducible Lelong sets of the closed positive $(1,1)$-current
defined by $(3.2.2)$ consists of
$\left\{\left(0,\frac{1}{m}\right)\right\}$ for $1\leq m\leq N-1$
plus the set $\left\{z=0\right\}$, where the Lelong number at the
point $\left(0,\frac{1}{m}\right)$ is $\frac{1}{m}$ and the Lelong
number at a generic point of $\left\{z=0\right\}$ is $\frac{1}{N}$.

\medbreak The Lelong set $\left\{\left(0,\frac{1}{m}\right)\right\}$
with Lelong number $\frac{1}{m}$ for $m\geq N$ does not occur for
the closed positive $(1,1)$-current defined by $(3.2.2)$, because it
is hidden inside and incorporated into the irreducible Lelong set
$\left\{z=0\right\}$ where the generic Lelong number is
$\frac{1}{N}$.  The Noetherian argument can be applied to the case
of the finite sum, ruling out an infinite number of irreducible
Lelong sets (because we can raise each individual term to some
finite common power to get rid of the fractional powers of all the
individual terms), but cannot be applied to the case of the infinite
sum. The situation is like the following picture.

\medbreak For an infinite number of pebbles arranged at different
heights on a slope immersed in water, when the water recedes slowly,
gradually more and more (but only a finite number of) pebbles are
visible at any given time until the water completely drains out to
make all the infinite number of pebbles visible.

\medbreak An important tool for our analytic proof of the finite
generation of the canonical ring for the case of general type is the
following decomposition theorem for closed positive
$(1,1)$-currents.  (See [Siu 1974] and [Kiselman 1979].)

\bigbreak\noindent(3.3) {\it Theorem (Decomposition of Closed
Positive (1,1)-Currents).}  Let $\Theta$ be a closed positive
$(1,1)$-current on a complex manifold $X$.  Then $\Theta$ admits a
unique decomposition of the following form
$$ \Theta=\sum_{j=1}^J\gamma_j\left[V_j\right]+R,
$$
where $\gamma_j>0$, $J\in{\mathbb N}\cup\left\{0,\infty\right\}$,
$V_j$ is a complex hypersurface and the Lelong number of the
remainder $R$ is zero outside a countable union of subvarieties of
codimension $\geq 2$ in $X$.

\bigbreak\noindent(3.4) {\it Pullbacks of Closed Positive
(1,1)-Currents.}  Though in general currents can only be pushed
forward and cannot be pulled back, yet for the case of a closed
positive $(1,1)$-current $\Theta$ on a complex manifold $G$ it is
possible to pull back by a surjective holomorphic map $\pi: D\to G$
from $G$ to a manifold $D$, because we can locally write
$\Theta=\sqrt{-1}\partial\bar\partial\varphi$ for some local
plurisubharmonic function on $G$ and pull back $\varphi$ to
$\pi^*\varphi$ and then form
$\sqrt{-1}\partial\bar\partial\pi^*\varphi$ as the pullback of
$\Theta$ to $D$.

\medbreak The assumption of the surjectivity of $\pi$ is just to
make sure that the pullback $\pi^*\varphi$ is not identically
$-\infty$ on $D$.  Pullbacks can also be defined for the case when
the holomorphic map $\pi: D\to G$ is not surjective as long as
locally the image is not contained in the ($-\infty$)-set of the
local plurisubharmonic function $\varphi$.  If a submanifold $V$ of
$G$ is not entirely contained in the ($-\infty$)-set of the
plurisubharmonic potential $\varphi$, we can define the restriction
of the closed positive $(1,1)$-current
$\Theta=\sqrt{-1}\partial\bar\partial\varphi$ to a submanifold $V$
of $G$ (or even a subvariety $V$ of $G$) as
$\sqrt{-1}\partial\bar\partial\left(\varphi\big|_V\right)$.

\bigbreak\noindent(3.5) {\it Closed Positive (1,1)-Currents on
Complex Spaces and Their Decomposition.}  We can also define a
closed positive $(1,1)$-current $\Theta$ on a reduced complex space
$Y$ with singularities by defining it as
$\sqrt{-1}\partial\bar\partial\varphi$ for some local
plurisubharmonic function $\varphi$ on $Y$ (in the sense that
$\varphi$ can be extended to some plurisubharmonic function on some
complex manifold in which locally $Y$ is a subvariety).  We can also
consider the decomposition of the closed positive $(1,1)$-current
$\Theta$ on $Y$ by pulling it back to a complex manifold $\tilde Y$
which is a desingularization of $Y$ and then push forward to $Y$ the
decomposition on $\tilde Y$ of the pullback of $\Theta$.

\bigbreak\noindent(3.6) {\it Finite Generation from Precise
Achievement of Stable Vanishing Order.}  Recall that a vanishing
order of $\Phi$ at a point we mean the Lelong number of
$$2\Theta_{\log\Phi}:=2\,\frac{\sqrt{-1}}{2\pi}\partial\bar\partial\log\Phi.$$
By the generic vanishing order of $\Phi$ at a subvariety $V$ we mean
the Lelong number of $2\Theta_{\log\Phi}$ at a generic point of $V$.

\medbreak We say that at a point $P$ of $X$ the vanishing order of
$\Phi$ is precisely achieved if for some $m_0\in{\mathbb N}$ the
function
$$
\Phi_{m_0}:=\sum_{m=1}^{m_0}\varepsilon_m\sum_{j=1}^{q_m}\left|s^{(m)}_j\right|^{\frac{2}{m}}
$$
is comparable to $\Phi$ on some open neighborhood of $P$ in $X$ in
the sense that there exist some open neighborhood $U$ of $P$ in $X$
and some positive number $C$ such that
$$
\frac{1}{C}\,\Phi_{m_0}\leq\Phi\leq C\Phi_{m_0}
$$
on $U$. To indicate the value $m_0$, we also say that at the point
$P$ of $X$ the vanishing order of $\Phi$ is precisely achieved by
the $m_0$-th partial sum of $\Phi$ (or for the finite number $m_0$).
Sometimes we drop the adverb ``precisely'' and simply say that the
vanishing order of $\Phi$ is achieved. Sometimes we also
alternatively say that the precise vanishing order of $\Phi$ is
achieved.

\medbreak When we say that at a point the stable vanishing order is
precisely achieved, we mean that the vanishing order of $\Phi$ is
precisely achieved at that point.  We use the adjective ``stable''
when $\Phi$ is not explicitly used, because $\Phi$ involves an
infinite sum and involves all $\Gamma\left(X, mK_X\right)$ as
$m\to\infty$ so that we refer to the vanishing order of $\Phi$ as
the stable vanishing order.

\medbreak Note that if the vanishing order of $\Phi$ is precisely
achieved at $P$ by the $m_0$-th partial sum of $\Phi$ and if
$m_1=\left(m_0\right)!$, then
$$
\frac{1}{C^\prime}\Phi\leq\sum_{j=1}^{q_{m_1}}\left|s^{(m_1)}_j\right|^{\frac{2}{m_1}}\leq
C^\prime\Phi
$$
on some open neighborhood of $P$ for some positive number
$C^\prime$.

\bigbreak\noindent(3.7) {\it Theorem (Finite Generation as
Consequence of Precise Achievement of Stable Vanishing Order).}
Suppose the stable vanishing orders are precisely achieved at every
point of $X$ for some $m_0\in{\mathbb N}$. Denote
$\left(m_0\right)!$ by $m_1$. Then the canonical ring
$$
\bigoplus_{m=1}^\infty\Gamma\left(X, mK_X\right)
$$
is generated by
$$
\bigoplus_{m=1}^{\left(n+2\right)m_1}\Gamma\left(X, mK_X\right)
$$
and hence is finitely generated by the finite set of elements
$$
\left\{s^{(m)}_j\right\}_{1\leq m\leq m_1,\,1\leq j\leq q_m}.
$$

\medbreak\noindent{\it Proof.}  Let $e^{-\varphi}=\frac{1}{\Phi}$.
For $m>(n+2)m_1$ and any $s\in\Gamma\left(X, mK_X\right)$ we have
$$
\int_X\frac{\left|s\right|^2e^{-\left(m-\left(n+2\right)m_1-1\right)\varphi}}
{\left(\sum_{j=1}^{q_{m_1}}\left|s^{\left(m_1\right)}_j\right|^2\right)^{n+2}}<\infty,
$$
because $\left|s\right|^2\leq\tilde C\Phi^m$ on $X$ for some $\tilde
C$. By Skoda's theorem on ideal generation ((1.1) and (2.4)) there
exist
$$
h_1,\cdots,h_{q_{m_1}}\in\Gamma\left(X,\left(m-m_1\right)K_X\right)
$$
such that $s=\sum_{j=1}^{q_{m_1}}h_j s^{\left(m_1\right)}_j$.   If
$m-\left(n+2\right)m_1$ is still greater than $\left(n+2\right)m_1$,
we can apply the argument to each $h_j$ instead of $s$ until we get
$$
h^{(j_1,\cdots,j_\nu)}_1,\cdots,h^{(j_1,\cdots,j_\nu)}_{q_{m_1}}\in\Gamma\left(X,\left(m-m_1\left(\nu+1\right)\right)K_X\right)
$$
for $1\leq j_1,\cdots,j_\nu\leq q_{m_1}$ with $0\leq\nu< N$, where
$N=\left\lfloor\frac{m}{m_1}\right\rfloor$, such that
$$
s=\sum_{1\leq j_1,\cdots, j_N\leq
q_{m_1}}h^{(j_1,\cdots,j_{N-1})}_{j_N}\prod_{\lambda=1}^Ns^{\left(m_1\right)}_{j_\lambda}.
$$
Q.E.D.

\bigbreak Once we have the precise achievement of stable vanishing
orders, Skoda's theorem on ideal generation ((1.1) and (2.4)) can
also be applied with coefficients in a line bundle with a metric of
positive curvature current.  The following theorem is Theorem (3.7)
with such a twisting added.

\bigbreak\noindent(3.8) {\it Theorem (Twisted Finite Generation as
Consequence of Precise Achievement of Stable Vanishing Order).} Let
$E$ be a line bundle on $X$ with a metric $e^{-\chi}$ of positive
curvature current.  Suppose the vanishing orders of $\Phi$ are
precisely achieved at every point of $X$ for some $m_0\in{\mathbb
N}$.  Denote $\left(m_0\right)!$ by $m_1$. Then $\Gamma\left(X,
{\mathcal I}_{m\varphi+\chi}\left(mK_X+E\right)\right)$ is equal to
$$
\left(\Gamma\left(X,m_1 K_X\right)\right)^{p_m}\Gamma\left(X,
{\mathcal
I}_{\left(m-m_1p_m\right)\varphi+\chi}\left(\left(m-m_1p_m\right)K_X+E\right)\right)$$
for $m\geq\left(n+2\right)m_1$, where
$p_m=\left\lfloor\frac{m}{m_1}\right\rfloor-(n+2)$ and
$\varphi=\log\Phi$.

\bigbreak This theorem on twisted finite generation is for later use
in the discussion in (8.1)(vi) about the alternative approach to the
problem of the finite generation of the canonical ring by using the
extension techniques from the proof of the deformational invariance
of plurigenera.

\bigbreak

\bigbreak\noindent{\bf \S4. Reduction of Achievement of Vanishing
Order to Non-Vanishing Theorem on Hypersurface by Fujita Conjecture
Type Techniques}

\bigbreak\noindent(4.1) {\it Strategy to Prove Finite Generation.}
Our strategy to prove the finite generation of the canonical ring is
to show that the vanishing orders of $\Phi$ are achieved for some
finite $m_0$. We implement this strategy by descending induction on
the dimension of the subvariety $V$ of $X$ such that at points of
$X-V$ the vanishing orders of $\Phi$ are achieved in the sense (as
given in (3.6)) that at every point $P\in X-V$ there exist some open
neighborhood $U$ of $P$ in $X-V$ and some positive number $C$ and
some $m_0$ (which may depend on $P$) such that
$$
\frac{1}{C}\,\Phi_{m_0}\leq\Phi\leq C\Phi_{m_0}
$$
on $U$, where
$$
\Phi_{m_0}:=\sum_{m=1}^{m_0}\sum_{j=1}^{q_m}\left|s^{(m)}_j\right|^{\frac{2}{m}}.
$$
When we have $V$, we will proceed to prove that for each branch
$V_0$ of $V$ there is a subvariety $Z$ of codimension $\geq 1$ in
$V_0$ such that at every point of $V_0-Z$ the vanishing orders of
$\Phi$ are achieved for some finite $m^\prime$.  At the end of the
induction process we then invoke Theorem (3.7) to finish off the
proof of the finite generation of the canonical ring.

\bigbreak\noindent(4.2) {\it Difference Between the Hypersurface
Case and the Case of a Subvariety of Higher Codimension.} We will
first deal with the case when $V$ is a hypersurface. Though we can
always blow up a lower-dimensional $V$ to a hypersurface $\tilde
V$, yet there is the difficulty that the proper subvariety $\tilde
Z$ of $\tilde V$ (where stable vanishing orders are not yet
precisely achieved) may be projected onto $V$. So the
higher-dimensional case calls for an approach other than blow-up.
There is more in-depth discussion about this in (7.8).

\medbreak Across a hypersurface the vanishing order of each finite
sum $\Phi_m$ can be described by a single number, but across a local
submanifold $W$ of higher codimension in $X$ the vanishing orders of
each $\Phi_m$ at a point $P$ of $W$ have to be described by an
Artinian subscheme in the normal directions of $W$ at $P$.  The
Artinian subschemes used in these notes are unreduced complex
subspaces of ${\mathbb C}^n$ supported at a singular point.

\medbreak We will use the moduli space of Artinian subschemes in
the normal directions of $V$ and consider the points of $V$ where
these Artinian subschemes vary without jump in order to locate, in
an {\it a priori} manner, a countable union $E$ of subvarieties of
codimension $\geq 1$ in $V$ so that $\tilde Z$ is always projected
to inside $E$.  The notion of the continuous variation of an
Artinian subscheme without jump is introduced in (7.5.1) below.

\bigbreak\noindent(4.3) {\it Sketch of the Idea of Using Fujita
Conjecture Type Techniques.}  We now consider the case of precisely
achieving vanishing orders of $\Phi$ by some finite partial sum at a
generic point of an irreducible hypersurface $Y$ of $X$. Let
$\gamma$ be the generic Lelong number of $\Phi$ at $Y$.  To make the
argument more transparent, we break it up into two parts. First, we
assume that $\gamma$ is rational and show that stable vanishing
orders of $\Phi$ are precisely achieved by some finite partial sum
at a generic point of the irreducible hypersurface $Y$. This is done
by Fujita conjecture type techniques. Kawamata already used this
kind of techniques for the special case when $K_X$ is numerically
effective so that $\gamma=0$ in his paper [Kawamata 1985] before the
introduction of Fujita conjecture type techniques. Then we assume
that $\gamma$ is irrational and show how to modify the argument for
the case of a rational $\gamma$ to handle the case of an irrational
$\gamma$.

\medbreak We first sketch the main ideas here and then give the
details of the argument. Let us assume $\gamma$ rational (which
may be zero). When the stable vanishing order $\gamma$ at a
generic point of $Y$ is not precisely achieved by some finite sum
$\Phi_{m_0}$, the vanishing order $\gamma_{m_0}$ of the finite sum
$\Phi_{m_0}$ is always $>\gamma$ but approaches $\gamma$ as
$m_0\to\infty$. We are going to modify the metric $\frac{1}{\Phi}$
for $K_X$. The Fujita conjecture type techniques involve
constructing another metric $\frac{1}{\tilde\Phi}$ of $K_X$ so
that for some $m$ the multiplier ideal sheaf of the metric
$\frac{1}{\tilde\Phi^{m-1}}$ of $(m-1)K_X$ at a generic point of
$Y$ is the ideal sheaf of $\left(m\gamma+1\right)Y$.  This
actually means that instead of getting the vanishing order
$\gamma_{m_0}$ to be equal to $\gamma$ for some finite $m_0$ which
we cannot do, we settle for some extra order across $Y$, namely
for holomorphic sections of $(m-1)K_X$ over $X$ instead of the
expected vanishing order $(m-1)\gamma$ we allow the higher
vanishing $m\gamma+1$. The will be done in (4.4.1).  This extra
vanishing order enables us to use the theorem of
Kawamata-Viehweg-Nadel ([Kawamata 1982], [Viehweg 1982], [Nadel
1990]) with the use of an additional $K_X$ to conclude that
$$
H^1\left(X, {\mathcal J}\left(mK_X-m\gamma Y-Y\right)\right)=0,
$$
where ${\mathcal J}$ is a coherent ideal sheaf whose zero-set does
not contain $Y$.  This means that if we are able to get a
holomorphic section of $m\left(K_X-\gamma Y\right)$ over $Y$ which
is nonzero at some point of $Y$ (which locally belongs to the ideal
sheaf ${\mathcal J}$), then we can extend it to a holomorphic
section of $m\left(K_X-\gamma Y\right)$ over $X$, which implies that
the precise vanishing order of $\gamma$ at a generic point of $Y$ is
achieved by the finite sum $\Phi_m$.

\medbreak In the case of an irrational $\gamma$ we will use
Kronecker's theorem on diophantine approximation to make
$m\gamma-\left\lfloor m\gamma\right\rfloor$ small so that the
discrepancy can be handled by using the bigness of $K_X$ and the
Kodaira's decomposition of a big line bundle into an effective
${\mathbb Q}$-line-bundle and an ample ${\mathbb Q}$-line-bundle.
This will be done in (4.4.6).

\bigbreak\noindent(4.4) {\it Details of the Use of Fujita
Conjecture Type Techniques.}  The induction process to prove the
finite generation of the canonical ring is a descending induction
on the dimension of the subvariety in $X$ outside of which the
vanishing order of $\Phi$ is achieved by some finite partial sum
$\Phi_{m_0}$. The zero-set of $\Phi$ is a subvariety of $X$.
Clearly at points where $\Phi$ is nonzero, the vanishing order of
$\Phi$ is achieved by some finite partial sum $\Phi_{m_0}$. For
the initial induction step we take an irreducible hypersurface $Y$
in the zero-set of $\Phi$ if there is any.  By resolving the
singularities of $Y$, we can assume without loss of generality
that $Y$ is a nonsingular complex hypersurface of $X$. Let
$\gamma\geq 0$ be the generic Lelong number of $\Phi$ at $Y$. We
would like to show that at a generic point of $Y$ we can find some
$m$-canonical section over $X$ whose vanishing order at a generic
point of $Y$ is precisely $m\gamma$.  Suppose the contrary and we
are going to derive a contradiction.  We assume first that
$\gamma$ is rational. Choose a positive integer $m_1$ such that
$m_1\gamma$ is an integer. Let $L=m_1\left(K_X-\gamma Y\right)$.
Let $s_Y$ be the canonical section of the line bundle $Y$ on $X$.
We introduce the metric
$e^{-\tilde\varphi}=\left(\frac{\left|s_Y\right|^{2\gamma}}{\Phi}\right)^{m_1}$
of the line bundle $L=m_1\left(K_X-\gamma Y\right)$. The curvature
current $\Theta_{\tilde\varphi}$ is a closed positive
$(1,1)$-current on $X$.

\medbreak To reduce to a non-vanishing theorem the problem of
achieving precisely the stable vanishing order at a generic point of
a hypersurface $Y$, we have to introduce two metrics $e^{-\chi}$ and
$e^{-\xi}$ of $p_0L-K_X$ on $X$, which we will do respectively in
(4.4.1) and (4.4.2). The first metric $e^{-\chi}$ will have enough
singularity across $Y$ so that its multiplier ideal sheaf ${\mathcal
I}_\chi$ is contained in the ideal sheaf of $Y$ and at some generic
point of $Y$ equals the ideal sheaf of $Y$. The second metric
$e^{-\xi}$ will essentially have minimal singularity at a generic
point of $Y$ so that its multiplier ideal sheaf ${\mathcal I}_\xi$
has no zero at a generic point of $Y$. Moreover, the curvature
currents $\Theta_\chi$ and $\Theta_\xi$ of both metrics dominate
some smooth positive $(1,1)$-form on $X$. We apply the vanishing
theorem of Kawamata-Viehweg-Nadel to the line bundle
$\left(p+p_0\right)L-K_X$ with the metric $e^{-p\tilde\varphi-\chi}$
to get the surjectivity of
$$\Gamma\left(X,{\mathcal
I}_{p\tilde\varphi+\xi}\left(\left(p+p_0\right)L\right)\right)\to
\Gamma\left(X,\left({\mathcal
I}_{p\tilde\varphi+\xi}\left/{\mathcal
I}_{p\tilde\varphi+\chi}\right.\right)\left(\left(p+p_0\right)L\right)\right).
$$
If we have an element of
$$\Gamma\left(X,\left({\mathcal
I}_{p\tilde\varphi+\xi}\left/{\mathcal
I}_{p\tilde\varphi+\chi}\right.\right)\left(\left(p+p_0\right)L\right)\right)$$
which as a local holomorphic function on $Y$ is nonzero at some
generic point $P_0$ of $Y$, then the stable vanishing order is
precisely achieved at $P_0$.  This would finish the reduction to a
non-vanishing theorem the problem of precisely achieving the
stable vanishing order at a generic point of $Y$. However, there
is a small technical problem.

\medbreak The problem is that a nonzero element of
$$\Gamma\left(X,\left({\mathcal
I}_{p\tilde\varphi+\xi}\left/{\mathcal
I}_{p\tilde\varphi+\chi}\right.\right)\left(\left(p+p_0\right)L\right)\right)$$
may not be nonzero at some generic point $P_0$ of $Y$ as a local
holomorphic function on $Y$ when the sheaf ${\mathcal
I}_{p\tilde\varphi+\xi}\left/{\mathcal
I}_{p\tilde\varphi+\chi}\right.$ cannot be regarded as an ideal
sheaf on $Y$ but can only be regarded as an ideal sheaf over some
unreduced structure of $Y$.  (An unreduced structure means that
there are nonzero nilpotent elements in the structure sheaf.)  In
that case we have to replace the metric $e^{-\chi}$ with an
appropriate interpolation between $e^{-\chi}$ and $e^{-\xi}$ (and
possibly also apply some slight modification) to make the quotient
${\mathcal I}_{p\tilde\varphi+\xi}\left/{\mathcal
I}_{p\tilde\varphi+\chi}\right.$ the analog of a minimal center of
log canonical singularities.  The procedure of interpolation and
slight modification of metrics will be explained in (4.4.4).
Concerning the analog of minimal center of log canonical
singularities, we will discuss in (4.4.5) the r\^oles of unreduced
subspace and minimal centers of log canonical singularities in
Fujita conjecture type problems.

\medbreak After we solve the technical problem of unreduced
structures using the interpolation of metrics and the analog of
minimal centers of log canonical singularities, we can regard the
sheaf ${\mathcal I}_{p\tilde\varphi+\xi}\left/{\mathcal
I}_{p\tilde\varphi+\chi}\right.$ as an ideal sheaf on $Y$ or on some
complex subspace of $Y$ with reduced structure, then we need only
consider the problem of the existence of nonzero elements in
$$\Gamma\left(X,\left({\mathcal
I}_{p\tilde\varphi+\xi}\left/{\mathcal
I}_{p\tilde\varphi+\chi}\right.\right)\left(\left(p+p_0\right)L\right)\right).$$
This will be handled by our general non-vanishing theorem (6.2).

\bigbreak\noindent(4.4.1) {\it Lemma.} For a sufficiently positive
integer $p_0$ there exists some metric $e^{-\chi}$ of the line
bundle $p_0L-K_X$ on $X$ such that the curvature current
$\Theta_\chi=\frac{\sqrt{-1}}{2\pi}\partial\bar\partial\chi$ of
$e^{-\chi}$ dominates some smooth strictly positive $(1,1)$-form on
$X$ and the multiplier ideal sheaf ${\mathcal I}_\chi$ of the metric
$e^{-\chi}$ is contained in the ideal sheaf of $Y$ and at some
generic point of $Y$ equals the ideal sheaf of $Y$.

\medbreak\noindent{\it Proof.}  Since $X$ is of general type, the
canonical line bundle $K_X$ is big and we can write $K_X$ in the
form $K_X=E+A$, where $E$ is an effective ${\mathbb Q}$-divisor and
$A$ is an ample ${\mathbb Q}$-line-bundle.  Let $s_E$ be the
canonical divisor of the line bundle defined by $E$. Let $\eta$ be
the coefficient of $Y$ in $E$. The number $\eta$ must be greater
than $\gamma$, otherwise for some $m$ sufficiently large
$\Gamma\left(X, mA\right)$ generates $A$ and the subset
$\left(s_E\right)^m\Gamma\left(X, A\right)$ of $\Gamma\left(X,
mK_X\right)$ contains an element whose vanishing order at a generic
point of $Y$ is $m\eta$ across $Y$, which contradicts the assumption
that at a generic point of $Y$ there does not exist any
$m$-canonical section over $X$ with vanishing order $\leq m\gamma$.
Let $h_A$ be a smooth metric of $A$ with strictly positive curvature
form $\Theta_A$. Write $E=\eta Y+F$, where $F$ is an effective
${\mathbb Q}$-divisor. Let $h_A$ be a smooth metric of $A$ with
strictly positive curvature form $\Theta_A$.

\medbreak Let $p_0$ be a positive integer on which we will later
impose more conditions.  For any rational number $0<\beta\leq
p_0m_1-1$, the line bundle $p_0L-K_Y$ can be rewritten as
$$\displaylines{(4.4.1.1)\qquad\qquad p_0 L-K_Y=p_0m_1\left(K_X-\gamma Y\right)-\left(K_X+Y\right)
\qquad\qquad\qquad\qquad\hfill\cr=
\left(p_0m_1-1\right)\left(K_X-\gamma
Y\right)-\left(\gamma+1\right)Y\cr=
\left(p_0m_1-1-\beta\right)\left(K_X-\gamma Y\right)+\beta A+\beta
F-\left(\gamma+1-\beta\left(\eta-\gamma\right)\right)Y.\cr}
$$ Now choose $p_0$ and $0<\beta\leq p_0m_1-1$ such that $\gamma+1-\beta\left(\eta-\gamma\right)=0$.
Since $\eta>\gamma$, this choice of $p_0$ and $\beta$ is possible
when $p_0$ is sufficiently large.  Finally we define the metric
$$
e^{-\chi}=\frac{e^{-\frac{p_0m_1-1-\beta}{m_1}\,\tilde\varphi}\left(h_A\right)^\beta}{\left|s_Y\right|^2\left|s_F\right|^{2\beta}}
$$
which satisfies our requirement.  Q.E.D.

\bigbreak\noindent(4.4.2) {\it Lemma.} For a sufficiently large
positive integer $p_0$ there exists some metric $e^{-\xi}$ of the
line bundle $p_0L-K_X$ on $X$ such that the curvature current
$\Theta_\xi=\frac{\sqrt{-1}}{2\pi}\partial\bar\partial\xi$ of
$e^{-\xi}$ dominates some smooth strictly positive $(1,1)$-form on
$X$ and the multiplier ideal sheaf ${\mathcal I}_\xi$ is equal to
the structure sheaf ${\mathcal O}_X$ of $X$ at a generic point of
$Y$.

\medbreak\noindent{\it Proof.}  To construct the metric $e^{-\xi}$
of $p_0L-K_X$, we modify the construction of $e^{-\chi}$ in the
proof of (4.4.1) as follows.  We add $Y$ to both sides of (4.4.1.1)
and choose $0\leq\beta_1\leq p_0m_1-1$ with
$\gamma-\beta_1\left(\eta-\gamma\right)=0$ so that we can write
$$\displaylines{p_0 L-K_X=p_0m_1\left(K_X-\gamma Y\right)-K_X
\cr= \left(p_0m_1-1\right)\left(K_X-\gamma Y\right)-\gamma Y\cr=
\left(p_0m_1-1-\beta_1\right)\left(K_X-\gamma Y\right)+\beta_1
A+\beta_1
F-\left(\gamma-\beta_1\left(\eta-\gamma\right)\right)Y\cr
=\left(p_0m_1-1-\beta_1\right)\left(K_X-\gamma Y\right)+\beta_1
A+\beta_1 F.\cr}
$$
We now define the metric
$$
e^{-\xi}=
\frac{e^{-\frac{p_0m_1-1-\beta_1}{m_1}\,\tilde\varphi}\left(h_A\right)^\beta}{\left|s_F\right|^{2\beta_1}}
$$
which satisfies our requirement.  Q.E.D.

\bigbreak\noindent(4.4.3) {\it Interpolation of Metrics, Slight
Modification, and Minimal Center of Log Canonical Singularities.}
Suppose $M$ is a compact complex algebraic manifold and $E$ is a
holomorphic line bundle with two metrics $e^{-\kappa_1}$ and
$e^{-\kappa_2}$ such that each $\kappa_j$ is locally
plurisubharmonic (for $j=1,2$).   By an {\it interpolation of the
two metrics} $e^{-\kappa_1}$ and $e^{-\kappa_2}$ of $E$ we mean a
metric $e^{-\kappa_\eta}$ of $E$ of the form
$=e^{-\eta\kappa_1-\left(1-\eta\right)\kappa_2}$ with $0<\eta<1$.
Suppose $A$ is an ample ${\mathbb Q}$-line-bundle over $M$ with
metric $h_A$ and positive curvature form $\omega_A$ and $s_A$ is a
multi-valued holomorphic section of $A$ over $M$.  By a {\it
slight modification} of the interpolated metric $e^{-\kappa_\eta}$
of $E$ we mean a metric $e^{-\kappa_{\eta,s}}$ of $E$ either of
the form
$$\frac{e^{-\kappa_\eta}}{h_A\left|s\right|^2}$$ if the curvature
current $\Theta_{\kappa_\eta}$ of $e^{-\kappa_\eta}$ dominates
$(1+\varepsilon)\omega_A$ for some $\varepsilon>0$ or of the form
$e^{-\kappa_\eta}h_A\left|s\right|^2$ with
$\kappa_\eta-\log\left|s\right|^2$ locally plurisubharmonic on $M$
so that in both cases the curvature current of the new metric
$e^{-\kappa_{\eta,s}}$ of $E$ still dominates some positive multiple
of $\omega_A$.

\medbreak An interpolation of metrics is used in the following
context. Suppose the multiplier ideal sheaf ${\mathcal
I}_{\kappa_1}$ is properly contained in the multiplier ideal sheaf
${\mathcal I}_{\kappa_2}$.  We choose $0<\eta_0<1$ as the smallest
so that ${\mathcal I}_{\kappa_{\eta_0}}$ is different from
${\mathcal I}_{\kappa_2}$.  In other words, we choose $0<\eta_0<1$
as the smallest so that the support $Z_{\eta_0}$ of ${\mathcal
I}_{\kappa_2}\left/{\mathcal I}_{\kappa_{\eta_0}}\right.$ is still
nonempty.

\medbreak A slight modification of an interpolated metric is used in
the following context.  Suppose either $\Theta_{\kappa_1}$ or
$\Theta_{\kappa_2}$ dominates some smooth positive $(1,1)$-form on
$M$.  We choose the multi-valued holomorphic section $s$ of the
ample ${\mathbb Q}$-line-bundle $A$ so that the support
$Z_{\eta_0,s}$ of ${\mathcal I}_{\kappa_2}\left/{\mathcal
I}_{\kappa_{\eta_0,s}}\right.$ satisfies some additional properties.
One of such additional properties which is commonly used is that
there is some $\pi:\tilde M\to M$ obtained by a finite number of
successive monoidal transformations with nonsingular center such
that $Z_{\kappa_\eta,s}$ is the image under $\pi$ of some
nonsingular exceptional divisor $H$ in $\tilde M$.  The slight
modification helps us to go from a subvariety $Z_{\eta,s}$ to the
case of a complex manifold $H$.

\medbreak The use of interpolation and slight modification of
metrics corresponds to the process of getting the minimal center of
log canonical singularities in the techniques for Fujita conjecture
type problems and Shokurov's non-vanishing theorem (see {\it e.g.,}
[Shokurov 1985] and [Angehrn-Siu 1996] and [Kawamata 1997]).

\medbreak One important use of the interpolation and slight
modification process is to go from the case of a sheaf over an
unreduced space to the case of a sheaf over a reduced space which is
a subspace of the unreduced space (equivalently the structure sheaf
of the reduced subspace is a quotient of the structure sheaf of the
unreduced space). When it is a matter of finding a non-vanishing
section, we can start with a holomorphic section on a smaller
reduced space and extend it back. For example, when
$\Theta_{\kappa_1}$ dominates some smooth positive $(1,1)$-form on
$M$, we have the vanishing of $H^1\left(M,{\mathcal
I}_{\kappa_1}\left(E+K_M\right)\right)$ and $H^1\left(M,{\mathcal
I}_{\kappa_{\eta_0,s}}\left(E+K_M\right)\right)$ and therefore the
surjectivity of both restriction maps
$$
\displaylines{\Gamma\left(M,{\mathcal
I}_{\kappa_1}\left(E+K_M\right)\right)\to
\Gamma\left(M,\left({\mathcal I}_{\kappa_1}\left/{\mathcal
I}_{\kappa_2}\right.\right)\left(E+K_M\right)\right),\cr
\Gamma\left(M,{\mathcal I}_{\kappa_1}\left(E+K_M\right)\right)\to
\Gamma\left(M,\left({\mathcal I}_{\kappa_1}\left/{\mathcal
I}_{\kappa_{\eta_0,s}}\right.\right)\left(E+K_M\right)\right).\cr}
$$
When we are interested in getting some non identically zero
element of $\Gamma\left(M,{\mathcal
I}_{\kappa_1}\left(E+K_M\right)\right)$, we can either start with
a non identically zero element of $\Gamma\left(M,\left({\mathcal
I}_{\kappa_1}\left/{\mathcal
I}_{\kappa_2}\right.\right)\left(E+K_M\right)\right)$ or a non
identically zero element of $\Gamma\left(M,\left({\mathcal
I}_{\kappa_1}\left/{\mathcal
I}_{\kappa_{\eta_0,s}}\right.\right)\left(E+K_M\right)\right)$. By
replacing the former by the latter, we can work with a reduced
space (or even a manifold without multiplicity after using
blow-up). This way of handling the problem will also be used in
(6.10).

\bigbreak\noindent(4.4.4) {\it Remark on the R\^oles of Unreduced
Subspace and Minimal Centers of Log Canonical Singularities in
Fujita Conjecture Type Techniques.}  We would like to remark
further that the idea of unreduced structures and the
interpolations of metrics and slight modifications described in
(4.4.3) is actually the key point in the techniques for Fujita
conjecture type problems. Suppose $F$ is a holomorphic line bundle
with metric $e^{-\psi}$ over a compact complex algebraic manifold
$M$ so that the curvature current of $e^{-\psi}$ dominates some
smooth positive $(1,1)$-form on $M$.  The Fujita conjecture type
problem for this context is to study the question of the
generation of the multiplier ideal sheaf ${\mathcal I}_{m\psi}$ by
$\Gamma\left(M, {\mathcal I}_{m\psi}\left(mF+K_M\right)\right)$
for $m$ sufficiently large.

\medbreak While it is easy to get the vanishing of cohomology
$$
H^\nu\left(M, {\mathcal
I}_{m\psi}\left(mF+K_M\right)\right)=0\quad{\rm for\ \ }\nu\geq 1
$$
by using the theorem of Kawamata-Viehweg-Nadel which requires only
some small positive lower bound of the curvature current of the
metric $e^{-\psi}$.  It is difficult to construct elements of
$\Gamma\left(M, {\mathcal I}_{m\psi}\left(mF+K_M\right)\right)$, for
which a sufficiently large positive lower bound of the curvature
current of the metric $e^{-\psi}$ is usually needed (see (9.1.1)).

\medbreak The way to produce sections is to find another metric
$e^{-\psi_1}$ of $F$ (again with curvature current dominating some
smooth positive $(1,1)$-form on $M$) such that the multiplier ideal
sheaf ${\mathcal I}_{m\psi_1}$ is contains in the multiplier ideal
sheaf ${\mathcal I}_{m\psi}$.  Then the vanishing of positive-degree
cohomology for $mF+K_M$ with coefficients in both multiplier ideal
sheaves ${\mathcal I}_{m\psi}$ and ${\mathcal I}_{m\psi_1}$ gives
$$
H^1\left(M, \left({\mathcal I}_{m\psi}\left/{\mathcal
I}_{m\psi_1}\right.\right)\left(mF+K_M\right)\right)=0
$$
and the surjectivity of
$$
\Gamma\left(M, {\mathcal I}_{m\psi}\left(mF+K_M\right)\right)\to
\Gamma\left(M, \left({\mathcal I}_{m\psi}\left/{\mathcal
I}_{m\psi_1}\right.\right)\left(mF+K_M\right)\right).
$$
Still it is difficult to produce elements of $\Gamma\left(M,
\left({\mathcal I}_{m\psi}\left/{\mathcal
I}_{m\psi_1}\right.\right)\left(mF+K_M\right)\right)$ unless the
support of ${\mathcal I}_{m\psi}\left/{\mathcal I}_{m\psi_1}\right.$
is isolated.  If for some point $P_0\in M$ the stalk of ${\mathcal
I}_{m\psi_1}$ at $P_0$ belongs to a sufficiently high power of the
maximum ideal ${\mathfrak m}_{M,P_0}$ of $M$ at $P_0$, then when we
change the metric $e^{-\psi_1}$ by the interpolation of metrics and
slight modification to make the support of ${\mathcal
I}_{m\psi}\left/{\mathcal I}_{m\psi_1}\right.$ minimal, we can end
up with the support of ${\mathcal I}_{m\psi}\left/{\mathcal
I}_{m\psi_1}\right.$ isolated at $P_0$.

\medbreak The usual method to approach Fujita conjecture type
problems is to use the theorem of Riemann-Roch to produce
multi-valued holomorphic sections of $F$ vanishing to high order and
then use them to modify the metric $e^{-\psi}$ to construct
$e^{-\psi_1}$.

\medbreak For our notes, as we will see in the proof of our general
non-vanishing theorem (6.2), we are following the same strategy with
$mF$ replaced by $mF-K_M$ to help create some positive lower bound
for the curvature current and, much more importantly, to incorporate
Shokurov's technique [Shokurov 1985] of comparing the application of
the theorem of Riemann-Roch to a line bundle and its twisting by a
flat line bundle when we fail to produce a multi-valued holomorphic
section of $mF-K_X$ vanishing to high order.  We will inductively
reduce the dimension of the support of ${\mathcal
I}_{m\psi}\left/{\mathcal I}_{m\psi_1}\right.$ until we are forced
to use Shokurov's technique.  During the induction process, we do
not have control over where the support of ${\mathcal
I}_{m\psi}\left/{\mathcal I}_{m\psi_1}\right.$ is when we make it
minimal, except that it must contain the point where the
multi-valued holomorphic section of $mF-K_X$ vanishes to high order.

\bigbreak\noindent(4.4.6) {\it Modification for an Irrational
Stable Vanishing Order.}  We now explain what modification is
needed to get rid of the additional assumption in the first
paragraph of (4.4) that $\gamma$ is rational. In that case we use
Kronecker's diophantine approximation and replace
$m\left(K_X-\gamma Y\right)$ by $mK_X-\left\lfloor
m\gamma\right\rfloor Y$.  For any $m\in{\mathbb N}$ let
$\varepsilon_m$ denote the positive number $m\gamma-\left\lfloor
m\gamma\right\rfloor$.  First of all, by Kronecker's diophantine
approximation (see (5.2.1.1)), for any given $\varepsilon>0$ there
exists some $m$ such that $\varepsilon_m<\varepsilon$.

\medbreak Let $h_Y$ be any smooth metric for the line bundle $Y$
over $X$.  For the modification we let $m=\left(p+p_0\right)m_1$ and
replace the metric
$$
e^{-p\tilde\varphi-\chi}=\frac{e^{-\left(p+\frac{p_0m_1-1-\beta}{m_1}\right)\tilde\varphi}\left(h_A\right)^\beta}{\left|s_Y\right|^2\left|s_F\right|^{2\beta}}
$$
of $\left(p+p_0\right)L-K_X-Y$ given above in (4.1.1) for the case
of a rational $\gamma$ by the new metric
$$
e^{-p\tilde\varphi-\chi}=\frac{e^{-\left(p+\frac{p_0m_1-1-\beta}{m_1}\right)\tilde\varphi}\left(h_A\right)^\beta\,\left(h_Y\right)^{\varepsilon_m}}{\left|s_Y\right|^2\left|s_F\right|^{2\beta}}
$$
of $\left(mK_X-\left\lfloor m\gamma\right\rfloor Y\right)-K_X-Y$.
Likewise, we replace the metric
$$
e^{-p\tilde\varphi-\xi}=\frac{e^{-\left(p+\frac{p_0m_1-1-\beta_1}{m_1}\right)\tilde\varphi}\left(h_A\right)^{\beta_1}}{\left|s_Y\right|^2\left|s_F\right|^{2\beta_1}}
$$
of $\left(p+p_0\right)L-K_X-Y$ given above in (4.1.2) for the case
of a rational $\gamma$ by the new metric
$$
e^{-p\tilde\varphi-\xi}=\frac{e^{-\left(p+\frac{p_0m_1-1-\beta_1}{m_1}\right)\tilde\varphi}\left(h_A\right)^{\beta_1}\,\left(h_Y\right)^{\varepsilon_m}}{\left|s_Y\right|^2\left|s_F\right|^{2\beta_1}}
$$
of $\left(mK_X-\left\lfloor m\gamma\right\rfloor Y\right)-K_X-Y$.
All we have to do is to choose $m$ so that the consequence
(5.2.1.1) of Kronecker's theorem on diophantine approximation can
be applied to make $\varepsilon_m$ small enough for the curvature
form of the metric
$\left(h_A\right)^\beta\,\left(h_Y\right)^{\varepsilon_m}$ to be
strictly positive on $X$.

\medbreak The main point of the modification is to use the small
positivity squeezed out from the bigness of $K_Y$ to absorb the
contribution from the small left-over non-integral part
$\varepsilon_m Y$ of $m\gamma Y$ when $\gamma$ is irrational.

\medbreak After the use of the non-vanishing theorem to be
presented below (with the modification discussed in (4.4.3) and
(4.4.4)), the final result is that there exists some element of
$\Gamma\left(X,mK_X\right)$ whose vanishing order across $Y$ at a
generic point of $Y$ is $\left\lfloor m\gamma\right\rfloor$, which
would imply that the case of an irrational $\gamma$ cannot occur.

\bigbreak

\bigbreak\noindent{\bf \S5. Diophantine Approximation of Kronecker}

\bigbreak Before we present our general non-vanishing theorem (6.2),
we first derive some results on diophantine approximation by using
the diophantine approximation theorem of Kronecker which will be
needed in the presentation of our general non-vanishing theorem in
\S6. The key result due to Kronecker which we will use is the
following.  A reference is Theorem 444 on p. 382 of [Hardy-Wright
1960].

\bigbreak\noindent(5.1) {\it Theorem (Kronecker).}  Let
$a_1,\cdots,a_N$ be ${\mathbb Q}$-linearly independent real numbers.
Let $b_1,\cdots,b_N\in{\mathbb R}$.  Let $\varepsilon, T$ be
positive numbers. Then we can find $t>T$ and integers
$x_1,\cdots,x_N$ such that
$\left|ta_j-b_j-x_j\right|\leq\varepsilon$ for $1\leq j\leq N$.

\bigbreak\noindent(5.2) {\it Remark.} Note that the result implies
that for the case of a single irrational real number $\gamma$, the
image of the set $\gamma{\mathbb Z}$ in ${\mathbb R}\left/{\mathbb
Z}\right.$ is dense, because we can choose
$$N=2,\quad a_1=1,\quad a_2=\gamma,\quad b_1=b_2=0,\quad T=1$$
and conclude that given any $0<\varepsilon<1$ there exists some
$t>T$ such that $\left|t
a_j-p_j\right|<\frac{\varepsilon}{1+\left|\gamma\right|}$ ($j=1,2$)
for some $p_1, p_2\in{\mathbb Z}$ and, as a consequence,
$$\left|p_1\gamma-p_2\right|=\left|t\gamma-p_2-\left(t-p_1\right)\gamma\right|
\leq\frac{\varepsilon}{1+\left|\gamma\right|}+\left|\gamma\right|\frac{\varepsilon}{1+\left|\gamma\right|}=\varepsilon.$$
Of course, geometrically it is simply the well-known statement that
the straight line in ${\mathbb R}^2$ whose slope is an irrational
real number $\gamma$ has dense image in $\left({\mathbb
R}\left/{\mathbb Z}\right.\right)^2$ (as one can see by looking at
the straight lines passing through the points $(x,0)$ with the image
of $x\in{\mathbb R}$ in ${\mathbb R}\left/{\mathbb Z}\right.$
belonging to the image of $\gamma{\mathbb Z}$ in ${\mathbb
R}\left/{\mathbb Z}\right.$).

\bigbreak\noindent(5.2.1)  When we choose
$b_2=\frac{2\varepsilon}{3}$ instead of $b_2=0$, we can find $t>T$
such that $\left|t
a_j-p_j\right|<\frac{\varepsilon}{3\left(1+\left|\gamma\right|\right)}$
($j=1,2$) for some $p_1, p_2\in{\mathbb Z}$ and, as a consequence,
$$\left|p_1\gamma-p_2-b_2\right|=\left|t\gamma-p_2-\left(t-p_1\right)\gamma-b_2\right|
\leq\frac{\varepsilon}{3\left(1+\left|\gamma\right|\right)}+
\left|\gamma\right|\frac{\varepsilon}{3\left(1+\left|\gamma\right|\right)}=\frac{\varepsilon}{3}.$$
Thus
$$p_1\gamma\geq
p_2+b_2-\frac{\varepsilon}{3}=p_2+\frac{2\varepsilon}{3}-\frac{\varepsilon}{3}=p_2+\frac{\varepsilon}{3}
$$
and $\left\lfloor p_1\gamma\right\rfloor\geq p_2$ and
$$p_1\gamma-\left\lfloor p_1\gamma\right\rfloor\leq
p_1\gamma-p_2<b_2+\frac{\varepsilon}{3}=\frac{2\varepsilon}{3}+\frac{\varepsilon}{3}=\varepsilon.$$
We conclude that

\bigbreak\noindent(5.2.1.1) for any $\varepsilon>0$, any irrational
real number $\gamma$, and any $T>0$, there exists some integer $m>T$
such that $m\gamma-\left\lfloor m\gamma\right\rfloor<\varepsilon$.

\bigbreak\noindent(5.3) {\it Corollary.} Let $1=a_0, a_1,\cdots,a_N$
be ${\mathbb Q}$-linearly independent real numbers. Let
$0<\varepsilon<\frac{1}{2}$. Then we can find a positive integer $m$
such that $ma_j-\left\lfloor
ma_j\right\rfloor>\frac{1}{2}-\varepsilon$ for $1\leq j\leq N$.

\medbreak\noindent{\it Proof.} Let $\eta=\max_{0\leq j\leq
N}\left|a_j\right|$. Choose $b_0=0$ and $b_j=\frac{1}{2}$ for $1\leq
j\leq N$ and $\frac{1}{q}<\frac{\varepsilon}{2\eta}$ and $T=1$.  By
Kronecker's theorem applied to $1=a_0, a_1,\cdots,a_N$ we can find a
real number $t>T$ and integers $x_0, x_1,\cdots,x_N$ such that
$\left|ta_j-b_j-x_j\right|\leq\frac{1}{q}$ for $0\leq j\leq N$. We
now set $m=x_0$.  Then
$\left|t-m\right|<\frac{1}{q}<\frac{\varepsilon}{2\eta}$ so that
$$
\displaylines{\left|ma_j-\frac{1}{2}-x_j\right|\leq
\left|t-m\right|\left|a_j\right|+\left|ta_j-\frac{1}{2}-x_j\right|\cr
\leq \frac{\varepsilon}{2\eta}\left|a_j\right|+\frac{1}{q} \leq
\frac{\varepsilon}{2}+\frac{\varepsilon}{2}=\varepsilon.\cr}
$$
Since $x_j$ is an integer for $1\leq j\leq N$, we have
$ma_j-\left\lfloor ma_j\right\rfloor>\frac{1}{2}-\varepsilon$ for
$1\leq j\leq N$.  Since $m$ is an integer and
$\left|t-m\right|<\frac{1}{q}$ and $t>1$, it follows that $m$ is
positive. Q.E.D.

\bigbreak\noindent(5.4) {\it Lemma.} Let $\gamma$ be an irrational
positive number and let $a_1,\cdots,a_N$ be positive rational
numbers. Then there exists some positive integer $m$ such that
$m\gamma a_j-\left\lfloor m\gamma a_j\right\rfloor\geq\frac{1}{2}$
for $1\leq j\leq N$.

\medbreak\noindent{\it Proof.}  Let $p$ be the least common multiple
of the denominators of $a_1,\cdots,a_N$.  By replacing each $a_j$ by
$pa_j$, we can assume without loss of generality that each $a_j$ is
a positive integer.  Choose $\varepsilon>0$ so that
$\varepsilon\left|a_j\right|<\frac{1}{2}$ for $1\leq j\leq N$. Since
$\gamma$ is irrational, the set $\gamma{\mathbb Z}$ has dense image
in ${\mathbb R}\left/{\mathbb Z}\right.$. It follows that there
exists some $m,n\in{\mathbb Z}$ such that
$m\gamma-n\in\left(1-\varepsilon,\,1\right)$. Since $\gamma$ is
positive, the integer $m$ must be positive. From
$1-\varepsilon+n<m\gamma<1+n$ it follows that
$\left(n+1\right)a_j-\varepsilon a_j <m\gamma
a_j<\left(n+1\right)a_j$.  Then $\left\lfloor m\gamma
a_j\right\rfloor=\left(n+1\right)a_j-1$ and $$m\gamma
a_j-\left\lfloor m\gamma a_j\right\rfloor >1-\varepsilon a_j\geq
1-\varepsilon\left|a_j\right|>1-\frac{1}{2}=\frac{1}{2}.
$$
Q.E.D.

\bigbreak\noindent(5.5) {\it Lemma.}  Let $\gamma_j$ ($1\leq
j<\infty$) be a sequence of positive numbers and $\Lambda$ be a
positive integer such that $1,\gamma_1,\cdots,\gamma_\Lambda$ are
${\mathbb Q}$-linearly independent and
$$
\gamma_j=\sum_{\lambda=1}^{\Lambda}c_{j,\lambda}\gamma_\lambda\leqno{(5.5.1)}
$$
for $\Lambda<j<\infty$, where $c_{j,k}\in{\mathbb Q}$.  For any
positive integer $N$ there exists some positive integer $m$ such
that $$m \gamma_j- \left\lfloor m
\gamma_j\right\rfloor\geq\frac{1}{4}\quad{\rm for\ \ }1\leq j\leq
N.$$

\medbreak\noindent{\it Proof.} Choose a positive numbers
$\varepsilon$ which we will specify more precisely later. Choose a
rational number $r_\lambda$ so close to $\frac{1}{\gamma_\lambda}$
such that $\left|r_\lambda\gamma_\lambda-1\right|<\varepsilon$ and
$\left|\frac{1}{r_\lambda}-\gamma_\lambda\right|<\varepsilon$ for
$1\leq\lambda\leq\Lambda$. Let
$\gamma_\lambda^\prime=r_\lambda\gamma_\lambda$ for
$1\leq\lambda\leq\Lambda$ and let
$c_{j,\lambda}^\prime=\frac{c_{j,\lambda}}{r_\lambda}$ so that by
(5.5.1) we have
$$
\gamma_j=\sum_{\lambda=1}^{\Lambda}c^\prime_{j,\lambda}\gamma^\prime_\lambda.
$$
We have
$$
\displaylines{\sum_{\lambda=1}^{\Lambda}c^\prime_{j,\lambda}=
\sum_{\lambda=1}^{\Lambda}c^\prime_{j,\lambda}\gamma^\prime_\lambda
+\sum_{\lambda=1}^{\Lambda}c^\prime_{j,\lambda}\left(1-\gamma^\prime_\lambda\right)
\cr=\gamma_j
+\sum_{\lambda=1}^{\Lambda}c^\prime_{j,\lambda}\left(1-\gamma^\prime_\lambda\right)
\geq \gamma_j-\varepsilon
\sum_{\lambda=1}^{\Lambda}\left|c^\prime_{j,\lambda}\right| \cr=
\gamma_j-\varepsilon
\sum_{\lambda=1}^{\Lambda}\left|c_{j,\lambda}\frac{1}{r_\lambda}\right|\geq
\gamma_j-\varepsilon
\sum_{\lambda=1}^{\Lambda}\left|c_{j,\lambda}\right|\left(\varepsilon+\gamma_\lambda\right).\cr}
$$
We can choose $\varepsilon$ sufficiently small so that
$$
\varepsilon
\sum_{\lambda=1}^{\Lambda}\left|c_{j,\lambda}\right|\left(\varepsilon+\gamma_\lambda\right)<\frac{\gamma_j}{2}
\quad{\rm for\ \ }1\leq j\leq N.
$$
Thus
$$
\sum_{\lambda=1}^{\Lambda}c^\prime_{j,\lambda}\geq\frac{\gamma_j}{2}>0\quad{\rm
for\ \ }1\leq j\leq N.
$$
Let $q_N$ be the least common multiple of the denominator of
$c^\prime_{j,\lambda}$ for $1\leq\lambda\leq\Lambda$ and $1\leq
j\leq N$.  Then
$$
q_N\gamma_j=\sum_{\lambda=1}^{\Lambda}\left(q_Nc^\prime_{j,\lambda}\right)\gamma^\prime_\lambda,\leqno{(5.5.2)}
$$
where $q_Nc^\prime_{j,\lambda}$ is an integer.

\medbreak We choose any irrational number $0<\eta<1$. By Lemma (5.4)
there exists a positive integer $m_1$ such that
$$
m_1\eta\sum_{\lambda=1}^{\Lambda}\left(q_Nc^\prime_{j,\lambda}\right)-
\left\lfloor
m_1\eta\sum_{\lambda=1}^{\Lambda}\left(q_Nc^\prime_{j,\lambda}\right)\right\rfloor\geq\frac{1}{2}\leqno{(5.5.3)}
$$
for $1\leq j\leq N$.  Let $0<\delta<\frac{1}{4}$ such that $\delta$
is less than
$$
\left\lceil
m_1\eta\sum_{\lambda=1}^{\Lambda}\left(q_Nc^\prime_{j,\lambda}\right)\right\rceil-m_1\eta\sum_{\lambda=1}^{\Lambda}\left(q_Nc^\prime_{j,\lambda}\right)
$$
for $1\leq j\leq N$.  Let $C$ be the maximum of
$m_1\sum_{\lambda=1}^{\Lambda}\left|q_Nc^\prime_{j,\lambda}\right|$
for $1\leq j\leq N$.

\medbreak We apply Kronecker's theorem in (5.1) to the ${\mathbb
Q}$-linearly independent irrational numbers
$\gamma^\prime_1,\cdots,\gamma^\prime_\Lambda$ so that there exist
some $m_2\in{\mathbb Z}$ and integers $x_\lambda$ for $1\leq
\lambda\leq\Lambda$ such that
$$\left|m_2\gamma^\prime_\lambda-\left(x_\lambda+\eta\right)\right|
<\frac{\delta}{C}\quad{\rm for\ \ }1\leq\lambda\leq\Lambda,$$ from
which we obtain
$$\left|\sum_{\lambda=1}^{\Lambda}q_Nc^\prime_{j,\lambda}
\left(m_2\gamma^\prime_\lambda-\left(x_\lambda+\eta\right)\right)\right|
<\frac{\delta}{C}\,\sum_{\lambda=1}^{\Lambda}q_N\left|c^\prime_{j,\lambda}\right|.$$
Using $(5.5.2)$, we get
$$\left|m_2
q_N\gamma_j-\sum_{\lambda=1}^{\Lambda}\left(q_Nc^\prime_{j,\lambda}x_\lambda\right)-
\eta\sum_{\lambda=1}^{\Lambda}\left(q_Nc^\prime_{j,\lambda}\right)\right|
<\frac{\delta}{C}\sum_{\lambda=1}^{\Lambda}\left|q_Nc^\prime_{j,\lambda}\right|\leq\frac{\delta}{m_1}$$
for $1\leq j\leq N$ and
$$\left|m_1 m_2 q_N\gamma_j
-m_1\sum_{\lambda=1}^{\Lambda}\left(q_Nc^\prime_{j,\lambda}x_\lambda\right)-m_1
\eta\sum_{\lambda=1}^{\Lambda}\left(q_Nc^\prime_{j,\lambda}\right)\right|
\leq\delta$$ for $1\leq j\leq N$.  From $(5.5.3)$ and the choice of
$\delta$ it follows that
$$m_1 m_2 q_N\gamma_j-
\left\lfloor m_1 m_2
q_N\gamma_j\right\rfloor\geq\frac{1}{4}\quad{\rm for\ \ }1\leq j\leq
N.$$ The required $m$ can be set to be $m_1 m_2 q_N$.  Q.E.D.

\bigbreak

\bigbreak\noindent{\bf \S6 A General Non-Vanishing Theorem.}

\bigbreak Before we introduce our general non-vanishing theorem, we
first introduce some terminology.  It is to avoid the situation of a
change of the multiplier ideal sheaf by some small perturbation of
the metric, for example in the case of a metric like
$\frac{1}{\left|z_1\right|^2}$.  The terminology is that of a stable
metric for which a sufficiently small perturbation would not change
its multiplier ideal sheaf.

\bigbreak\noindent(6.1) {\it Definition.}  A metric $e^{-\varphi}$
of a line bundle $L$ on a compact K\"ahler manifold $M$ is said to
be {\it stable} if there exists some $\varepsilon>0$ with the
following property.  If $U$ is an open neighborhood of a point $P\in
M$, and $\kappa$ and $\psi$ are plurisubharmonic functions on $U$
such that the total mass, with respect to the K\"ahler form of $M$,
of the difference of the two closed positive $(1,1)$-currents
$\Theta_\varphi$ and $\Theta_\psi$ is less than $\epsilon$, then
there exists an open neighborhood $U^\prime$ of $P$ in $M$ such that
the multiplier ideal sheaf ${\mathcal I}_{\kappa+\psi}$ of the
metric $e^{-\kappa-\psi}$ on $U^\prime$ is equal to the multiplier
ideal sheaf ${\mathcal I}_{\kappa+\varphi}$ of the metric
$e^{-\varphi}$ on $U^\prime$, where
$\Theta_\varphi:=\frac{\sqrt{-1}}{2\pi}\partial\bar\partial\varphi$
and $\Theta_\psi:=\frac{\sqrt{-1}}{2\pi}\partial\bar\partial\psi$
are respectively the curvature currents of the metrics
$e^{-\varphi}$ and $e^{-\psi}$.

\bigbreak\noindent(6.2) {\it Theorem (Non-Vanishing).}  Let $M$ be a
compact complex projective algebraic manifold.  Let $L$ be a
holomorphic line bundle on $M$ with a (possibly singular) metric
$e^{-\varphi}$ along its fibers whose curvature current
$\Theta_\varphi$ is a closed positive $(1,1)$-current.  Assume that
for some positive integer $p_0$ there is a (possibly singular)
metric $e^{-\chi}$ along the fibers of $p_0L-K_M$ which is stable
and whose curvature current $\Theta_\chi$ is a closed positive
$(1,1)$-current which dominates some strictly positive smooth
$(1,1)$-form on $M$.  Assume that for any nonnegative integer $p$
the multiplier ideal sheaf ${\mathcal I}_{p\varphi+\chi}$ of the
metric $e^{-p\varphi-\chi}$ of the line bundle
$\left(p+p_0\right)L-K_M$ contains the multiplier ideal sheaf
${\mathcal I}_{\left(p+p_0\right)\varphi}$ of the metric
$e^{-\left(p+p_0\right)\varphi}$ of the line bundle
$\left(p+p_0\right)L$. Let $\varepsilon>0$. Then for some
sufficiently divisible $m\in{\mathbb N}$ the line bundle
$\left(m+p_0\right)L$ admits a non identically zero holomorphic
section on $M$ which belongs to the multiplier ideal of
$e^{-m\varphi-\chi}$. Moreover, $m$ can be chosen to satisfy the
stronger condition that global holomorphic sections of ${\mathcal
I}_{m\varphi+\chi}\left(mL\right)$ over $M$ generate ${\mathcal
I}_{m\varphi+\chi}$ outside some subvariety of codimension at least
$2$ in $M$.

\bigbreak\noindent(6.3) {\it Remark on the Assumption of the
Multiplier Ideal Sheaf of $e^{-p\varphi-\chi}$ to Contain That of
$e^{-\left(p+p_0\right)\varphi}$.}  The reason for this assumption
is to tailor our general non-vanishing theorem to its application to
the proof of the finite generation of the canonical ring for the
case of general type.  To illustrate the reason for this assumption,
let us consider the case of $M=X$ and $L=K_X$ and
$e^{-\varphi}=\frac{1}{\Phi}$, where $\Phi$ is as defined in the
beginning part of \S3.  In this case we choose $p_0>1$ and
$p_0L-K_M$ is simply $\left(p_0-1\right)K_X$ and we can choose a
sufficiently small $\varepsilon>0$ and set
$$
e^{-\chi}=\frac{e^{-\left(p_0-1-\varepsilon\right)\varphi}\left(h_A\right)^\varepsilon}
{\left|s_D\right|^{2\varepsilon}},
$$
where $K_X=A+D$ is the decomposition into an effective ${\mathbb
Q}$-divisor $D$ and an ample ${\mathbb Q}$-divisor $A$.  When
$\varepsilon$ is sufficiently small, for any nonnegative integer $p$
the multiplier ideal sheaf ${\mathcal I}_{p\varphi+\chi}$ of the
metric
$$
e^{-p\varphi-\chi}=\frac{e^{-\left(p+p_0-1-\varepsilon\right)\varphi}\left(h_A\right)^\varepsilon}
{\left|s_D\right|^{2\varepsilon}}
$$
contains the multiplier ideal sheaf ${\mathcal
I}_{\left(p+p_0\right)\varphi}$ of the metric
$e^{-\left(p+p_0\right)\varphi}$ of the line bundle
$\left(p+p_0\right)L$.  The reason of this relation between the two
multiplier ideal sheaves is, of course, due to the fact that the
metric $e^{-\varphi}$ of $K_X$ is the least singular of all the
metrics of $K_X$ with positive curvature current and, in particular,
the metric $e^{-\varphi}$ of $K_X$ is no more singular than the
metric $\frac{h_A}{\left|s_D\right|^2}$ of $K_X$.

\medbreak This condition is the minimum singularity condition for
the metric $e^{-\varphi}$.  In the induction steps to prove the
general non-vanishing theorem (.62) and also in the induction steps
in our analytic proof of the finite generation of the canonical ring
for the case of general type, we have to use induction.  In each
induction step we use the Fujita conjecture type techniques and the
analog of minimal centers of log canonical singularities (or after
removing the necessary minimum vanishing orders) to get to a
subspace.  Since the choice of the subspace is to guarantee the
extendability of sections from the subspace to the original ambient
space, the induced metric also inherits the minimum singularity
condition.

\bigbreak\noindent(6.4) {\it Remark on Requiring Constructed
Sections to Belong to Multiplier Ideal Sheaves.} The conclusion of
the general non-vanishing theorem in (6.2) requires more than the
mere existence of a non identically zero holomorphic section of the
line bundle in question.  It requires the constructed holomorphic
sections to belong to the appropriate multiplier ideal sheaves for
two reasons.

\medbreak The first reason is to enable us to carry out the
induction argument because to extend a holomorphic section from a
lower-dimensional subspace to a higher dimensional space, we need to
make sure that it belongs to the appropriate multiplier ideal sheaf
in order to use the theorem of Kawamata-Viehweg-Nadel for the
extension.

\medbreak The second reason is for the application of the general
non-vanishing theorem and the extension of sections to conclude that
generic stable vanishing orders are achieved by finite sums of the
absolute-value squares of fractional holomorphic pluricanonical
sections.  Again we need to make sure that the constructed sections
belong to the appropriate multiplier ideal sheaf in order to use the
theorem of Kawamata-Viehweg-Nadel for the extension to the entire
manifold.

\bigbreak\noindent(6.5) {\it Beginning of the Proof and the
Dichotomy According to the Structure of the Curvature Current.} We
now start with the proof of the general non-vanishing theorem stated
in (6.2). We will differentiate among several cases according to the
structure of the curvature current $\Theta_\varphi$ of the metric
$e^{-\varphi}$ of $L$.  According to Theorem (3.3), we have the
following decomposition of the curvature current $\Theta_\varphi$ of
the metric $e^{-\varphi}$ of $L$.
$$
\Theta_\varphi=\sum_{j=1}^J\tau_j\left[V_j\right]+R,
$$
where $J\in{\mathbb N}\cup\left\{0,\infty\right\}$ and the Lelong
number of $R$ is zero outside a countable union $Z$ of subvarieties
of codimension at least two in $M$ and $V_j$ is an irreducible
hypersurface in $M$ and $\tau_j>0$. The assumption concerning the
metric $e^{-\chi}$ of $p_0L-K_M$ whose curvature current has a
positive lower bound is introduced in order to be able to apply the
vanishing theorem of Kawamata-Viehweg-Nadel.

\medbreak First we make some remarks about the proof.  The main idea
of the proof is to use \begin{itemize}\item[(i)] the techniques for
Fujita's conjecture (see {\it e.g.,} [Angehrn-Siu 1995]) together
with \item[(ii)] the technique of Shokurov of using the theorem of
Riemann-Roch to compare the arithmetic genus of the line bundle and
that of its twisting by a flat line bundle [Shokurov
1985].\end{itemize} For the arguments of the Fujita conjecture,
given an ample line bundle $E$ over a compact complex algebraic
manifold $Y$ we use the theorem of Riemann-Roch to get the lower
bound of $\dim_{\mathbb C}\Gamma\left(Y, mE\right)$ as a function of
$m$ for $m$ sufficiently large to produce a multi-valued holomorphic
section of a multiple of $E$ vanishing to high order.  Now in our
case our line bundle $L$ is not ample but with a weaker condition of
having a metric $e^{-\varphi}$ of nonnegative curvature current. We
will, however, still follow the path of getting a good lower bound
of the dimension of the section module.

\medbreak The new ingredient is the dichotomy into the case of $R=0$
and the case of $R\not=0$ and, in the case of $R=0$, the further
dichotomy into the case of $J=\infty$ and the case of
$J\not=\infty$. For the case of $R\not=0$ or $J=\infty$ we still get
a good lower bound for the dimension of the section module. For the
case of $R=0$ and $J\not=\infty$ we explicitly construct by hand a
holomorphic section by following the technique of Shokurov of using
the theorem of Riemann-Roch to compare the arithmetic genus of the
line bundle and that of its twisting by a flat line bundle. For our
purpose we are only interested in producing a single nontrivial
holomorphic section (which belongs to the appropriate multiplier
ideal sheaf).

\bigbreak\noindent(6.6) {\it Slicing by an Ample Divisor.}  Let $A$
be a very ample line bundle over $M$ such that $A-K_M$ is ample. Let
$h_A$ be a smooth metric of $A$ whose curvature form $\omega_A$ is
positive on $M$.  We assume that $A$ is chosen to be sufficiently
ample so that for each point $P\in M$ the proper transform of $A$ in
the manifold obtained from $M$ by blowing up $P$ is still very
ample. This technical assumption will enable us to choose a generic
element of $\Gamma\left(M,\,A\right)$ vanishing at $P_0$ which is
not a zero-divisor of a prescribed coherent ideal sheaf.

\medbreak Let $p$ and $k$ be positive integers and we will impose
more conditions on $p$ and $k$ later. Take $P_0\in M$ and we will
also impose more conditions on $P_0$ later. Let $s_1$ be a generic
element of $\Gamma\left(M,\,A\right)$ vanishing at $P_0$ so that the
short exact sequence
$$
\displaylines{\qquad 0\to{\mathcal
I}_{p\varphi}\left(pL+kA\right)\stackrel{\theta_{s_1}}{\longrightarrow}{\mathcal
I}_{p\varphi}\left(pL+\left(k+1\right)A\right)\hfill\cr\hfill\to\left({\mathcal
I}_{p\varphi}\left/s_1{\mathcal
I}_{p\varphi}\right.\right)\left(pL+\left(k+1\right)A\right)\to
0\qquad\cr}
$$
is exact, where $\theta_{s_1}$ is defined by multiplication by
$s_1$. For this step we have to make sure that the maximum ideal
${\mathfrak m}_{M,P_0}$ of $M$ at $P_0$ is not an associated prime
ideal in the primary decomposition of the stalk of the ideal sheaf
${\mathcal I}_{p\varphi}$ at $P_0$. This means that for each $p$ we
have to impose the condition that $P_0$ does not belong to some
finite subset $Z_0$ of $M$.  Let $M_1$ be the zero-set of $s_1$ and
$$
{\mathcal O}_{M_1}=\left({\mathcal O}_M\left/s_1{\mathcal
O}_M\right.\right)|_{M_1},
$$
which we can assume to be regular with ideal sheaf equal to
$s_1{\mathcal O}_M$ because $s_1$ is generic element of
$\Gamma\left(M,A\right)$ vanishing at $P_0$.  By choosing $s_1$
generically we can also assume that ${\mathcal
I}_{\left(p\varphi|_{M_1}\right)}={\mathcal
I}_{p\varphi}\left/s_1{\mathcal I}_{p\varphi}\right.$. We use
$\chi\left(\cdot,\,\cdot\right)$ to denote the arithmetic genus
which means
$$
\chi\left(\cdot,\,\cdot\right)=\sum_{\nu=0}^\infty(-1)^\nu\dim_{\mathbb
C}H^\nu\left(\cdot,\,\cdot\right).
$$
From the long cohomology exact sequence of the above short exact
sequence we obtain
$$
\displaylines{\chi\left(M,\,{\mathcal
I}_{p\varphi}\left(pL+\left(k+1\right)A\right)\right)=\cr
\chi\left(M,\,{\mathcal
I}_{p\varphi}\left(pL+kA\right)\right)+\chi\left(M_1,\,{\mathcal
I}_{\left(p\varphi|_{M_1}\right)}\left(pL+\left(k+1\right)A\right)|_{M_1}\right).\cr}
$$
Since $A-K_M$ is ample and $2A-K_{M_1}=A-K_M$ is also ample, when we
assume $k\geq 1$, by the theorem of Kawamata-Viehweg-Nadel
$$
\displaylines{H^\nu\left(M,{\mathcal
I}_{p\varphi}\left(pL+kA\right)\right)=0\quad{\rm for\ }\nu\geq
1,\cr H^\nu\left(M_1,{\mathcal
I}_{\left(p\varphi|_{M_1}\right)}\left(\left(pL+\left(k+1\right)A\right)|_{M_1}\right)\right)=0\quad{\rm
for\ }\nu\geq 1.\cr}
$$
so that
$$
\displaylines{\Gamma\left(M,\,{\mathcal
I}_{p\varphi}\left(pL+\left(k+1\right)A\right)\right)=\cr
\Gamma\left(M,\,{\mathcal
I}_{p\varphi}\left(pL+kA\right)\right)+\Gamma\left(M_1,\,{\mathcal
I}_{\left(p\varphi|_{M_1}\right)}\left(\left(pL+\left(k+1\right)A\right)|_{M_1}\right)\right)\cr}
$$

\bigbreak\noindent(6.7) {\it Slicing by Ample Divisors Down to a
Curve.} Instead of one single element $s\in\Gamma\left(M,A\right)$,
we can choose generically
$$
s_1, \cdots, s_{n-1} \in\Gamma\left(M,\,A\right)
$$
all vanishing at $P_0$ so that inductively for $1\leq\nu\leq n-1$
the common zero-set $M_\nu$ of $s_1,\cdots,s_\nu$ with the structure
sheaf
$$
{\mathcal O}_{M_\nu}:=\left({\mathcal O}_M\left/\sum_{j=1}^\nu
s_j{\mathcal O}_M\right)\right|_{M_\nu}
$$
is regular and we end up with the inequality
$$
\displaylines{\dim_{\mathbb C}\Gamma\left(M,\,{\mathcal
I}_{p\varphi}\left(pL+\left(k+n-1\right)A\right)\right)\cr\geq
\dim_{\mathbb C}\Gamma\left(M_{n-1},\,{\mathcal
I}_{\left(p\varphi|M_{n-1}\right)}\left(\left(pL+\left(k+n-1\right)A\right)|_{M_{n-1}}\right)\right).\cr}
$$
For this step we have to exclude $P_0$ from a subvariety $Z_{n-2}$
of dimension $\leq n-2$ in $M$, because we have to exclude a finite
set in each $M_1$ which would come together as the hypersurface
$M_1$ varies to form a subvariety $Z_1$ of dimension $\leq 1$ in $M$
(as one can argue with the quotients of coherent ideal sheaves by
non zero-divisors and with the primary decompositions for coherent
ideal sheaves). Likewise we have a subvariety $Z_k$ of dimension
$\leq k$ in $M$ so that $Z_k$ intersects $M_k$ in a finite number of
points and finally we have end up with a subvariety $Z_{n-2}$ of
dimension $\leq n-2$ in $M$ which intersects $M_{n-2}$ in a finite
number of points and we impose the condition that $P_0$ does not
belong to $Z_{n-2}$.

\medbreak Since $M_{n-1}$ is a curve, all coherent ideal sheaves on
it are principal and are locally free and they come from holomorphic
line bundles. We can choose $s_1,\cdots,s_{n-1}$ so generically that
$M_{n-1}$ is disjoint from $Z$.  For this step we need to make sure
that $P_0$ does not belong to $Z$.

\medbreak We remark that actually we can also accommodate the case
of $P_0\in Z$ as long as the number of times $M_{n-1}$ intersects
$Z$ is far less than $\int_{M_{n-1}}R$.  For the proof of this
non-vanishing theorem we are not interested in the problem of
accommodating the case $P_0\in Z$.  We simply agree to choose $P_0$
outside of $Z$.  We add this remark to indicate how our general
non-vanishing theorem can be strengthened.

\medbreak We would like to remark also that this particular step of
slicing by $n-1$ ample divisors to get down to a curve roughly
corresponds to the step in Shokurov's proof of his non-vanishing
theorem [Shokurov 1985] where he takes the product of his
numerically effective divisor in his $n$-dimensional manifold with
the $(n-1)$-th power of a numerically effective big line bundle.

\bigbreak\noindent(6.8) {\it Application of the Theorem of
Riemann-Roch to a Curve, Comparing Contributions from the Curvature
Current and the Multiplier Ideal Sheaves, and the Unbounded Lower
Bound Condition.} Let $c$ be the nonnegative number
$$\int_{M_{n-1}}R=\int_M
R\wedge\left(\omega_A\right)^{n-1}.$$ Then
$$
\displaylines{(6.8.1)\qquad\qquad\dim_{\mathbb
C}\Gamma\left(M,\,{\mathcal
I}_{p\varphi}\left(pL+\left(k+n-1\right)A\right)\right)\hfill\cr\geq
\dim_{\mathbb C}\Gamma\left(M_{n-1},\,{\mathcal
I}_{\left(\left.p\varphi\right|M_{n-1}\right)}\left(\left(pL+\left(k+n-1\right)A\right)|_{M_{n-1}}\right)\right)\cr
\geq 1-{\rm genus}\left(M_{n-1}\right)+
\left(k+n-1\right)A^{n-1}M_{n-1}
\cr+\sum_{j=1}^J\left(p\tau_j-\left\lfloor
p\tau_j\right\rfloor\right) V_j\cdot A^{n-1}+p\int_{M_{n-1}}R,\cr}
$$
where the last identity is from the theorem of Riemann-Roch applied
to the regular curve $M_{n-1}$ and the locally free sheaf
$${\mathcal
I}_{\left(\left.p\varphi\right|M_{n-1}\right)}\left(\left(pL+\left(k+n-1\right)A\right)|_{M_{n-1}}\right)$$
on $M_{n-1}$. Let $\Xi_p$ be the right-hand side of $(6.8.1)$. We
introduce the following {\it unbounded-lower-bound condition}
$(6.8.2)$.
$$\lim_{\nu\to\infty}\Xi_{p_\nu}=\infty\quad{\rm for\ some\ subsequence\ }\left\{p_\nu\right\}_{\nu=1}^\infty.\leqno{(6.8.2)}
$$
For example, $(6.8.2)$ holds when $c>0$ (and we can choose
$\left\{p_\nu\right\}_\nu$ to be the same as
$\left\{p\right\}_{p\in{\mathbb N}}$).  Later we will show by using
diophantine approximation that $(6.8.2)$ also holds when $J=\infty$.
So we are left with the case when $R=0$ and $J<\infty$, which we
will handle later with Shokurov's technique of using the theorem of
Riemann-Roch to compare the arithmetic genus of a line bundle and
that of its twisting by a flat line bundle.  At this point let us
assume that $(6.8.2)$ holds and continue with the proof under such
an assumption.

\bigbreak\noindent(6.9) {\it Construction of Sections with Extra
Vanishing Order from Dimension Counting and Construction of Metrics
by Canceling Contributions from Ample Divisors by Using the General
Type Property.} For any $\ell\in{\mathbb N}$ the number of terms in
a polynomial of degree $\ell$ in $d$ variables is ${d+\ell\choose
d}$. Take any positive integer $q$. Take a positive integer $N$ and
we will impose more condition on $N$ later.  By Condition $(6.8.2)$,
there exists $p\in{\mathbb Z}$ such that
$$
\dim_{\mathbb C}\Gamma\left(M, \,{\mathcal
I}_{p\varphi}\left(pL+\left(k+n-1\right)A\right)\right)\geq
1+{n+N\left(k+n-1\right)q\choose n}
$$
and we can find some non identically zero element $s$ of
$$\Gamma\left(M, \,{\mathcal
I}_{p\varphi}\left(pL+\left(k+n-1\right)A\right)\right)$$ which
vanishes to order at least $N\left(k+n-1\right)q$ at $P_0$ so that
$s^{\frac{1}{N\left(k+n-1\right)}}$ is a multi-valued holomorphic
section of the ${\mathbb Q}$-line-bundle
$\frac{p}{N\left(k+n-1\right)}\,L+\frac{1}{N}\,A$ over $M$ which
vanishes to order at least $q$ at $P_0$.  We assume that $N$ is
chosen so large that the curvature current $\Theta_\chi$ dominates
$\frac{2}{N}\omega_A$.  Let $\hat p$ to be the round-up of
$\frac{p}{N\left(k+n-1\right)}$ and $\delta_p=\hat
p-\frac{p}{N\left(k+n-1\right)}$.  We introduce the metric
$$
e^{-\tilde\chi}:=\frac{e^{-\chi-\delta_p\varphi}}{\left(h_A\right)^{\frac{1}{N}}\left|s\right|^{\frac{2}{Nk}}}
$$
of $\left(p+p_0\right)L-K_M$ so that the multiplier ideal at $P_0$
is contained in $\left({\mathfrak
m}_{M,P_0}\right)^{{}^{\left\lfloor\frac{q}{n}\right\rfloor}}$.

\medbreak Since the metric $e^{-\chi}$ is assumed to be stable, for
any given $\hat q$ we can conclude by choosing $N$ and $q$
sufficiently large that the multiplier ideal sheaf of
$e^{-\tilde\chi}$ is contained in $\left({\mathfrak
m}_{M,P_0}\right)^{\hat q}{\mathcal I}_{\hat p\varphi+\chi}.$

\bigbreak\noindent(6.10) {\it Interpolation of Metrics to Get the
Analog of Minimal Center of Log Canonical Singularities.}  We now
use the techniques introduced for the proof of Fujita conjecture
type results.  We interpolate the two metrics $e^{-\hat
p\varphi-\chi}$ and $e^{-\tilde\chi}$ for the line bundle
$\left(\hat p+p_0\right)L-K_M$ on $M$ and use a slight
modification as described in (4.4.3) to get the analog of a
minimal center of log canonical singularities.  We end up with a
metric $e^{-\kappa}$ of $\left(\hat p+p_0\right)L-K_M$ with the
following properties.

\begin{itemize}\item[(i)] The curvature current of the metric $e^{-\kappa}$ of
$\left(\hat p+p_0\right)L-K_M$ dominates some smooth positive
$(1,1)$-form on $M$.
\item[(ii)] The multiplier ideal sheaf
${\mathcal I}_\kappa$ of the metric $e^{-\kappa}$ of $\left(\hat
p+p_0\right)L-K_M$ is contained in the multiplier ideal sheaf
${\mathcal I}_{\hat p\varphi+\chi}$ of the metric $e^{-\hat
p\varphi-\chi}$ of $\left(\hat p+p_0\right)L-K_M$.
\item[(iii)] The support of
${\mathcal I}_{\hat p\varphi+\chi}\left/{\mathcal I}_\kappa\right.$
is an irreducible subvariety $\check M$ of codimension at least one
in $M$ such that after replacing $M$ by the result of applying to
$M$ some finite number of successive monoidal transformations with
nonsingular centers we can assume without loss of generality that
$\check M$ is a nonsingular hypersurface in $M$.
\item[(iv)] There is a metric $e^{-\check\chi}$ of $\left(\hat p+p_0\right)L-K_{\check M}$
on $\check M$ whose curvature current dominates some smooth positive
$(1,1)$-form on $\check M$ such that the restriction of ${\mathcal
I}_{\hat p\varphi+\chi}\left/{\mathcal I}_\kappa\right.$ to $\check
M$ agrees with ${\mathcal I}_{\check\chi}$ on $\check
M$.
\item[(v)] If $\check M$ is not some $V_{j_0}$,
the metric $e^{-\check\chi}$ of $\left(\hat p+p_0\right)L-K_{\check
M}$ on $\check M$ and the metric $e^{-\varphi}$ of $L|_{\check M}$
on $\check M$ satisfy the conditions of Theorem (6.1) when the
metric $e^{-\check\chi}$ of $\left(\hat p+p_0\right)L-K_{\check M}$
on $\check M$ replaces the metric $e^{-\chi}$ of $p_0L-K_M$ on $M$
and $L|_{\check M}$ on $\check M$ replaces $L$ on $M$ with the
restriction of the metric $e^{-\varphi}$.  If $\check M$ is some
$V_{j_0}$, in the above description we have to replace $L|_{\check
M}$ by $\left(L-\tau_{j_0}\check M\right)|_{\check M}$ with the
corresponding natural modification of the metric.
\end{itemize}

\bigbreak\noindent(6.11) {\it Use of Induction Hypothesis for One
Dimension Less.}  Let us assume that $\check M$ is not some
$V_{j_0}$ and proceed, otherwise we simply have to make a natural
modification as indicated above.  Since the dimension of $\check M$
is less than $n$, by induction hypothesis, outside some subvariety
of codimension at least two in $\check M$, for $p$ sufficiently
large the sheaf ${\mathcal
I}_{p\varphi+\check\chi}\left(\left(p+\hat p+p_0\right)L\right)$ is
generated by its global sections on $\check M$.  We are going to
extend such sections to all of $M$ in the standard way by using a
short exact sequence and the vanishing theorem of
Kawamata-Viehweg-Nadel as follows. From the short exact sequence
$$ \displaylines{\qquad 0\to{\mathcal I}_{p\varphi+\kappa}\left(\left(p+\hat
p\right)L\right)\to{\mathcal I}_{\left(p+\hat
p\right)\varphi+\chi}\left(\left(p+\hat
p\right)L\right)\hfill\cr\hfill\to \left({\mathcal I}_{\left(p+\hat
p\right)\varphi+\chi}\left/{\mathcal
I}_{p\varphi+\kappa}\right.\right)\left(\left(p+\hat
p\right)L\right)\to 0\qquad\cr}
$$
and the vanishing of
$$
H^1\left(M,\,{\mathcal I}_{p\varphi+\kappa}\left(\left(p+\hat
p\right)L\right)\right)
$$
it follows that the map
$$
\displaylines{\Gamma\left(M,\,{\mathcal I}_{\left(p+\hat
p\right)\varphi+\chi}\left(\left(p+\hat
p\right)L\right)\right)\to\Gamma\left(M,\,\left({\mathcal
I}_{\left(p+\hat p\right)\varphi+\chi}\left/{\mathcal
I}_{p\varphi+\kappa}\right.\right)\left(\left(p+\hat
p\right)L\right)\right)\hfill\cr\hfill=\Gamma\left(\check
M,\,{\mathcal I}_{p\varphi+\check\chi}\left(\left(p+\hat
p\right)L\right)\right)\cr}
$$
is surjective.

\bigbreak We now go back to Condition $(6.8.2)$.  A main tool is the
results given in \S5 on diophantine approximation which are derived
from Kronecker's theorem on diophantine approximation.

\bigbreak\noindent(6.12) {\it Different Cases for the
Unbounded-Lower-Bound Condition.}  For the unbounded-lower-bound
condition $(6.8.2)$ we earlier introduced a dichotomy into two cases
and again for one case a further dichotomy.  It is the same as to
differentiate among three cases. The first and the easiest case is
$R\not=0$ and we have seen that it guarantees the
unbounded-lower-bound condition $(6.8.2)$. When $R=0$, the other two
cases are (i) $J=\infty$ and (ii) $J<\infty$.

\bigbreak\noindent(6.12.1) {\it Case of Infinite Number of
Irreducible Hypersurface Lelong Sets.} We now look at the case
$J=\infty$ and we are going to verify the unbounded-lower-bound
condition $(6.8.2)$ for this case.  For the verification we
differentiate among three possibilities.  Recall the following
decomposition of the curvature current $\Theta_\varphi$ of the
metric $e^{-\varphi}$ of the line bundle $L$ as a closed positive
$(1,1)$-current.
$$
\Theta_\varphi=\sum_{j=1}^\infty\tau_j\left[V_j\right]+R.
$$
We have the following three cases.
\begin{itemize}
\item[(i)] for any positive number $N$ one can find $N$ irrational
elements $\left\{\tau_{j_\nu}\right\}_{\nu=1}^N$ in
$\left\{\tau_j\right\}_{j=1}^\infty$ which are ${\mathbb
Q}$-linearly independent.
\item[(ii)] there is some positive integer $N_0$ such that any subsets
of $N_0+1$ elements of $\left\{\tau_j\right\}_{j=1}^\infty$ are
${\mathbb Q}$-linearly dependent, but at least one element of
$\left\{\tau_j\right\}_{j=1}^\infty$ is irrational.
\item[(iii)] every element of
$\left\{\tau_j\right\}_{j=1}^\infty$ is rational.
\end{itemize}
For Case (i) we can simply apply Corollary (5.3) above to get the
unbounded-lower-bound condition $(6.8.2)$.

\medbreak For Case (ii) we can simply apply Lemma (5.5) above to get
the unbounded-lower-bound condition $(6.8.2)$.

\medbreak Case (iii) is more complicated.  We have to handle it by
modifying somewhat the above main argument of (6.8).

\bigbreak\noindent(6.12.2) {\it Modification of Main Argument for
the Case of All Lelong Numbers of Hypersurface Lelong Sets Being
Rational.}  For any prescribed positive integer $N$, let $p$ be the
least common multiple for the denominators of the rational numbers
$\gamma_1,\cdots,\gamma_N$ so that each $p\gamma_j$ is a positive
integer for $1\leq j\leq N$.

\medbreak We choose $\eta_0>0$ sufficiently small so that
$\frac{1}{2}\,A+\sum_{j=1}^N\eta p\,V_j$ is an ample ${\mathbb
Q}$-line-bundle over $X$ with positively curved smooth metric
$h_\eta$ for any positive rational number $\eta\leq\eta_0$. Now
choose a positive integer $\tau $ so large that $\frac{1}{\tau
}<\eta_0$ and $\tau
>2p\tau_j$ for every $1\leq j\leq N$.  Then from $\frac{\tau
-1}{\tau }\,p\tau_j<p\tau_j$ and
$$\frac{\tau -1}{\tau }\,p\tau_j-\left(p\tau_j-1\right)=1-\frac{1}{\tau }\,p\tau_j>\frac{1}{2}
\quad{\rm for\ \ }1\leq j\leq N$$ it follows that
$$\sum_{j=1}^N\left(\frac{\tau -1}{\tau }\,p\tau_j-
\left\lfloor \frac{\tau -1}{\tau }\,p\tau_j\right\rfloor\right) \geq
\sum_{j=1}^N\frac{1}{2}\geq\frac{N}{2}.\leqno{(6.12.2.1)}$$  We have
a metric $h_{\frac{1}{\tau }}$ for the ${\mathbb Q}$-line-bundle
$\frac{1}{2}\,A+\sum_{j=1}^N\frac{p}{\tau }\,V_j$ (which is the
metric $h_\eta$ for the ${\mathbb Q}$-line-bundle
$\frac{1}{2}\,A+\sum_{j=1}^N\eta\,p\,V_j$ when $\eta=\frac{1}{\tau
}$).  We now define the metric
$$
e^{-\tilde\varphi_p}=e^{-\frac{\left(\tau -1\right)p}{\tau
}\varphi}\,h_{\frac{1}{\tau }}\left(h_A\right)^{\frac{1}{2}}
$$
of $pL+A$.

\medbreak At this point we introduce our modification of the above
main argument of (6.8).  Instead of using the metric
$e^{-p\varphi}h_A$ for $pL+A$, we use the metric
$e^{-\tilde\varphi_p}$ of $pL+A$ and we have the following two
conclusions, the second of which comes from $(6.13.2.1)$.
\begin{itemize}
\item[(a)] The curvature current of the metric
$e^{-\tilde\varphi_p}$ of $pL+A$ dominates that of the metric
$\left(h_A\right)^{\frac{1}{2}}$ of $\frac{1}{2}\,A$.
\item[(b)] For a generic curve $C$ in $M$ the Chern class of the
line bundle associated to ${\mathcal
I}_{\tilde\varphi_p|_C}\left(pL+A\right)$ is at least $\left\lfloor
\frac{N}{2}\right\rfloor$.
\end{itemize}

We repeat the main argument of (6.8) with the metric
$e^{-\tilde\varphi_p}$ of $pL+A$ instead of the metric
$e^{-p\varphi}h_A$ of $pL+A$. Let $s_1$ be a generic element of
$\Gamma\left(M,\,A\right)$ vanishing at $P_0$ so that the short
exact sequence
$$
0\to{\mathcal
I}_{\tilde\varphi_p}\left(pL+kA\right)\stackrel{\theta_{s_1}}{\longrightarrow}{\mathcal
I}_{\tilde\varphi_p}\left(pL+\left(k+1\right)A\right)\to{\mathcal
I}_{\tilde\varphi_p}\left(pL+\left(k+1\right)A\right)\left/s_1\,{\mathcal
I}_{\tilde\varphi_p}\right.\to 0
$$
is exact, where $\theta_s$ is defined by multiplication by $s_1$.
Let $M_1$ be the zero-set of $s_1$ and we have
$$
\displaylines{\chi\left(M,\,{\mathcal
I}_{\tilde\varphi_p}\left(pL+\left(k+1\right)A\right)\right)=\cr
\chi\left(M,\,{\mathcal
I}_{\tilde\varphi_p}\left(pL+kA\right)\right)+\chi\left(M_1,\,{\mathcal
I}_{\left(\tilde\varphi_p|Y\right)}\left(pL+\left(k+1\right)A\right)|_{M_1}\right).\cr}
$$
Again, instead of one single element $s\in\Gamma\left(M,A\right)$,
we can use generically
$$
s_1, \cdots, s_{n-1} \in\Gamma\left(M,\,A\right)
$$
all vanishing at $P_0$ so that inductively for $1\leq\nu\leq n-1$ we
have $M_\nu$ equal to the common zero-set of $s_1,\cdots,s_\nu$ and
we end up with the inequality
$$
\displaylines{\dim_{\mathbb C}\Gamma\left(M,\,{\mathcal
I}_{\tilde\varphi_p}\left(pL+\left(k+n-1\right)A\right)\right)\cr\geq
\dim_{\mathbb C}\Gamma\left(M_{n-1},\,{\mathcal
I}_{\left(\tilde\varphi_p|M_{n-1}\right)}\left(\left(pL+\left(k+n-1\right)A\right)|_{M_{n-1}}\right)\right).\cr}
$$
Since $M_{n-1}$ is a curve, all coherent ideal sheaves on it are
principal and are locally free and they come from holomorphic line
bundles.  We can choose $s_1,\cdots,s_{n-1}$ so generically that
$M_{n-1}$ is disjoint from $Z$. Then
$$
\displaylines{\dim_{\mathbb C}\Gamma\left(M,\,{\mathcal
I}_{\tilde\varphi_p}\left(pL+\left(k+n-1\right)A\right)\right)\cr\geq
\dim_{\mathbb C}\Gamma\left(M_{n-1},\,{\mathcal
I}_{\left(\left.\tilde\varphi_p\right|M_{n-1}\right)}\left(\left(pL+\left(k+n-1\right)A\right)|_{M_{n-1}}\right)\right)\cr
\geq 1-{\rm genus}\left(M_{n-1}\right)+
\left(k+n-3\right)A^{n-1}M_{n-1}+\frac{N}{2},\cr}
$$
because of Property (b) of the metric $e^{-\tilde\varphi_p}$ of
$pL+A$. Thus
$$
\dim_{\mathbb C}\Gamma\left(M,\,{\mathcal
I}_{\tilde\varphi_p}\left(pL+\left(k+n-1\right)A\right)\right)\geq
\frac{N}{3}$$ for $N$ sufficiently large (relative to the genus of
the curve $M_{n-1}$).  From this point on we just follow the rest of
the main argument.

\bigbreak\noindent(6.13) {\it Comparing the Use of the Theorem of
Riemann-Roch for a Line Bundle and Its Twisting by a Flat Line
Bundle.} We now discuss the remaining case of $R\not=0$ and
$J<\infty$.  For this case we cannot continue the next inductive
step of reducing the dimension of the new subspace defined by the
multiplier ideal sheaf. However, for this case there is no need to
continue with the next step, because the purpose of inductive step
is to produce a holomorphic section at the end when the dimension of
the subspace defined by the multiplier ideal sheaf becomes isolated.
For this case we can construct a holomorphic section directly from
the curvature current without continuing with the inductive process.
The construction comes from the technique of modification by a flat
bundle and comparing the use of the theorem of Riemann-Roch for the
original line bundle and its twisting by the flat bundle.

\medbreak For this case our curvature current is a finite ${\mathbb
R}$-linear combination of integration over (the regular part of)
irreducible hypersurfaces.  Our decomposition of the curvature
current $\Theta_\varphi$ of the metric $e^{-\varphi}$ of the line
bundle $L$ as a closed positive $(1,1)$-current becomes
$$
\Theta_\varphi=\sum_{j=1}^J\tau_j\left[V_j\right]
$$
with $J<\infty$ and $\tau_j\in{\mathbb R}$ and $\tau_j>0$.  If
$\tau_1,\cdots,\tau_J$ are all rational, we can simply introduce the
${\mathbb Q}$-line-bundle $\sum_{j=1}^J\tau_j\left[V_j\right]$ over
$X$ whose curvature current is the same as that of $L$, implying
that ${\mathbb Q}$-line-bundle $\sum_{j=1}^J\tau_j\left[V_j\right]$
is equal to the tensor product of some flat ${\mathbb
Q}$-line-bundle with $L$.  This flat ${\mathbb Q}$-line-bundle would
then be used in our twisting of $L$ and in comparing the application
of the theorem of Riemann-Roch to a high multiple of $L$ and its
twisting by this flat ${\mathbb Q}$-line-bundle (the high multiple
being used to make sure that such a high multiple of the ${\mathbb
Q}$-line-bundle would make it a usual holomorphic line bundle). The
strict positive lower bound for the curvature current of the metric
$e^{-\chi}$ of $p_0L-K_X$ makes the application of the theorem of
Kawamata-Viehweg-Nadel possible so that the dimension of the section
module is the same as the full arithmetic genus.  Now we have to
deal with the situation when some of the real numbers
$\tau_1,\cdots,\tau_J$ are not rational.

\medbreak Our strategy is to use the fact that $L$ is a usual
holomorphic line bundle.  If the classes in $H^{1,1}(M)\cap
H^2\left(M,{\mathbb R}\right)$ defined by $V_1,\cdots,V_J$ are all
${\mathbb R}$-linearly independent, then the real numbers
$\tau_1,\cdots,\tau_J$ must be all rational.  Hence there is some
${\mathbb R}$-linearly dependence among the classes in
$H^{1,1}(M)\cap H^2\left(M,{\mathbb R}\right)$ defined by
$V_1,\cdots,V_J$.  We use such ${\mathbb R}$-linearly dependency
relations to change each $\tau_1,\cdots,\tau_J$ slightly without
changing the class in $H^{1,1}(M)\cap H^2\left(M,{\mathbb R}\right)$
defined by $$ \Theta_\varphi=\sum_{j=1}^J\tau_j\left[V_j\right].
$$ This changes are so slight that the new numbers which replace
$\tau_1,\cdots,\tau_J$ are still positive.  Here is the detailed
implementation of this strategy.

\bigbreak\noindent(6.13.1) {\it Implementation of Strategy of
Twisting by a Flat Bundle.} We denote by ${\mathfrak V}$ the
${\mathbb Q}$-linear subspace ${\mathfrak V}$ in $H^{1,1}(M)\cap
H^2\left(M,{\mathbb R}\right)$ spanned by $V_1,\cdots,V_J$.  Let
${\rm Class}\left(V_j\right)$ be the class in $H^{1,1}(M)\cap
H^2\left(M,{\mathbb R}\right)$ defined by $V_j$ for $1\leq j\leq J$
and let ${\rm Class}\left(\Theta_\varphi\right)$ be the class in
$H^{1,1}(M)\cap H^2\left(M,{\mathbb R}\right)$ defined by
$\Theta_j$.  After relabeling we can assume without loss of
generality that ${\rm Class}\left(V_1\right),\cdots,{\rm
Class}\left(V_{J_0}\right)$ form a ${\mathbb R}$-basis of
${\mathfrak V}$ for some $1\leq J_0\leq J$. If $J_0=J$, since the
class defined $L$ in $H^{1,1}(M)\cap H^2\left(M,{\mathbb R}\right)$
lies in $H^{1,1}(M)\cap H^2\left(M,{\mathbb Q}\right)$, it follows
that $\tau_1,\cdots,\tau_J$ must be all rational.  In this case we
set $\rho_j$ for $1\leq j\leq J$.

\medbreak Now assume that $J_0<J$.  Inside ${\mathfrak V}$ we can
write
$$
{\rm Class}\left(V_j\right)=\sum_{k=1}^{J_0}a_{j,k}\,{\rm
Class}\left(V_k\right)
$$
for some $a_{j,k}\in{\mathbb Q}\ $ ($J_0+1\leq j\leq J$, $1\leq
k\leq J_0$) which may be positive or zero or negative.  Choose
$p\in{\mathbb N}$ such that\begin{itemize}
\item[(i)] the denominator of each $a_{j,k}$ is a factor of $p$
for $J_0+1\leq j\leq J$ and $1\leq k\leq J_0$,
\item[(ii)] the denominator of every rational member of
$\tau_1,\cdots,\tau_J$ is a factor of $p$,
\item[(iii)] if $a_1,\cdot,a_{J_0}\in{\mathbb Q}$ and $\sum_{j=1}^{J_0}a_j\,{\rm Class}\left(V_j\right)$
belongs to $H^{1,1}(M)\cap H^2\left(M,{\mathbb Z}\right)$, then
$pa_j\in{\mathbb Z}$ for $1\leq j\leq J_0$.\end{itemize}

\medbreak Choose $C>1$ such that $C>p\,\tau_j$ for $1\leq j\leq J_0$
and $C>\left|a_{jk}\right|$ for $J_0+1\leq j\leq J$ and $1\leq k\leq
J_0$.  Let $s_{V_j}$ be the canonical section of the line bundle
$V_j$ over $M$ for $1\leq j\leq J$.   Choose $\delta>0$ such that
the local function
$$\frac{1}{\left|s_{V_j}\right|^{2\delta}}$$ is locally integrable on
$M$ for $1\leq j\leq J$.  We choose a number
$$0<\eta<\frac{\delta}{\left(J-J_0\right)C^2p}\leqno{(6.13.1.1)}$$ which is less than
the minimum of $\tau_1,\cdots,\tau_J$ such that the total mass of
the closed positive $(1,1)$-current
$\left(J-J_0\right)C^2p\,\eta\sum_{j=1}^J\left[V_j\right]$ is less
than the positive number in the definition (6.1) for the stability
of the metric $e^{-\chi}$. Choose $0<\eta<\frac{1}{2}$ and we will
impose more conditions on $\eta$ later.  If none of
$\tau_1,\cdots,\tau_{J_1}$ is rational, we apply Kronecker's theorem
(5.1) to the ${\mathbb Q}$-linearly independent set
$1,\tau_1,\cdots,\tau_{J_1}$ to find $t>1$ such that
$\left|t-1-x_0\right|<\eta$ and $\left|t\tau_j-x_j\right|<\eta$ for
$1\leq j\leq J_0$ for some integers $x_0,\cdots,x_{J_0}$.  If one of
$\tau_1,\cdots,\tau_{J_1}$ is rational and we can assume without
loss of generality that $\tau_1$ is rational, we apply Kronecker's
theorem (5.1) to the ${\mathbb Q}$-linearly independent set
$p\,\tau_1,\cdots,p\,\tau_{J_1}$ to find $t>1$ such that
$\left|t-p\,\tau_j-x_j\right|<\eta$ for $1\leq j\leq J_0$ for some
integers $x_1,\cdots,x_{J_0}$.  Let $q$ be the integer closest to
$t$.  Then $\left|t-q\right|<\eta$.  Moreover, by our choice of $t$
we have $\left|qp\,\tau_j-\rho_j\right|<C\eta$ for every $1\leq
j\leq J_0$ and for some integers $\rho_j$. Let
$$
\rho_k=qp\,\tau_k+\sum_{j=J_0+1}^J\left(qp\,\tau_j-\rho_j\right)a_{j,k}
$$
for $1\leq j\leq J_0$.  Then
$$
\left|qp\,\tau_k-\rho_k\right|<\left(J-J_0\right)C^2\eta\quad{\rm
for\ \ }1\leq j\leq J_0.$$  Inside ${\mathfrak V}$ we have
$$qp\,{\rm
Class}\left(\Theta_\varphi\right)-\sum_{j=J_0+1}^J\rho_j\,{\rm
Class}\left(V_j\right)=\sum_{k=1}^{J_0}\rho_k\,{\rm
Class}\left(V_k\right).\leqno{(6.13.1.2)}
$$
Since $$qp\,{\rm
Class}\left(\Theta_\varphi\right)-\sum_{j=J_0+1}^J\rho_j\,{\rm
Class}\left(V_j\right)$$ belongs to $H^{1,1}(M)\cap
H^2\left(M,{\mathbb Z}\right)$, it follows from the choice of $p$
and $(6.13.1.2)$ that $p\,\rho_k$ is rational for $1\leq k\leq J_0$.
We now have
$$
qp^2{\rm Class}\left(\Theta_\varphi\right)=\sum_{j=1}^J
p\,\rho_j\,{\rm Class}\left(V_j\right)
$$
with each $p\,\rho_j\in{\mathbb Z}$ and
$$\left|p\,\rho_j-qp^2\tau_j\right|<\left(J-J_0\right)C^2p\eta\quad{\rm
for\ \ }1\leq j\leq J.\leqno{(6.13.1.3)}$$ Let $F$ be the flat
holomorphic line bundle $L-\sum_{j=1}^J p\,\rho_j V_j$ on $M$.
Consider the metric
$$
e^{-\tilde\psi}=\frac{1}{\left|\prod_{j=1}^J\left(s_{V_j}\right)^{p\,\rho_j}\right|^{\frac{2}{qp^2}}}
$$
of $L$.  By the theorem of Kawamata-Viehweg-Nadel we have
$$
H^\nu\left(M,\,{\mathcal
I}_{\left(qp^2-p_0\right)\tilde\psi+\chi}\left(qp^2L\right)\right)=0
$$
and $$H^\nu\left(M,\,{\mathcal
I}_{\left(qp^2-p_0\right)\tilde\psi+\chi}\left(qp^2\left(L+F\right)\right)\right)=0
$$
for $\nu\geq 1$.  By using the theorem of Riemann-Roch and the fact
that $F$ is flat, we conclude that
$$
\dim_{\mathbb C}\Gamma\left(M,\,{\mathcal
I}_{\left(qp^2-p_0\right)\tilde\psi+\chi}\left(qp^2L\right)\right)=
\dim_{\mathbb C}\Gamma\left(M,\,{\mathcal
I}_{\left(qp^2-p_0\right)\tilde\psi+\chi}\left(qp^2\left(L+F\right)\right)\right).
$$
Let $s=\prod_{j=1}^J\left(s_{V_j}\right)^{p\,\rho_j}$.  Then $s$ is
a non identically zero holomorphic section of $\sum_{j=1}^J p\rho_j
V_j=qp^2L+F$ over $M$.  From $(6.13.1.1)$ and $(6.13.1.3)$ it
follows that $s$ locally belongs to the multiplier ideal sheaf
${\mathcal I}_{\left(qp^2-p_0\right)\tilde\psi+\chi}$ and as a
consequence
$$0\not\equiv s\in\Gamma\left(M,\,{\mathcal
I}_{\left(qp^2-p_0\right)\tilde\psi+\chi}\left(qp^2\left(L+F\right)\right)\right)
$$
and $$\dim_{\mathbb C}\Gamma\left(M,\,{\mathcal
I}_{\left(qp^2-p_0\right)\tilde\psi+\chi}\left(qp^2\left(L+F\right)\right)\right)\geq
1
$$
and
$$
\dim_{\mathbb C}\Gamma\left(M,\,{\mathcal
I}_{\left(qp^2-p_0\right)\tilde\psi+\chi}\left(qp^2L\right)\right)\geq
1.
$$
This concludes the proof of our general non-vanishing theorem (6.2)
by induction.

\medbreak The next step for the proof of the finite generation of
the canonical ring for the case of the general type is to verify the
achievement of stable vanishing order at a generic point of a
subvariety of codimension at least two by a finite sum of the
absolute-value squares of multi-valued holomorphic pluricanonical
sections.  This we are going to do in the next section using
holomorphic families of Artinian subschemes.

\bigbreak

\bigbreak\noindent{\bf \S7. Holomorphic Family of Artinian
Subschemes and Achievement of Stable Vanishing Orders for the Case
of Higher Codimension.}

\bigbreak\noindent(7.1) {\it Graded Coherent Ideal Sheaves.}  To
study the problem of achieving the stable vanishing orders across a
higher codimensional subvariety at a generic point, we have to find
the analog for the following two procedures for the codimension one
case.

\bigbreak\noindent(i) When $Y$ is an irreducible codimension-one
Lelong set in $X$ for the metric $e^{-\varphi}=\frac{1}{\Phi}$ of
$K_X$ with generic Lelong number $\gamma$, we have to modify the
metric $e^{-\varphi}$ of $K_X$ to get the metric
$\frac{e^{-\varphi}}{\left|s_Y\right|^{2\gamma}}$ of $K_X-\gamma Y$
before restricting it to $Y$, where $s_Y$ is the canonical section
for the line bundle $Y$ on $X$. This is the procedure to remove the
vanishing along $Y$ by dividing by $\left(s_Y\right)^\gamma$.  The
problem is to find the corresponding procedure in the case of the
higher-codimensional subvariety $V$ to remove the vanishing along
the subvariety $V$.  Of course, the simplest way is to blow up to
replace the higher-codimensional subvariety $V$ by a hypersurface,
but the difficulty of using blow-ups is the inability to know
definitely that we can stop in a finite number of steps. Without
blowing up we can use descending induction on the dimension of the
higher-codimensional subvariety to finish in a finite number of
stops.  The procedure we are going to remove the vanishing along a
higher-codimensional subvariety $V$ involves graded coherent ideal
sheaves, which we will explain in details shortly.

\bigbreak\noindent(ii) The proof of the finite generation of the
canonical ring uses descending induction on the dimension of the
subvariety where the stable vanishing order is yet achieved by a
finite partial sum of $\Phi$.  For each subvariety $V$ in the
induction process we have to identify a subvariety $Z$ of
codimension at least $1$ in $V$ so that at points of $V-Z$ we can
prove that the stable vanishing order is achieved by a finite
partial sum of $\Phi$.  In the case of a hypersurface $Y$ with
generic Lelong number $\gamma$ of $\Theta_\varphi$ (with
$e^{-\varphi}=\frac{1}{\Phi}$) instead of a higher-codimensional
subvariety $V$, we use the decomposition of
$\Theta_\varphi-\gamma\left[Y\right]=\sum_{j=1}^J\tau_j V_j+R$ and
show that the stable vanishing order is achieved at some point
outside of any $V_j$ where the Lelong number of
$\Theta_\varphi-\gamma\left[Y\right]$ (and hence of the remainder
$R$ also) is zero.  In the case of a higher-codimensional $V$
instead of a hypersurface $Y$, the vanishing orders across $V$ in
the normal directions of $V$ are no longer given by a single number.
Instead the analog is an Artinian subscheme in the normal directions
of a generic point of $V$ (or more precisely an infinite sequence of
Artinian subschemes because, unlike the situation with nonnegative
numbers in the codimension one case, it makes no sense for us to
take roots of Artinian subschemes and then go to the limit). To
locate {\it a priori} a subvariety $Z$ of codimension at least $1$
in $V$ for proving the achievement of stable vanishing order by a
finite partial sum of $\Phi$ at points of $V-Z$, we will introduce
the procedure of detecting a jump in the structure of an Artinian
subscheme in a holomorphic family.  We will explain later this
procedure.

\bigbreak\noindent(7.2) {\it Definition of Graded Coherent Ideal
Sheaves, Their Finite Generation, Conductors, and Order
Functions.} By {\it a sequence of graded coherent ideal sheaves}
we mean a sequence of coherent ideal sheaves ${\mathcal
J}^{(\nu)}$ indexed by $\nu\in{\mathbb N}$ such that ${\mathcal
J}^{(\lambda)}{\mathcal J}^{(\nu)}\subset {\mathcal
J}^{(\lambda+\nu)}$ for $\lambda,\,\nu\in{\mathbb N}$.

\medbreak For a sequence of graded coherent ideal sheaves ${\mathcal
J}^{(\nu)}$ indexed by $\nu\in{\mathbb N}$ we define its order
function as
$$\sum_{\nu=1}^\infty\varepsilon_\nu\sum_{j=1}^{p_\nu}\left|g^{(\nu)}\right|^{\frac{2}{\nu}},
$$
where $g^{(\nu)}_1,\cdots,g^{(\nu)}_{p_\nu}$ are local generators of
the coherent ideal sheaf ${\mathcal J}^{(\nu)}$ and
$\left\{\varepsilon_\nu\right\}_{\nu\in{\mathbb N}}$ is a sequence
of positive numbers decreasing to $0$ so fast that the infinite
sequence locally converges.  An order function is only locally
defined and is not unique and we are interested only in its
vanishing order.  For our purpose we can use any function which is
comparable to an order function in the sense that locally one
function is some positive constant times the other.

\medbreak A sequence of graded coherent ideal sheaves
$\left\{{\mathcal J}^{(\nu)}\right\}_{\nu\in{\mathbb N}}$ is said to
be {\it finitely generated} on an open subset $U$ if there exists
some $\nu_0$ such that every ${\mathcal J}^{(\nu)}$ is generated on
$U$ by elements $\Gamma\left(U,{\mathcal J}^{(\lambda)}\right)$ for
$1\leq\lambda\leq\nu_0$ in the sense of a ring.

\medbreak Let $\left\{{\mathcal J}^{(\nu)}\right\}_{\nu\in{\mathbb
N}}$ and $\left\{{\mathcal K}^{(\nu)}\right\}_{\nu\in{\mathbb N}}$
be two sequences of graded coherent ideal sheaves.  By the {\it
conductor\;} from $\left\{{\mathcal
J}^{(\nu)}\right\}_{\nu\in{\mathbb N}}$ into $\left\{{\mathcal
K}^{(\nu)}\right\}_{\nu\in{\mathbb N}}$ we mean the sequence of
graded coherent ideal sheaves $\left\{{\mathcal
L}^{(\nu)}\right\}_{\nu\in{\mathbb N}}$ defined as follows.  A
holomorphic function germ $f$ at a point $P_0$ belongs to the stalk
$\left({\mathcal L}^{(\nu)}\right)_{P_0}$ of ${\mathcal L}^{(\nu)}$
at $P_0$ if and only if $f\left({\mathcal
J}^{(\lambda)}\right)_{P_0}\subset\left({\mathcal
K}^{(\lambda+\nu)}\right)_{P_0}$ for any $\lambda\in{\mathbb N}$.

\bigbreak\noindent(7.3) {\it Examples of Graded Coherent Ideal
Sheaves and the Motivations for Them.}  The most important example
which motivates the introduction of a sequence of graded coherent
ideal sheaves comes from our compact complex algebraic manifold
$X$ of general type.  Let ${\mathcal J}_{K_X}^{(\nu)}$ be the
coherent ideal sheaf on $X$ generated locally by elements of
$\Gamma\left(X,\nu K_X\right)$ for $\nu\in{\mathbb N}$.  The
collection $\left\{{\mathcal
J}_{K_X}^{(\nu)}\right\}_{\nu\in{\mathbb N}}$ so defined is a
sequence of graded coherent ideal sheaves on $X$. The function
$\Phi$ is an order function for this sequence of graded coherent
ideal sheaves on $X$.  The main purpose of introducing sequences
of graded coherent ideal sheaves is to study the achievement of
the vanishing orders of an order function by its finite partial
sum.

\medbreak Suppose $Y$ is a Lelong hypersurface with vanishing order
$\gamma$ for $\Phi$. We introduce the sequence of graded coherent
ideal sheaves $\left\{{\mathcal
J}_{Y,\gamma}^{(\nu)}\right\}_{\nu\in{\mathbb N}}$ so that
${\mathcal J}_{Y,\gamma}^{(\nu)}$ is the coherent ideal sheaf on $X$
locally generated by
$\left(s_Y\right)^{\left\lceil\nu\gamma\right\rceil}$. Let
$\left\{{\mathcal J}_{K_X,Y,\gamma}^{(\nu)}\right\}_{\nu\in{\mathbb
N}}$ be the conductor from $\left\{{\mathcal
J}_{Y,\gamma}^{(\nu)}\right\}_{\nu\in{\mathbb N}}$ into
$\left\{{\mathcal J}_{K_X}^{(\nu)}\right\}_{\nu\in{\mathbb N}}$.
This conductor $\left\{{\mathcal
J}_{K_X,Y,\gamma}^{(\nu)}\right\}_{\nu\in{\mathbb N}}$ is what is
left from removing the vanishing order along $Y$ from ${\mathcal
J}_{K_X}^{(\nu)}$.  The function
$\frac{\Phi}{\left|s_Y\right|^{2\gamma}}$ is comparable to an order
function of the conductor $\left\{{\mathcal
J}_{K_X,Y,\gamma}^{(\nu)}\right\}_{\nu\in{\mathbb N}}$.

\medbreak We now come to the higher codimensional case which is
actually the real reason to introduce the notion of a sequence of
graded coherent ideal sheaves and its order function.  Suppose the
precise achievement of vanishing order of $\Phi$ by one of its
finite partial sum is known at points outside a countable union of
subvarieties of codimension $\ell$ in $X$ (and hence outside a
single subvariety of codimension $\ell$ in $X$, because the
Noetherian argument applies when the precise achievement of stable
vanishing order is known for subvarieties of lower codimension). For
example, in the case of $\ell=2$ we know the stable vanishing order
$\gamma_j$ for a hypersurface $Y_j$ is known to be achieved and the
set of all such $Y_j$ is indexed by $1\leq j\leq N$.  Just as above
for the case of $N=1$, we can introduce the sequence of graded
coherent ideal sheaves $\left\{{\mathcal
J}_{Y_1,\gamma_1,\cdots,Y_N,\gamma_N}^{(\nu)}\right\}_{\nu\in{\mathbb
N}}$ so that ${\mathcal
J}_{Y_1,\gamma_1,\cdots,Y_N,\gamma_N}^{(\nu)}$ is the coherent ideal
sheaf on $X$ locally generated by
$\prod_{j=1}^N\left(s_{Y_j}\right)^{\left\lceil\nu\gamma_j\right\rceil}$.

\medbreak Under the assumption of the precise achievement of stable
vanishing order at points outside a subvariety of codimension $\ell$
in $X$, in the case of $\ell\geq 3$, the analog of $\left\{{\mathcal
J}_{Y_1,\gamma_1,\cdots,Y_N,\gamma_N}^{(\nu)}\right\}_{\nu\in{\mathbb
N}}$ will be a finitely generated sequence of graded coherent ideal
sheaves $\left\{{\mathcal J}_{\left({\rm codim}\geq
\ell-1\right)}^{(\nu)}\right\}_{\nu\in{\mathbb N}}$ which will
generate the precise vanishing orders of $\left\{{\mathcal
J}_{K_X}^{(\nu)}\right\}_{\nu\in{\mathbb N}}$ across subvarieties of
codimension $\geq\ell-1$ in $X$, both isolated and embedded.  Let
$\left\{{\mathcal J}_{\left(K_X,\,{\rm
codim}\geq\ell-1\right)}^{(\nu)}\right\}_{\nu\in{\mathbb N}}$ be the
conductor from $\left\{{\mathcal J}_{\left({\rm codim}\geq
\ell-1\right)}^{(\nu)}\right\}_{\nu\in{\mathbb N}}$ into
$\left\{{\mathcal J}_{K_X}^{(\nu)}\right\}_{\nu\in{\mathbb N}}$.
This conductor $\left\{{\mathcal J}_{\left(K_X,\,{\rm codim}\geq
\ell-1\right)}^{(\nu)}\right\}_{\nu\in{\mathbb N}}$ is what is left
from removing the vanishing order along all subvarieties of
codimension $\geq\ell-1$ from ${\mathcal J}_{K_X}^{(\nu)}$.  An
order function $\Phi_{K_X,\ell}$ of the conductor $\left\{{\mathcal
J}_{\left(K_X,\,{\rm
codim}\geq\ell-1\right)}^{(\nu)}\right\}_{\nu\in{\mathbb N}}$ is for
us to study the achievement of the vanishing orders of $\Phi$ by one
of its finite partial sums for the next step in the descending
induction process for codimension at least $\ell$.

\bigbreak\noindent(7.4) {\it Artinian Subschemes Defined by
Sequences of Graded Coherent Ideal Sheaves.}  We continue with the
notations introduced in (7.4.1).  Suppose there is an irreducible
subvariety $V$ of codimension $\ell$ which is a branch of the
zero-set of an order function $\Phi_{K_X,\ell}$ of the sequence
$\left\{{\mathcal J}_{\left(K_X,\,{\rm
codim}\geq\ell-1\right)}^{(\nu)}\right\}_{\nu\in{\mathbb N}}$ of
graded coherent ideal sheaves on $X$.   We are going to introduce
Artinian subschemes transversal to a generic point of $V$.  For
the subvariety $V$ of codimension $\ell$ in $X$ these Artinian
subschemes will play the r\^ole which the Lelong number plays for
the case of a hypersurface $Y$.  The Artinian subschemes used at
this point are unreduced complex subspaces of ${\mathbb C}^\ell$
supported at the origin.

\bigbreak\noindent(7.4.1) {\it Construction of Artinian
Subschemes.} The sequence of Artinian subschemes are introduced as
follows. What we consider now is only locally in $X$ at a generic
point $P_0$ of $V$. (Later when we have to go back to the global
situation, we will at that point indicate our return to global
consideration.)  Since the environment now is local, we can
consider the situation in which $X$ is replaced by some relatively
compact connected open neighborhood $U$ of the origin in ${\mathbb
C}^n$ and $V$ is replaced by
$\left\{\,z_1=\cdots=z_\ell=0\,\right\}\cap U$ with
$z_j\left(P_0\right)=0$ for $1\leq j\leq n$.  Let
$g^{(\nu)}_1,\cdots,g^{(\nu)}_{p_\nu}$ be holomorphic functions on
$U$ which generate the coherent ideal sheaf ${\mathcal
J}_{\left(K_X,\,{\rm codim}\geq\ell-1\right)}^{(\nu)}$ on $U$.

\medbreak We will later let $P_0$ vary in $V$.  So instead of fixing
$P_0$ in $V$ we consider a variable point $P$ in a neighborhood of
$P_0$ in $V$. Fix $P\in V$. Let $W_P$ be the $\ell$-dimensional
submanifold in ${\mathbb C}^n$ defined by
$$z_{\ell+1}=z_1\left(P\right),\cdots,z_n=z_n\left(P\right)$$
which we identify with ${\mathbb C}^\ell$.  Let $N_0\in{\mathbb N}$
be chosen so that the common zero-set of
$$\left\{g^{(\nu)}_j\right\}_{1\leq j\leq p_\nu,\,1\leq\nu\leq N_0}$$
in $U$ is $V$. (Here in some steps we may have to replace $U$ by a
slightly smaller neighborhood of the origin in $U$ and for
notational simplicity when there is no risk of confusion we will
just use the same letter $U$ to denote this slightly smaller
neighborhood of the origin in $U$.)

\medbreak We consider for every $N\geq N_0$ and $P\in V$ the
Artinian subscheme ${\mathcal A}_{V,N,P}$ in ${\mathbb C}^\ell$ with
coordinates $z_1,\cdots,z_\ell$ defined by the ideal ${\mathcal
K}_{V,N,P}$ generated by
$$\left\{\left.\left(g^{(\nu)}_j\right)^{\frac{N!}{\nu}}\right|_{W_P\cap U}\right\}_{1\leq j\leq p_\nu,\,1\leq\nu\leq N}.$$
We raise the holomorphic function $g^{(\nu)}_j$ to the
$\left(\frac{N!}{\nu}\right)$-th power in the definition of
${\mathcal K}_{V,N,P}$ and ${\mathcal A}_{V,N,P}$ solely for the
technical reason of changing fractional powers of holomorphic
functions to just holomorphic functions.  Of course, later when we
consider an order function for $\left\{{\mathcal
J}_{\left(K_X,\,{\rm
codim}\geq\ell-1\right)}^{(\nu)}\right\}_{\nu\in{\mathbb N}}$ we
will have to return to the original situation by taking the
$\left(N!\right)$-th root.

\bigbreak\noindent(7.5) {\it Purpose for Introducing Artinian
Subschemes.} Our purpose of introducing the sequence of Artinian
subschemes ${\mathcal A}_{V,N,P}$ is to locate the countable union
of proper subvarieties $Z_N$ in $V$ so that the structure of the
Artinian subscheme ${\mathcal A}_{V,N,P}$ varies continuously
without jump when $P\in V-Z_N$.  Let us first give here the
precise definition of a family of Artinian scheme ``varying
continuously without jump.''

\bigbreak\noindent(7.5.1) {\it Definition of the Continuous
Variation of an Artinian Subscheme Without Jump.} Let ${\frak M}$
be the moduli space of all Artinian subschemes on $W_{P_0}$
supported at the single point $P_0$. We can decompose ${\frak M}$
as the disjoint union of local submanifolds ${\frak M}_\nu$
(indexed by $\nu\in{\frak J}$) inside ${\frak M}$. Let ${\mathcal
A}_P$ be an Artinian subscheme on $W_P$ (supported at the single
point $P$) parametrized holomorphically by $P\in V$. We identify
naturally $W_P$ with $W_{P_0}$ so that we regard ${\mathcal A}_P$
naturally as an Artinian subscheme on $W_{P_0}$ supported at the
single point $P_0$. Let $\check Z$ be a subvariety of $V$.  We say
that the Artinian subscheme ${\mathcal A}_P$ {\it varies
continuously without jump} for $P\in V-\check Z$ if for some
$\nu_0\in{\frak J}$ each ${\mathcal A}_P$ is an element in ${\frak
M}_{\nu_0}$ for $P\in V-\check Z$.  Equivalently we also refer to
it as the variation of the Artinian subscheme ${\mathcal A}_P$
being {\it continuous without jump} for $P\in V-\check Z$.

\bigbreak The reason for considering the proper subvarieties $Z_N$
in $V$ is as follows.  In the case of a hypersurface $Y$ in $X$ we
consider the vanishing order of $\Phi$ across $Y$ and after we
remove the generic vanishing order $\gamma$ we end up with the
current $\Theta_{\log\Phi}-\gamma Y$ whose restriction to $Y$ has
the decomposition
$$
\Theta_{\log\Phi}-\gamma Y=\sum_{j=1}^J\tau_j\left[V_j\right]+R,
$$ where the Lelong number of $R$ is zero except outside a countable
union of subvarieties of codimension at least two in $Y$.

\medbreak For the hypersurface case, each Artinian subscheme
${\mathcal A}_{Y,N,P}$ (or ${\mathcal K}_{Y,N,P}$) is just
represented by a single positive integer $\gamma_N$ and the limit of
the sequence of Artinian subschemes $\left\{{\mathcal
A}_{Y,N,P}\right\}_{N\geq N_0}$ (each one as a single positive
integer $\gamma_N$) after normalization by the factor $\frac{1}{N!}$
is the number $\gamma$ at a generic point $P$ of $Y$ so that
$\lim_{N\to\infty}\frac{\gamma_N}{N!}=\gamma$. We are not as
interested in the number $\gamma$ as in the union $Z$ of
$\bigcup_{j=1}^J V_j$ and the Lelong sets of $R$, which is a
countable union of proper subvarieties of $Y$.  The significance of
this countable union $Z$ of proper subvarieties of $Y$ is that at a
point $P$ of $Y$ outside this countable union $Z$ we can get a
holomorphic section of $m\left(K_X-\gamma Y\right)$ for some
sufficiently large $m$ which is nonzero at $P$.  For the
hypersurface case this means that at a point of $Y-Z$ the infinite
sum $\Phi$ is comparable to one of its finite truncations in some
open neighborhood of that point.

\medbreak In the case of higher codimension what goes into replacing
the set $Z$ is the union $\bigcup_{N=N_0}^\infty Z_N$. Since
excluding a countable union (instead of a finite union) of proper
subvarieties is good enough for the purpose of showing the existence
of some points where the infinite sum $\Phi$ is comparable to one of
its finite truncations, there is no need for us to consider the
analog of the limit $\gamma$ and we can just take the entire
sequence of Artinian subschemes $\left\{{\mathcal
A}_{V,N,P}\right\}_{N\geq N_0}$ and work with the union
$\bigcup_{N=N_0}^\infty Z_N$ without taking any limit.  For the case
of the hypersurface $Y$, the key point is to be able to find
beforehand a countable union $Z$ of proper subvarieties in $Y$ which
contain the Lelong sets of $\Theta_{\log\Phi}-\gamma Y$ without any
knowledge of $\gamma$, because in the case of higher codimension
there is no single number $\gamma$ and it is difficult to make any
good sense of the limit, as $N\to\infty$, of a sequence of Artinian
subschemes $\left\{{\mathcal A}_{V,N,P}\right\}_{N\geq N_0}$, even
after normalization by the factor $\frac{1}{N!}$.

\bigbreak\noindent(7.5.2) {\it Continuous Variation of Relative
Positions of Artinian Subschemes Without Jump.}  In the case of
handling the precise achievement of vanishing orders of $\Phi$ at
a generic point of a subvariety $V$ of higher codimension
$\ell>1$, there is also another contribution which goes into the
{\it a priori} ``bad'' set in $V$ other than the union
$\bigcup_{N=N_0}^\infty Z_N$.  The reason is as follows.  When we
restrict the curvature current $\Theta_\varphi$ of the metric
$e^{-\varphi}=\frac{1}{\Phi}$ of $K_X$ to an irreducible
subvariety $S$ of codimension $\ell-1$ in $X$ which contains $V$
and when, after subtracting an appropriate multiple of $V$ from
$\Theta_\varphi\Big|_S$, we restrict it to $V$, the resulting
closed positive $(1,1)$-current
$\left(\Theta_\varphi\big|_S-\gamma_S\left[V\right]\right)\big|_V$
on $V$ has Lelong sets $Z_S$ in $V$ which may move inside $V$
according to the subvariety $S$ of $X$. These Lelong sets $Z_S$
are also ``bad'' sets which we would like to locate in an {\it a
priori} manner to make sure that they are contained in some fixed
countable union of subvarieties of codimension $\geq 1$ in $V$ as
$S$ goes through a family of such
$\left(\ell-1\right)$-codimensional subvarieties $S$ in $X$
parametrized by the space ${\mathbb P}_{\ell-1}$ of normal
directions of $V$. In order to accomplish this, we have to
consider a generalization of the continuous variation of Artinian
subschemes without jump. The generalization is the continuous
variation of relative positions of Artinian subschemes whose
precise meaning is as follows.

\medbreak Instead of one single Artinian subscheme ${\mathcal A}_P$
on $W_P$ supported at a single point parametrized holomorphically by
$P\in V$, we consider several Artinian subscheme ${\mathcal
A}^{(\lambda)}_P$ on $W_P$ supported at a single point parametrized
holomorphically by $P\in V$ for $0\leq\lambda\leq\rho$.  We also
assume that we have the following inclusion relation
$${\mathcal A}^{(\rho)}_P\subset{\mathcal A}^{(\rho-1)}_P\subset{\mathcal A}^{(2)}_P
\subset\cdots\subset{\mathcal A}^{(1)}_P\subset {\mathcal
A}^{(0)}_P\leqno{(7.5.2.1)_P}$$ Again we naturally identify $W_P$
with $W_{P_0}$ so that each ${\mathcal A}^{(\lambda)}_P$ can be
naturally regarded as an Artinian subscheme on $W_0$ supported at
the single point $P_0$.

\medbreak Let ${\frak M}^{(\rho)}$ be the moduli space of all
nested sequences $(7.5.2.1)_P$ of $\rho+1$ Artinian subschemes on
$W_{P_0}$, all supported at the single point $P_0$. We can
decompose ${\frak M}^{(\rho)}$ as the disjoint union of local
submanifolds ${\frak M}^{(\rho)}_\nu$ (indexed by $\nu\in{\frak
J}^{(\rho)}$) inside ${\frak M}^{(\rho)}$.  Let $\check
Z^{(\rho)}$ be a subvariety of $V$. We say that the nested
sequence $(7.5.2.1)_P$ of $\rho+1$ Artinian subschemes {\it varies
continuously without jump} for $P\in V-\check Z^{(\ell)}$ if for
some $\nu_0\in{\frak J}^{(\rho)}$ each nested sequence
$(7.5.2.1)_P$ is an element in ${\frak M}^{(\rho)}_{\nu_0}$ for
$P\in V-\check Z^{(\rho)}$. Equivalently we also refer to it as
the variation of the nested sequence $(7.5.2.1)_P$ of $\rho+1$
Artinian subschemes being {\it continuous without jump} for $P\in
V-\check Z^{(\rho)}$.

\medbreak The notion of continuous nested sequences of Artinian
subschemes without jump will be used in the following context.
Recall that we are still in the local situation where $X$ is
replaced by an open subset $U$ of ${\mathbb C}^n$.  We assume that
$h_1,\cdots,h_\rho$ are holomorphic functions on $U$ so that their
common zero-set contains $V$ and is of codimension $\rho$ in $U$.
Let $h_{j,P}$ be the restriction of $h_j$ to $W_P\cap U$.  Let
${\mathcal A}_P$ be an Artinian subscheme on $W_P$ (supported at
the single point $P$) parametrized holomorphically by $P\in V$.
For $0\leq\lambda\leq\rho$ we introduce the Artinian subscheme
$$
{\mathcal A}^{(\lambda)}_P={\mathcal A}_P\left/\sum_{j=1}^\lambda
h_{j,P}{\mathcal A}_P\right.\leqno{(7.5.2.2)}
$$
so that ${\mathcal A}^{(0)}_P={\mathcal A}_P$.  The notion of
continuous nested sequences of Artinian subschemes without jump
will be applied to $(7.5.2.2)$ and, in particular, to the case
where ${\mathcal A}_P={\mathcal A}_{V,N,P}$ for each $N\in{\mathbb
N}$.

\medbreak When $\rho=\ell$ and ${\mathcal A}_P={\mathcal
A}_{V,N,P}$ for each $N\in{\mathbb N}$, let
$Z_{N,h_1,\cdots,h_\ell}$ be a subvariety of codimension $\geq 1$
in $V$ such that the nested sequence $(7.5.2.1)_P$ of $\rho+1$
Artinian subschemes {\it varies continuously without jump} for
$P\in V-\check Z_{N,h_1,\cdots,h_\ell}$.  Let $Z^*$ be the union
of $Z_N$ and $Z_{N,h_1,\cdots,h_\ell}$ for all $N\geq N_0$.

\bigbreak\noindent(7.6) {\it Family of Subvarieties of One
Dimension Higher Containing Embedded Component of Stable
Base-Set.} After using the local situation to discuss sequences of
graded coherent ideal sheaves, their order functions, and the
continuous variation of a nested sequence of Artinian subschemes
without jump, we now come back to the global situation.  Assume
that vanishing orders of $\Phi$ have been achieved outside a
subvariety of codimension $\geq\ell$ in $X$.  That is, there
exists some subvariety $\hat V$ of codimension $\geq\ell$ such
that for every point $P\in X-\hat V$ there exist some open
neighborhood $U_P$ of $P$ in $X-\hat V$ and some positive number
$C_P$ and some $\hat m_\ell\in{\mathbb N}$ such that
$$
\frac{1}{C_P}\Phi_{\hat m_\ell}\leq\Phi\leq C_P\Phi_{\hat
m_\ell}\quad{\rm on\ \ }U_P,
$$
where $\Phi_{\hat m_\ell}$ is the $\hat m_\ell$-th partial sum of
$\Phi$ as explained in (3.6).  Denote $\left((n+2)\hat
m_\ell\right)!$ by $m_\ell$.  (The factor $n+2$ comes from the
theorem of Skoda on ideal generation (1.1) and the factorial is to
get a uniform grading for the pluricanonical sections used in ideal
generation.)

\medbreak Let ${\mathcal J}_\ell$ be the ideal sheaf on $X$
generated by $\Gamma\left(X, m_\ell K_X\right)$.  We blow up $X$ by
using monoidal transformations with nonsingular centers inside the
zero-set of ${\mathcal J}_\ell$ to get $\pi:\tilde X\to X$ so that
\begin{itemize}\item[(i)] $\pi^{-1}\left({\mathcal
J}_\ell\right)=\prod_{j=1}^N\left({\mathcal I}_{Y_j}\right)^{b_j}$
and \item[(ii)] $K_{\tilde X}=\pi^*K_X+\sum_{j=1}^N b_j^\prime
Y_j$,\end{itemize} where $b_j, b_j^\prime$ are nonnegative integers
and $\left\{Y_j\right\}_{1\leq j\leq N}$ is a collection of
nonsingular hypersurfaces in $\tilde X$ in normal crossing and
${\mathcal I}_{Y_j}$ is the ideal sheaf of $Y_j$.  Let $\tilde V$ be
the common zero-set of $$\Gamma\left(\tilde X, k\left(m_\ell
K_{\tilde X}-\sum_{j=1}^N\left(b_j+m_\ell
b_j^\prime\right)Y_j\right)\right)$$ for all $k\in{\mathbb N}$.
Without loss of generality we can assume that $\hat
V=\pi\left(\tilde V\right)$, otherwise we just simply replace $\hat
V$ by $\pi\left(\tilde V\right)$ and replace $\ell$ by the
codimension of $\pi\left(\tilde V\right)$ in $X$.  For some
$k_\ell\in{\mathbb N}$ the common zero-set of
$$\Gamma\left(\tilde X, k_\ell\left(m_\ell
K_{\tilde X}-\sum_{j=1}^N\left(b_j+m_\ell
b_j^\prime\right)Y_j\right)\right)$$ is $\tilde V$.  We take $\ell$
generic elements $\tilde\sigma_1,\cdots,\tilde\sigma_\ell$ of
$$\Gamma\left(\tilde X, k_\ell\left(m_\ell
K_{\tilde X}-\sum_{j=1}^N\left(b_j+m_\ell
b_j^\prime\right)Y_j\right)\right)$$ and let
$\sigma_1,\cdots,\sigma_\ell$ be the elements of $\Gamma\left(X,
k_\ell m_\ell K_X\right)$ corresponding to
$\tilde\sigma_1,\cdots,\tilde\sigma_\ell$.  Note that by blowing up
$\tilde X$ further we can also assume without loss of generality
that $\tilde V$ is a hypersurface and that all the branches of
$\tilde V$ and $Y_1,\cdots,Y_N$ together are in normal crossing.

\medbreak Let $V$ be a branch of $\hat V$ of codimension $\ell$ in
$X$.  In the notations of (7.5.2) we let $\rho=\ell$ and
$h_j=\sigma_j$ for $1\leq j\leq\ell$ and we get the countable union
$Z^*$ of subvarieties of codimension $\geq 1$ in $V$.  Let
$c^{(k)}=\left(c^{(k)}_1,\cdots,c^{(k)}_\ell\right)$ ($1\leq
k\leq\ell-1$) be ${\mathbb C}$-independent $\ell$-tuples of complex
numbers and we denote $\left(c^{(1)},\cdots,c^{(\ell-1)}\right)$ by
${\mathfrak c}$. Let $\tilde S_{\mathfrak c}\subset\tilde X$ be the
common zero-set of $\sum_{j=1}^\ell c^{(k)}_j\sigma_j$ in $X$ for
$1\leq k\leq\ell-1$.  Let $\check S_{\mathfrak c}$ be the
topological closure of $\tilde S_{\mathfrak c}-\tilde V$ in $\tilde
X$.  Let $S_{\mathfrak c}\subset X$ be the $\pi$-image of $\check
S_{\mathfrak c}$.  We can consider the element of the Grassmannian
of all $(\ell-2)$-dimensional linear subspaces ${\mathbb
P}_{\ell-2}$ in ${\mathbb P}_{\ell-1}$ defined by ${\mathfrak c}$
and naturally regard ${\mathfrak c}$ as an element of ${\mathbb
P}_{\ell-1}$ so that the family of $S_{\mathfrak c}$ is parametrized
by ${\mathfrak c}\in{\mathbb P}_{\ell-1}$.  Each $S_{\mathfrak c}$
is a subvariety of pure codimension $\ell-1$ in $X$ for ${\mathfrak
c}\in{\mathbb P}_{\ell-1}$, possibly after excluding a proper
subvariety of ${\mathbb P}_{\ell-1}$, in which case we use a finite
number of $\ell$-tuples of elements of
$$\Gamma\left(\tilde X, k_\ell\left(m_\ell K_{\tilde
X}-\sum_{j=1}^N\left(b_j+m_\ell
b_j^\prime\right)Y_j\right)\right)$$ instead of a single
$\ell$-tuple
$\left(\tilde\sigma_1,\cdots,\tilde\sigma_\ell\right)$ so that the
intersection of the excluded proper subvarieties of ${\mathbb
P}_{\ell-1}$ is empty.  In order not to distracted from the
important points of the argument, let us suppress the mention of
such a proper subvariety of ${\mathbb P}_{\ell-1}$.

\medbreak Let
$$\Theta_{K_X,\ell}:=\Theta_{\left(-\log\Phi_{K_X,\ell}\right)}
=\frac{\sqrt{-1}}{2\pi}\partial\bar\partial\log\Phi_{K_X,\ell}$$
be the curvature current from the metric
$\frac{1}{\Phi_{K_X,\ell}}$ defined by the local order function
$\Phi_{K_X,\ell}$ of the conductor $\left\{{\mathcal
J}_{\left(K_X,\,{\rm
codim}\geq\ell-1\right)}^{(\nu)}\right\}_{\nu\in{\mathbb N}}$
introduced in (7.3).  Note that the order function
$\Phi_{K_X,\ell}$ is only locally defined, for example, on some
open neighborhood $U_P$ in $X$ of some point $P\in V$.  The
curvature current $\Theta_{K_X,\ell}$ is only defined on $U_P$.
Let $\gamma_{\mathfrak c}$ be the nonnegative number such that
$$\left(\left(\Theta_{K_X,\ell}|_{S_{\mathfrak c}}\right)-\gamma_{\mathfrak
c}\left[V\right]\right)_V\leqno{(7.6.0.1)}$$ is a closed positive
$(1,1)$-current on $V\cap U_P$.  This we can obtain by using the
decomposition of $\Theta_{K_X,\ell}|_{S_{\mathfrak c}}$ according
to (3.3), (3.4) and (3.5).  Though the order function
$\Phi_{K_X,\ell}$ and $\Theta_{K_X,\ell}$ are only defined on
$U_P$, yet the nonnegative number $\gamma_{\mathfrak c}$ is
globally defined and the irreducible Lelong sets of the closed
positive $(1,1)$-current $(7.6.0.1)$ $V\cap U_P$ are globally
defined for $V$.  We have the following.

\bigbreak\noindent(7.6.1) {\it Lemma.} The irreducible Lelong sets
of the closed positive $(1,1)$-current $(7.6.0.1)$ on $V\cap U_P$
are contained in the countable union $Z^*\cap U_P$ of subvarieties
of codimension $\geq 1$ in $U_P$.

\bigbreak\noindent(7.7) {\it Construction of Pluricanonical
Sections to Achieve Stable Vanishing Order at Codimension $\ell$.}
For each ${\mathfrak c}\in{\mathbb P}_{\ell-1}$ we are going to
construct a pluricanonical section of $X$ which will achieve the
stable vanishing order of $\Phi$ at a generic point of $V$ in the
direction of $S_{\mathfrak c}$.   Here we can specify precisely
what genericity is required for this conclusion.  The genericity
condition is that the point in question has to be in $V-Z^*$ and
the important point is that the genericity condition is
independent of ${\mathfrak c}$.  (Again we may have to confine
${\mathfrak c}$ to ${\mathbb P}_{\ell-1}$ minus a proper
subvariety of ${\mathbb P}_{\ell-1}$ first and then use different
proper subvarieties with empty intersection.)

\medbreak We will apply the general non-vanishing theorem (6.2).
The idea is to use, for an appropriate nonnegative number
$\check\gamma_{\mathfrak c}$ and some sufficiently large $m$, a
line bundle $m K_X-\left\lfloor m\check\gamma_{\mathfrak
c}\right\rfloor V$ on $S_{\mathfrak c}$ and get a holomorphic
section of it over $V$ with coefficients belonging to the
appropriate multiplier ideal sheaf and then extend it first to
$S_{\mathfrak c}$ and then to an element of $\Gamma\left(X,
mK_X\right)$.  The number $\check\gamma_{\mathfrak c}$ may be
greater than $\gamma_{\mathfrak c}$ because of the contribution
from the ideal sheaf ${\mathcal J}_\ell$ on $X$ generated by
$\Gamma\left(X, m_\ell K_X\right)$ introduced in (7.6) from the
induction assumption that the stable vanishing orders of $\Phi$
can be precisely achieved outside of a subvariety of codimension
$\geq\ell$.  Note that in order to go to a global closed positive
$(1,1)$-current on $V$ from the local closed positive
$(1,1)$-current $(7.6.0.1)$ only defined on $V\cap U_P$, we must
add the contribution from the ideal sheaf ${\mathcal J}_\ell$ on
$X$ generated by $\Gamma\left(X, m_\ell K_X\right)$.

\medbreak We have to worry about branches of $S_{\mathfrak c}$ which
are in the zero-set of ${\mathcal J}_\ell$.   The easiest way to
handle this is to go up to $\tilde X$ and do the extensions in
$\tilde X$ instead of in $X$, because the pullback
$\pi^{-1}\left({\mathcal J}_\ell\right)$ to $\tilde X$ of the ideal
sheaf ${\mathcal J}_\ell$ on $X$ is of the form
$\prod_{j=1}^N\left({\mathcal I}_{Y_j}\right)^{b_j}$ and can be
removed from $\tilde X$ as a line bundle.  Let $m=km_\ell$. Earlier,
before we worry about branches of $S_{\mathfrak c}$ which are in the
zero-set of ${\mathcal J}_\ell$,
\begin{itemize}\item[(i)] we apply the general non-vanishing theorem (6.2) to get a
holomorphic section of $m K_X-\left\lfloor m\check\gamma_{\mathfrak
c}\right\rfloor V$ over $V$ with coefficients in the appropriate
multiplier ideal sheaf and \item[(ii)] then extend it first to
$S_{\mathfrak c}$ and \item[(iii)] finally to an element of
$\Gamma\left(X, mK_X\right)$.\end{itemize} After we worry about
branches of $S_{\mathfrak c}$ which are in the zero-set of
${\mathcal J}_\ell$, we carry out the equivalent extensions in
$\tilde X$ in the following way.
\begin{itemize}\item[(i)] We apply the
general non-vanishing theorem (6.2) to get a holomorphic section of
$$k\left(m_\ell K_{\tilde X}-\sum_{j=1}^N\left(b_j+m_\ell
b_j^\prime\right)Y_j\right)-\left\lfloor km_\ell\gamma_{\mathfrak
c}\right\rfloor V_{\mathfrak c}$$ over $V_{\mathfrak c}$ with
coefficients in the appropriate multiplier ideal sheaf and
\item[(ii)] then extend it first to $\check S_{\mathfrak c}$ and
\item[(iii)] then to an element of
$$\Gamma\left(X, k\left(m_\ell K_{\tilde
X}-\sum_{j=1}^N\left(b_j+m_\ell b_j^\prime\right)Y_j\right)\right)$$
\item[(iv)] and finally multiply it by
$\prod_{j=1}^N\left(s_{Y_j}\right)^{b_j}$ to get an element of
$\Gamma\left(X, k\left(m_\ell K_{\tilde X}\right)\right)$, where
$s_{Y_j}$ is the canonical section of $Y_j$.
\end{itemize}

\medbreak Thus for each ${\mathfrak c}\in{\mathbb P}_{\ell-1}$ we
succeed in constructing a pluricanonical section of $X$ which will
achieve the stable vanishing order of $\Phi$ at a generic point of
$V$ in the direction of $S_{\mathfrak c}$ and, moreover, the
generic point of $V$ can be chosen to be independent of
${\mathfrak c}$.  This finishes the descending induction on the
dimension of the subvariety where the vanishing order of $\Phi$ is
not achieved. It finishes also the proof of the finite generation
of the canonical ring for the case of general type.

\medbreak Since in order to rigorously apply the general
non-vanishing theorem (6.2) we have to go to the blown-up manifold
$\tilde X$ where the inverse image $\tilde V$ of $V$ under the
blow-up map is a nonsingular hypersurface, the question naturally
arises whether we are simply applying the hypersurface case
argument to $\tilde V$ in $\tilde X$ to handle the problem of
achieving the precise vanishing order of $\Phi$ at a generic point
of $V$. What is precisely the r\^ole played by the Artinian
subschemes? Is it necessary to introduce the Artinian subschemes
and to use the theory of the continuous variation of an Artinian
subscheme without jump? We are going to answer these questions in
the next paragraph and highlight the r\^ole played by the Artinian
subschemes and their continuous variation without jump.

\bigbreak\noindent(7.8) {\it Difference Between Using Artinian
Subschemes and Simply Blowing Up to Reduce to the Hypersurface
Case.}  If we just apply the hypersurface argument to $\tilde V$
in $\tilde X$, we would consider the closed positive
$(1,1)$-current
$$
\left(\Theta_{\log\tilde\Phi}-\tilde\gamma_{\tilde V}\right)_{\tilde
V}
$$
where the metric $\frac{1}{\tilde\Phi}$ of $K_{\tilde X}$ is to
$\tilde X$ as the metric $\frac{1}{\Phi}$ of $K_X$ is to $X$.  The
nonnegative number $\tilde\gamma_{\tilde V}$ is the Lelong number of
the closed positive $(1,1)$-current
$\Theta_{\left(-\log\tilde\Phi\right)}$ on $\tilde X$ at a generic
point of $\tilde V$.  We can get the precise achievement of the
vanishing order of $\tilde\Phi$ at points of $\tilde V-\tilde Z$ for
some subvariety $\tilde Z$ of codimension $\geq 1$ in $\tilde V$.
However, in general the $\pi$-image $\pi\left(\tilde Z\right)$ in
$X$ contains $V$ and we cannot draw any conclusion about the precise
achievement of vanishing order of $\Phi$ at a generic point of $V$.

\medbreak When we use $\check S_{\mathfrak c}$ and consider the line
bundle $$k\left(m_\ell K_{\tilde X}-\sum_{j=1}^N\left(b_j+m_\ell
b_j^\prime\right)Y_j\right)-\left\lfloor km_\ell\gamma_{\mathfrak
c}\right\rfloor V_{\mathfrak c}$$ over $V_{\mathfrak c}$ and the
closed positive $(1,1)$-current $$\left(k\left(m_\ell
\Theta_{\left(-\log\tilde\Phi\right)}-\sum_{j=1}^N\left(b_j+m_\ell
b_j^\prime\right)\left[Y_j\right]\right)-\left\lfloor
km_\ell\gamma_{\mathfrak c}\right\rfloor\left[V_{\mathfrak
c}\right]\right)\bigg|_{V_{\mathfrak c}}$$ on $V_{\mathfrak c}$, the
nonnegative number $\gamma_{\mathfrak c}$ depends on ${\mathfrak c}$
and is in general different from $\tilde\gamma_{\tilde V}$ which is
independent of ${\mathfrak c}$, even after we take into account the
contributions from $\sum_{j=1}^N\left(b_j+m_\ell
b_j^\prime\right)Y_j$.  Let us use the picture in $X$ instead of in
$\tilde X$ and discount the contribution from the stable vanishing
order at a generic point of a subvariety of codimension
$\leq\ell-1$.  In this picture what is going on is that
$\gamma_{\mathfrak c}$ is from the computation of the vanishing
order of $\Phi$ at a generic point of $V$ along a generic curve
inside $S_{\mathfrak c}$, whereas $\tilde\gamma_{\tilde V}$ comes
from the computation of the vanishing order of $\Phi$ at a generic
point of $V$ along any generic curve in $X$ (without the additional
condition that the curve be on $S_{\mathfrak c}$).

\medbreak The use of Artinian subschemes and their continuous
variation without jump allows us to take away different generic
vanishing orders $\gamma_{\mathfrak c}$ along $V_{\mathfrak c}$
depending on ${\mathfrak c}$.  Each $V_{\mathfrak c}$ inside $\tilde
V$ is mapped onto $V$ by $\pi$, but the subvariety of codimension
$\geq 1$ in $V_{\mathfrak c}$ where the vanishing order along
$V_{\mathfrak c}$ is higher than at a generic point of $V_{\mathfrak
c}$ is mapped under $\pi$ to the countable union $Z^*$ of
subvarieties of codimension $\geq 1$ in $V$ which is independent of
${\mathfrak c}$.  The picture is that in the directions normal to
$V$ the vanishing order of $\Phi$ depends very much on the
direction.  The use of Artinian subschemes and their continuous
variation without jump allows us to identify precisely the generic
vanishing order at a point of $V$ along a generic curve inside
$S_{\mathfrak c}$ in order to achieve that particular vanishing
order and, moreover, to do it at points of $V-Z^*$ with the
countable union $Z^*$ of subvarieties of codimension $\geq 1$ in $V$
independent of ${\mathfrak c}$.

\bigbreak

\bigbreak\noindent{\bf \S8. Extension Techniques from the Proof of
the Invariance of Plurigenera and the Finite Generation of the
Canonical Ring.}  The development of extension techniques for the
problem of the deformational invariance of plurigenera was intended
for application to the problem of the finite generation of the
canonical ring.

\bigbreak\noindent(8.1) {\it Heuristic Discussion.}  The most
obvious way of applying the extension techniques from the proof of
the deformational invariance of plurigenera to the finite generation
of the canonical ring is to try to implement in a rigorous and
precise manner the following strategy.

\bigbreak\noindent(i) Choose a divisor $Y$ of some pluricanonical
section $s_0\in\Gamma\left(X,m_0K_X\right)$ with
$K_Y=\left(m_0+1\right)K_X|_Y$ and show that
$$
\rho: \Gamma\left(X,
m\left(m_0+1\right)K_X\right)\to\Gamma\left(Y,m\left(m_0+1\right)K_X\right)=\Gamma\left(Y,mK_Y\right)
$$
is surjective.

\bigbreak\noindent(ii) Since the dimension of $Y$ is one lower
than that of $X$, use the induction assumption to conclude that
the ring $\bigoplus_{m=1}^\infty\Gamma\left(Y,mK_Y\right)$ is
finitely generated so that for some $p_1\in{\mathbb N}$ the ring
$\bigoplus_{m=1}^\infty\Gamma\left(Y,mK_Y\right)$ is generated by
its finitely truncated part
$\bigoplus_{m=1}^{p_1}\Gamma\left(Y,mK_Y\right)$.

\bigbreak\noindent(iii)  Lift by $\rho$ the generators
$$
\sigma^{(m_\nu)}_\nu\in\Gamma\left(Y, m_\nu
K_Y\right)\quad{1\leq\nu\leq N}
$$
(from $\bigoplus_{m=1}^{p_1}\Gamma\left(Y,mK_Y\right)$) to the
following holomorphic pluricanonical sections on $X$.
$$
\tilde\sigma^{(m_\nu)}_\nu\in\Gamma\left(X,
m_\nu\left(m_0+1\right) K_X\right)\quad{1\leq\nu\leq N}.
$$

\bigbreak\noindent(iv) Take $q\in{\mathbb N}$ and any
$s\in\Gamma\left(X,q\left(m_0+1\right)K_X\right)$.  In terms of
the finite number of generators of
$\bigoplus_{m=1}^\infty\Gamma\left(Y,mK_Y\right)$ we can write
$$
\rho(s)=\sum_{j_1m_1+\cdots j_N
m_N=q}c_{j_1,\cdots,j_N}\prod_{\lambda=1}^N\left(\sigma^{\left(m_\nu\right)}_\nu\right)^{j_\nu}
$$
on $Y$ for some $c_{j_1,\cdots,j_N}\in{\mathbb C}$.

\bigbreak\noindent(v)  We lift the right-hand side of the above
equation via $\rho$ to $X$ and divide by $s_0$ to get
$$
s-\sum_{j_1m_1+\cdots j_N
m_N=q}c_{j_1,\cdots,j_N}\prod_{\lambda=1}^N\left(\tilde\sigma^{\left(m_\nu\right)}_\nu\right)^{j_\nu}
=s_0 s^\prime
$$
for some
$s^\prime\in\Gamma\left(X,\left(\left(q-1\right)\left(m_0+1\right)+1\right)K_X\right)$.

\bigbreak\noindent(vi)  The natural next step is to apply the
preceding argument to $s^\prime$ instead of $s$. However, there is
the difficulty of $s^\prime$ being in
$\Gamma\left(X,\left(\left(q-1\right)\left(m_0+1\right)+1\right)K_X\right)$
instead of in
$\Gamma\left(X,\left(q-1\right)\left(m_0+1\right)K_X\right)$.  To
overcome this difficulty, we use the theorem (3.8) on twisted
finite generation (after stable vanishing orders are known to be
achieved) with $E=K_X$ endowed with the metric
$\frac{1}{\Phi}=e^{-\varphi}$ to conclude that if the stable
vanishing orders for $Y$ is achieved by $\tilde
m_Y$-pluricanonical sections on $Y$, then $\Gamma\left(X,
\left(m\left(m_0+1\right)m_Y+1\right)K_X\right)|_Y$ is contained
in
$$
\left(\Gamma\left(Y,m_Y K_Y\right)\right)^{p_m}\Gamma\left(Y,
\left(\left(m-p_m\right)K_Y+K_X\right)\right)$$ for
$m\geq\left(n+1\right)(m_0+1)m_Y$, where $m_Y=\tilde m_Y !$ and
$p_m=\left\lfloor\frac{m}{\left(m_0+1\right)m_Y}\right\rfloor-(n+1)$.

\medbreak We now needs to extend pluricanonical sections on $Y$ with
twisting by $pK_X$ together with supremum bound condition with
respect to $\frac{1}{\Phi}$ for $p=1$ and later also for $1\leq
p\leq m_0$.  For the deformational invariance of plurigenera this
twisting extension result was proved in [Siu 2002].  (Recently Paun
[Paun 2005] generalized it to the case of $L^2$ bound.)  Thus we can
apply the modified argument to $s^\prime$.  After a finite number of
application of this modified argument we will get the finite
generation of the canonical ring of $X$.

\medbreak The problem arises concerning the problem of extending
pluricanonical sections of a hypersurface to pluricanonical sections
of the ambient space.  We are going to discuss the approach to this
problem.

\bigbreak\noindent(8.2) {\it Extension of Pluricanonical Sections
From a Hypersurface.}  Let $Y$ be a regular hypersurface in $X$ so
that the ideal sheaf ${\mathcal I}_Y$ on $X$ is equal to the
multiplier ideal sheaf ${\mathcal I}_{\varphi_Y}$ of some metric
$e^{-\varphi_Y}$ whose curvature current dominates some smooth
positive $(1,1)$-form on $X$.  Let $L$ be a line bundle over $X$
with a metric $e^{-\varphi_L}$ whose curvature current is
nonnegative.  Then the vanishing theorem of Kawamata-Viehweg-Nadel
gives
$$
H^1\left(X,{\mathcal
I}_{\varphi_L+\varphi_Y}\left(L+Y+K_X\right)\right)=0.\leqno{(8.2.1)}
$$
Take the short exact sequence
$$
0\to{\mathcal I}_{\varphi_L+\varphi_Y}\hookrightarrow{\mathcal
I}_{\varphi_L}\to {\mathcal I}_{\varphi_L}\left/{\mathcal
I}_{\varphi_L+\varphi_Y}\right.\to 0
$$
and tenor it with $L+Y+K_X$ to get the short exact sequence
$$
\displaylines{\qquad0\to{\mathcal
I}_{\varphi_L+\varphi_Y}\left(L+Y+K_X\right)\to {\mathcal
I}_{\varphi_L}\left(L+Y+K_X\right)\hfill\cr\hfill\to \left({\mathcal
I}_{\varphi_L}\left/{\mathcal
I}_{\varphi_L+\varphi_Y}\right.\right)\left(L+Y+K_X\right)\to
0.\qquad\cr}
$$
From its long cohomology exact sequence and the vanishing of first
cohomology group in $(8.2.1)$ it follows that
$$
\Gamma\left(X,{\mathcal I}_{\varphi_L}\left(L+Y+K_X\right)\right)\to
\Gamma\left(X,\left({\mathcal I}_{\varphi_L}\left/{\mathcal
I}_{\varphi_L+\varphi_Y}\right.\right)\left(L+Y+K_X\right)\right)
$$
is surjective.  This is an analog of the extension theorem of
Ohsawa-Takegoshi type.  The reason for dubbing it as an analog is
that we can write $Y+K_X$ as $K_Y$ so that we have the surjectivity
of
$$
\Gamma\left(X,{\mathcal I}_{\varphi_L}\left(L+K_Y\right)\right)\to
\Gamma\left(X,\left({\mathcal I}_{\varphi_L}\left/{\mathcal
I}_{\varphi_L+\varphi_Y}\right.\right)\left(L+K_Y\right)\right)
$$
and we can interpret
$$
\Gamma\left(X,\left({\mathcal I}_{\varphi_L}\left/{\mathcal
I}_{\varphi_L+\varphi_Y}\right.\right)\left(L+K_Y\right)\right)
$$
as the analog of all top-degree forms on $Y$ with coefficient in $L$
which is $L^2$ with respect to the metric $e^{-\varphi_L}$ of $L$.

\medbreak At this point the natural strategy is to apply the
``two-tower'' extension technique from the proof of the invariance
of plurigenera (at least for the case of general type) to our
situation at hand to get the surjectivity of
$$
\Gamma\left(X,{\mathcal I}_{\varphi_L}\left(L+mK_Y\right)\right)\to
\Gamma\left(X,\left({\mathcal I}_{\varphi_L}\left/{\mathcal
I}_{\varphi_L+\varphi_Y}\right.\right)\left(L+mK_Y\right)\right)
$$
for $m\geq 1$ (see [Siu 1998], [Siu 2002], [Siu 2003], [Siu 2005],
[Paun 2005]). The ``two-tower'' argument involves

\begin{itemize}\item[(i)] raising the given section to be extended first to a high power $N$,
\item[(ii)] adding the twisting by a sufficiently ample line
bundle $A$ for global generation of multiplier ideal sheaves,
\item[(iii)] adding the canonical line bundle one copy at one type
from the extension theorem of Ohsawa-Takegoshi type, , and
\item[(iv)] finally using the $N$-th root of the absolute value of
the extension of the twisted sections (or its limit for the case
of non general type) to apply the extension theorem of
Ohsawa-Takegoshi type to construct the extension of the given
section.
\end{itemize}
The main difficulty in implementing this precisely and rigorously
is the following. The coherent sheaf ${\mathcal
I}_{\varphi_L}\left/{\mathcal I}_{\varphi_L+\varphi_Y}\right.$
supported on $Y$ is in general over some unreduced structure sheaf
of $Y$ ({\it i.e.,} a structure sheaf with nonzero nilpotent
elements) and is not the multiplier ideal sheaf for some
$e^{-\kappa}$ with $\kappa$ locally plurisubharmonic on $Y$.  One
encounters formidable obstacles when one tries to implement the
step of taking roots of absolute value of a section over an
unreduced structure.  Earlier we have seen this kind of problem
with unreduced structures in (4.4) and (6.10), where the problem
is handled by the analog of minimal centers of log canonical
singularities, but this kind of handling cannot be used here.

\medbreak To implement precisely and rigorously this approach to the
finite generation by extension techniques from the proof of the
deformational invariance of plurigenera would involve a tremendous
amount of tedious geometrically-uninspiring messy bookkeeping and
the task of verifying the correctness of such an implementation
after it is done would also be daunting.

\medbreak For the finite generation of the canonical ring we put
aside the approach of extension techniques from the plurigenera
problem and use instead the approach presented in these notes for
the following reason.  The approach presented in these notes is
geometrically more enlightening.  It gives us a clear geometric
picture of how and why a general non-vanishing theorem can be proved
from techniques for Fujita conjecture type problems with input from
the diophantine-approximation contribution of infinite number of
irreducible Lelong sets and Shokurov's technique of comparing the
theorem of Riemann-Roch for a general line bundle and its twisting
by a flat line bundle.  It also provides a clear geometric picture
of how and why we can use Artinian subschemes to give a useful
analog of Lelong numbers for the hypersurface and locate an {\it a
priori} bad set in an embedded stable base-point set of higher
codimension by using the continuous deformation of Artinian
subschemes without jump so that such {\it a priori} bad sets can be
used to finish the induction process in a finite number of steps.

\bigbreak

\bigbreak\noindent{\bf \S9. Remark on Positive Lower Bound of
Curvature Current}

\bigbreak The proof of the finite generation of the canonical ring
for the case of general type hinges on the small positivity of the
curvature currents for certain line bundles.  This small
positivity is needed both for the general non-vanishing theorem
and for the vanishing theorem of Kawamata-Viehweg-Nadel.  The
question is whether it is possible to prove the finite generation
of the canonical ring without the assumption of general type by
first artificially introducing some positive line bundle and then
later using some limit process to get rid of the
artificially-introduced positivity of the curvature current.

\medbreak The proof of the deformational invariance of the
plurigenera for the general case without the assumption of general
type follows this particular strategy of artificially introducing
some positivity and then getting rid of it by taking limit.  The
challenge is to handle well the limiting process.  This challenge
manifests itself already in the proof of the deformational
invariance of the plurigenera for the general case without the
assumption of general type where the convergence of the metric in
taking the limit is the key difficulty which has to be handled.
For the generation of the canonical ring without the general type
assumption the situation of such a limiting process to get a
non-vanishing theorem to precisely achieve stable vanishing orders
would be far more involved than the situation for the
deformational invariance of plurigenera.

\medbreak In this section we will make some remarks concerning the
positive lower bound of the curvature current in two situations
related to the proof of the finite generation of the canonical ring.

\bigbreak\noindent(9.1) {\it General Non-Vanishing Theorem and
Kawamata-Viehweg-Nadel Vanishing as a Pair.}  There are two
fundamental theorems in the Oka-Cartan theory of several complex
variables, called Theorem A and Theorem B.  Theorem B states that
$H^\nu\left(S,{\mathcal F}\right)=0$ for $\nu\geq 1$ when $S$ is a
Stein manifold (or space in general) and ${\mathcal F}$ is a
coherent analytic sheaf on $S$.  Theorem A states that any
coherent analytic sheaf ${\mathcal F}$ on a Stein space $S$ is
generated at every point $P$ over the local ring ${\mathcal
O}_{S,P}$ by elements of $\Gamma\left(S,{\mathcal F}\right)$.
Theorem A and Theorem B come together as a pair of fundamental
results in the theory of Stein spaces.

\medbreak In the context of compact complex algebraic manifolds and
multiplier ideal sheaves there are two results in a pair which are
analogous to Theorem A and Theorem B in the theory of Stein spaces.
The first result in the pair is the vanishing theorem of
Kawamata-Viehweg-Nadel.  It applies to a holomorphic line bundle $L$
over a compact complex algebraic manifold $Y$ with metric
$e^{-\varphi}$ whose curvature current admits a strict positive
lower bound (in the sense that it dominates some positive smooth
$(1,1)$-form).  The other result in the pair is the theorem on the
global generation of multiplier ideal sheaves.  It was first
introduced to prove the deformational invariance of plurigenera [Siu
1998]. Its statement is as follows.

\bigbreak\noindent(9.1.1) {\it Theorem on Global Generation of
Multiplier Ideal Sheaves.} Let $Y$ be a compact complex projective
algebraic manifold of complex dimension $m$, and $L$ be a
holomorphic line bundle on $Y$ with metric $e^{-\varphi}$ where
$\varphi$ is plurisubharmonic. Let $E$ be a holomorphic line
bundle over $Y$ sufficiently ample in the sense that for every
$P\in Y$ there exist a finite number of elements of $\Gamma(Y,E)$
which vanish to order $\geq m+1$ at $P$ and do not vanish
simultaneously outside $P$. Then $\Gamma(Y,{\mathcal
I}_{\varphi}\otimes (L+E+K_Y))$ generates ${\mathcal
I}_\varphi\otimes (L+E+K_Y)$.

\medbreak For the global generation of multiplier ideal sheaves, a
sufficiently ample line bundle $E$ is required.  If we get rid of
$E$ by replacing $L$ by $L-E$, the global generation of multiplier
ideal sheaves can also be formulated with a trivial line bundle
$E$ but the new formulation would instead require the curvature
current $e^{-\varphi}$ of $L$ to dominate the curvature form of
some sufficiently ample line bundle (see the proof of (9.2.1)).

\medbreak Though the vanishing theorem of Kawamata-Viehweg-Nadel
and the global generation of multiplier ideal sheaves form a pair,
yet their assumptions on the positive lower bound for the
curvature current are not comparable.  For the vanishing theorem
of Kawamata-Viehweg-Nadel any small positive lower bound suffices.
However, for the global generation of multiplier ideal sheaves the
lower bound for the curvature current has to be sufficiently
positive.

\medbreak For the proof of the deformational invariance of
plurigenera, the difficulty from the undesirable assumption of
sufficient positivity of the curvature current for the global
generation of multiplier ideal sheaves is alleviated by taking roots
of sufficiently high power.  This technique of alleviating the
difficulty by taking roots suffices for the proof of the
deformational invariance of plurigenera, because the problem of the
deformational invariance of plurigenera is less delicate than the
problem of the finite generation of the canonical ring.

\medbreak The general non-vanishing theorem (6.2) which we introduce
for the purpose of proving the finite generation of the canonical
ring for the case of general type fits better with the theorem of
Kawamata-Viehweg-Nadel as a pair.  Like the theorem of
Kawamata-Viehweg-Nadel the general non-vanishing theorem (6.2) only
requires some positivity, no matter how small, for the curvature
current of $e^{-\chi}$.  Of course, in the proof of the general
non-vanishing theorem (6.2) we have to take some high-order root of
some holomorphic section vanishing to high order.  However, such a
root-taking step is done in the proof instead of in the application
of the theorem.

\bigbreak\noindent(9.2) {\it No Strict Positive Lower Bound for
Curvature Current of $\frac{1}{\Phi}$ with Hypersurface Lelong
Set.} The metric $e^{-\varphi}=\frac{1}{\Phi}$ of $K_X$ used in
the proof of the finite generation of the canonical ring is the
metric with the least singularity.  When $X$ is of general type,
it is natural to suspect that its curvature current
$\Theta_\varphi$ dominates some positive smooth $(1,1)$-form on
$X$.  However, this is not necessarily the case.  Here is a simple
statement which explains why the general type condition does not
in general imply strict positivity of the curvature current
$\Theta_\varphi$.

\bigbreak\noindent(9.2.1) {\it Proposition.}  If the curvature
current $\Theta_\varphi$ has a hypersurface Lelong set, then it
cannot dominate any positive smooth $(1,1)$-form $\omega_0$ on
$X$.

\medbreak\noindent{\it Proof.} Suppose the contrary and there is
such a positive smooth $(1,1)$-form $\omega_0$ on $X$ such that
$\Theta_\varphi\geq\omega_0$ as $(1,1)$-currents. By assumption we
have the following decomposition
$$
\Theta_\varphi=\sum_{j=1}^J\tau_j\left[V_j\right]+R,\leqno{(9.2.1.1)}
$$
where $J\in{\mathbb N}\cup\left\{\infty\right\}$, each $\tau_j$ is a
positive number,  each $V_j$ is an irreducible hypersurface in $X$,
and the Lelong number of $R$ is zero outside a countable union of
subvarieties of codimension at least two in $X$.

\medbreak Let $A$ be a holomorphic line bundle on $X$ which is
sufficiently ample in the sense that for every $P\in X$ there
exist a finite number of elements of $\Gamma(X,A)$ which vanish to
order $\geq n+1$ at $P$ and do not vanish simultaneously outside
$P$.  Let $h_A$ be a smooth metric of $A$ whose curvature form
$\Theta_A$ is positive. Choose $m_0\in{\mathbb N}$ such that
$m_0\omega_0-\Theta_A$ is positive on $X$.  For $m\geq m_0$ the
curvature current of the metric $\frac{e^{-\varphi}}{h_A}$ of the
line bundle $mK_X-A$ is positive. Since the multiplier ideal sheaf
of the metric $\frac{e^{-\varphi}}{h_A}$ is ${\mathcal
I}_{m\varphi}$, it follows from the theorem on the global
generation of the multiplier ideal sheaf that the multiplier ideal
sheaf ${\mathcal I}_{m\varphi}$ is generated by
$$\Gamma\left(X,{\mathcal I}_{m\varphi}\left(\left(m+1\right)K_X\right)\right)=\Gamma\left(X,{\mathcal I}_{m\varphi}
\left(\left(mK_X-A\right)+A+K_X\right)\right).$$ Pick a regular
point $P_0$ of $V_1$ which is not in any $V_j$ for $j>1$ such that
the Lelong number of $R$ is zero at $P_0$.  Then ${\mathcal
I}_{m\varphi}$ is equal to the ideal sheaf of $\left\lfloor
m\tau_1\right\rfloor V_1$ at $P_0$.  Since the multiplier ideal
sheaf ${\mathcal I}_{m\varphi}$ is generated by
$\Gamma\left(X,{\mathcal
I}_{m\varphi}\left(\left(m+1\right)K_X\right)\right)$, it follows
that at $P_0$ the vanishing order along $V_1$ of some element $s$
of $\Gamma\left(X,{\mathcal
I}_{m\varphi}\left(\left(m+1\right)K_X\right)\right)$ is equal to
$\left\lfloor m\tau_j\right\rfloor$ which is less than
$(m+1)\tau_1$, contradicting the decomposition $(9.2.1.1)$. Q.E.D.

\bigbreak\noindent{\it References}

\medbreak\noindent[Angehrn-Siu 1995] U. Angehrn and Y.-T. Siu,
Effective freeness and point separation for adjoint bundles. {\it
Invent. Math.} {\bf 122} (1995), 291--308.

\medbreak\noindent[Kawamata 1982] Y. Kawamata, A generalization of
Kodaira-Ramanujam's vanishing theorem. {\it Math. Ann.} {\bf 261}
(1982), 43-46.

\medbreak\noindent[Kawamata 1997] Y. Kawamata, On Fujita's freeness
conjecture for $3$-folds and $4$-folds. {\it Math. Ann.} {\bf 308}
(1997), 491-505.

\medbreak\noindent[Hardy-Wright 1960] G. H. Hardy and E. M. Wright,
{\it An Introduction to the Theory of Numbers}, 4th ed., Oxford
University Press 1960.

\medbreak\noindent[Kiselman 1979] C. Kiselman, Densit\'e des
fonctions plurisousharmoniques. {\it Bull. Soc. Math. France} {\bf
107} (1979), no. 3, 295--304.

\medbreak\noindent[Kawamata 1985] Y. Kawamata, Pluricanonical
systems on minimal algebraic varieties. {\it Invent. Math.} {\bf 79}
(1985), 567--588.

\medbreak\noindent[Nadel 1990] A. Nadel, Multiplier ideal sheaves
and K\"ahler-Einstein metrics of positive scalar curvature. {\it
Ann. of Math.} {\bf 132} (1990), 549-596.

\medbreak\noindent[Paun 2005] M. Paun, Siu's invariance of
plurigenera: a one-tower proof, preprint 2005.

\medbreak\noindent[Shokurov 1985] V.~V. Shokurov, A nonvanishing
theorem. {\it Izv. Akad. Nauk SSSR Ser. Mat.} {\bf 49} (1985),
635--651.

\medbreak\noindent[Siu 1974] Y.-T. Siu, Analyticity of sets
associated to Lelong numbers and the extension of closed positive
currents. {\it Invent. Math.} {\bf 27} (1974), 53-156.

\medbreak\noindent[Siu 1998]  Y.-T. Siu, Invariance of Plurigenera,
{\it Invent. Math.} {\bf 134} (1998), 661-673.

\medbreak\noindent[Siu 2002]  Y.-T. Siu, Extension of twisted
pluricanonical sections with plurisubharmonic weight and invariance
of semipositively twisted plurigenera for manifolds not necessarily
of general type. In: {\it Complex Geometry: Collection of Papers
Dedicated to Professor Hans Grauert}, Springer-Verlag 2002,
pp.223-277.

\medbreak\noindent[Siu 2003] Y.-T. Siu, Invariance of Plurigenera
and Torsion-Freeness of Direct Image Sheaves of Pluricanonical
Bundles, In: {\it Finite or Infinite Dimensional Complex Analysis
and Applications} (Proceedings of the 9th International Conference
on Finite or Infinite Dimensional Complex Analysis and Applications,
Hanoi, 2001), edited by Le Hung Son, W. Tutschke, C.C. Yang, Kluwer
Academic Publishers 2003, pp.45-83.

\medbreak\noindent[Siu 2005]  Y.-T. Siu, Multiplier ideal sheaves in
complex and algebraic geometry, {\it Science in China,
Ser.A:Math.}{\bf 48} (2005), 1-31 (arXiv:math.AG/0504259).

\medbreak\noindent[Skoda 1972] H. Skoda, Application des techniques
$L^2$ \`a la th\'eorie des id\'eaux d'une alg\`ebre de fonctions
holomorphes avec poids. {\it Ann. Sci. \'Ecole Norm. Sup.} {\bf 5}
(1972), 545-579.

\medbreak\noindent[Viehweg 1982] E. Viehweg, Vanishing theorems.
{\it J. Reine Angew. Math.} {\bf 335} (1982), 1-8.

\bigbreak\noindent{\it Author's mailing address}: Department of
Mathematics, Harvard University, Cambridge, MA 02138, U.S.A.

\medbreak\noindent {\it Author's e-mail address}:
siu@math.harvard.edu

\end{document}